        \input{epsf}   
   
        \newcommand\Z{{\mathbb Z}}   
	\newcommand\R{{\mathbb R}}   
	\newcommand\C{{\mathbb C}}   
        \renewcommand\H{{\mathbb H}}   
        \newcommand\K{{\mathbb K}} 
	\renewcommand\O{{\mathbb O}}

        \newcommand\RP{{\mathbb {RP}}}   
        \newcommand\CP{{\mathbb {CP}}}   
        \newcommand\HP{{\mathbb {HP}}}   
        \newcommand\KP{{\mathbb {KP}}}   
        \newcommand\OP{{\mathbb {OP}}}   
        \renewcommand\P{{\mathbb P}}   
   
        \newcommand{\Cliff}{{\rm Cliff}}   
        \newcommand{\J}{{\rm J}} 
        \newcommand{\M}{{\rm M}}
  
        \newcommand{\OO}{{\rm O}}   
	\newcommand{\SO}{{\rm SO}}   
	\newcommand{\SL}{{\rm SL}}   
        \newcommand{\PSL}{{\rm PSL}}   
	\newcommand{\SU}{{\rm SU}}   
        \newcommand{\Sp}{{\rm Sp}}
        \newcommand{\Spin}{{\rm Spin}}   
        \newcommand{\Pin}{{\rm Pin}}   
	\newcommand{\U}{{\rm U}}   
	\newcommand{\E}{{\rm E}}   
	\newcommand{\F}{{\rm F}}   
	\newcommand{\G}{{\rm G}}   
   
	\newcommand{\so}{{\mathfrak {so}}}   
	\newcommand{\Sl}{{\mathfrak {sl}}}   
	\newcommand{\symp}{{\mathfrak {sp}}}   
	\renewcommand{\u}{{\mathfrak {u}}}   
	\newcommand{\e}{{\mathfrak {e}}}   
	\newcommand{\f}{{\mathfrak {f}}}   
	\newcommand{\g}{{\mathfrak {g}}}   
	\newcommand{\su}{{\mathfrak {su}}}   
           
        \newcommand{\sa}{{\mathfrak {sa}}}   
        \newcommand{\h}{{\mathfrak {h}}}   
        \newcommand{\sh}{{\mathfrak {sh}}}   
           
        \newcommand{\isom}{{\mathfrak {isom}}}   
        \newcommand{\Isom}{{\rm Isom}} 
        \newcommand{\Der}{{\mathfrak {der}}}   
        \newcommand{\Tri}{{\mathfrak {tri}}}

	\newcommand{\et}{\hspace{-0.08in}{\bf .}\hspace{0.1in}}   
           
 	\newcommand{\BOX}{\hbox {$\sqcap$ \kern -1em $\sqcup$}}

        \renewcommand{\Re}{{\rm Re}}   
        \renewcommand{\Im}{{\rm Im}}   
	   
        \newcommand{\Aut}{{\rm Aut}}   
	   
	\newcommand{\tensor}{\otimes}   
        \newcommand{\iso}{\cong}   
        \newcommand{\implies}{\Longrightarrow}

	\newcommand{\be}{\begin{equation}}   
        \newcommand{\ee}{\end{equation}}   
        \newcommand{\ba}{\begin{eqnarray}}   
        \newcommand{\ea}{\end{eqnarray}}   
        \newcommand{\ban}{\begin{eqnarray*}}   
        \newcommand{\ean}{\end{eqnarray*}}   
  
        \newcommand{\ad}{{\rm ad}}

        \newcommand{\maps}{\colon}   
	\newcommand{\tr}{{\rm tr}}   
	\newcommand{\rank}{{\rm rank}}

	\newtheorem{thm}{Theorem}   
	\newtheorem{prop}{Proposition}

\documentclass {article} 
\usepackage {amsfonts}

\begin{document}   
 
	\begin{center}   
	{\bf  The Octonions \\}   
	{\em John\ C.\ Baez\\}   
	\vspace{0.3cm}   
	{\small Department of Mathematics \\   
	University of California\\   
        Riverside CA 92521\\}   
        {\small email:  baez@math.ucr.edu\\} 
	\vspace{0.3cm}   
	{May 16, 2001 \\}   
	\vspace{0.3cm}   
	\end{center}   

\begin{abstract} 
\noindent  
The octonions are the largest of the four normed division algebras.
While somewhat neglected due to their nonassociativity, they stand at
the crossroads of many interesting fields of mathematics.  Here we
describe them and their relation to Clifford algebras and spinors, Bott
periodicity, projective and Lorentzian geometry, Jordan algebras, and
the exceptional Lie groups.  We also touch upon their applications in
quantum logic, special relativity and supersymmetry. 
\end{abstract}

\section{Introduction}   
   
There are exactly four normed division algebras: the real numbers  
($\R$), complex numbers ($\C$), quaternions ($\H$), and octonions  
($\O$).  The real numbers are the dependable breadwinner of the family,  
the complete ordered field we all rely on.  The complex numbers are a  
slightly flashier but still respectable younger brother: not ordered,  
but algebraically complete.  The quaternions, being noncommutative, are  
the eccentric cousin who is shunned at important family gatherings.  But  
the octonions are the crazy old uncle nobody lets out of the attic: they  
are {\it nonassociative}.  
   
Most mathematicians have heard the story of how Hamilton invented the   
quaternions.  In 1835, at the age of 30, he had discovered how to treat   
complex numbers as pairs of real numbers.   Fascinated by the relation   
between $\C$ and 2-dimensional geometry, he tried for many years to   
invent a bigger algebra that would play a similar role in 3-dimensional   
geometry.  In modern language, it seems he was looking for a 3-dimensional   
normed division algebra.  His quest built to its climax in October 1843.   
He later wrote to his son, ``Every morning in the early part of the   
above-cited month, on my coming down to breakfast, your (then) little   
brother William Edwin, and yourself, used to ask me: `Well, Papa, can you   
{\it multiply} triplets?'  Whereto I was always obliged to reply, with a sad   
shake of the head: `No, I can only {\it add} and subtract them'.''     
The problem, of course, was that there exists no 3-dimensional normed    
division algebra.  He really needed a 4-dimensional algebra.   
   
Finally, on the 16th of October, 1843, while walking with his wife along   
the Royal Canal to a meeting of the Royal Irish Academy in Dublin, he made    
his momentous discovery.  ``That is to say, I then and there felt the    
galvanic circuit of thought {\it close}; and the sparks which fell from it    
were the {\it fundamental equations between $i,j,k$; exactly such} as I have    
used them ever since.''  And in a famous act of mathematical vandalism, he   
carved these equations into the stone of the Brougham Bridge:    
\[    i^2 = j^2 = k^2 = ijk = -1 .\]    
   
One reason this story is so well-known is that Hamilton spent the rest   
of his life obsessed with the quaternions and their applications to   
geometry \cite{Graves,Hankins}.  And for a while, quaternions  were   
fashionable.  They were made a mandatory examination topic in Dublin,   
and in some American universities they were the only advanced   
mathematics taught.  Much of what we now do with scalars and vectors in   
$\R^3$ was  then done using real and imaginary quaternions.   A school   
of `quaternionists' developed, which was led after Hamilton's death by   
Peter Tait of Edinburgh and Benjamin Peirce of Harvard.  Tait wrote 8   
books on the quaternions, emphasizing their applications to physics.    
When Gibbs invented the modern notation for the dot product and cross   
product, Tait condemned it as a ``hermaphrodite monstrosity''.  A war of   
polemics ensued, with such luminaries as Heaviside weighing   
in on the side of vectors.  Ultimately the quaternions lost, and   
acquired a slight taint of disgrace from which they have never fully   
recovered \cite{Crowe}.    
   
Less well-known is the discovery of the octonions by Hamilton's friend   
from college, John T.\ Graves.  It was Graves' interest in algebra that   
got Hamilton thinking about complex numbers and triplets in the first    
place.  The very day after his fateful walk, Hamilton sent an 8-page   
letter describing the quaternions to Graves.  Graves replied on October   
26th, complimenting Hamilton on the boldness of the idea, but adding   
``There is still something in the system which gravels me.  I have not   
yet any clear views as to the extent to which we are at liberty   
arbitrarily to create imaginaries, and to endow them with supernatural   
properties.''  And he asked: ``If with your alchemy you can make three   
pounds of gold, why should you stop there?''  

Graves then set to work on some gold of his own!  On December 26th,
he wrote to Hamilton describing a new 8-dimensional algebra, which he
called the `octaves'.   He showed that they were a normed division
algebra, and used this to express the product of two sums of eight
perfect squares as another sum of eight perfect squares: the `eight
squares theorem' \cite{Hamilton}. 

In January 1844, Graves sent three letters to Hamilton expanding on his   
discovery.  He considered the idea of a general theory of   
`$2^m$-ions', and tried to construct a 16-dimensional normed division   
algebra, but he ``met with an unexpected hitch'' and came to doubt that   
this was possible.  Hamilton offered to publicize Graves' discovery, but   
being busy with work on quaternions, he kept putting it off.  In July he   
wrote to Graves pointing out that the octonions were nonassociative:   
``$A \cdot BC = AB \cdot C = ABC$, if $A,B,C$ be quaternions, but not   
so, generally, with your octaves.''  In fact, Hamilton first invented   
the term `associative' at about this time, so the octonions may have   
played a role in clarifying the importance of this concept.    
   
Meanwhile the young Arthur Cayley, fresh out of Cambridge, had been
thinking about the quaternions ever since Hamilton announced their
existence.  He seemed to be seeking relationships between the
quaternions and hyperelliptic functions.  In March of 1845, he published
a paper in the {\it Philosophical Magazine} entitled `On Jacobi's
Elliptic Functions, in Reply to the Rev.\ B.\ Bronwin; and on Quaternions'
\cite{Cayley}.  The bulk of this paper was an attempt to rebut an
article pointing out mistakes in Cayley's work on elliptic functions.
Apparently as an afterthought, he tacked on a brief description of the
octonions.  In fact, this paper was so full of errors that it was 
omitted from his collected works --- except for the part about octonions
\cite{Cayley2}.

Upset at being beaten to publication, Graves attached a postscript to a
paper of his own which was to appear in the following issue of the same
journal, saying that he had known of the octonions ever since Christmas,
1843.  On June 14th, 1847, Hamilton contributed a short note to the
Transactions of the Royal Irish Academy, vouching for Graves' priority. 
But it was too late: the octonions became known as `Cayley numbers'. 
Still worse, Graves later found that his eight squares theorem had
already been discovered by C.\ F.\ Degen in 1818 \cite{Curtis,Dickson}.

Why have the octonions languished in such obscurity compared to the
quaternions?  Besides their rather inglorious birth, one reason is that
they lacked a tireless defender such as Hamilton.  But surely the reason
for {\it this} is that they lacked any clear application to geometry and
physics.  The unit quaternions form the group $\SU(2)$, which is the
double cover of the rotation group $\SO(3)$.  This makes them nicely
suited to the study of rotations and angular momentum, particularly in
the context of quantum mechanics.  These days we regard this phenomenon
as a special case of the theory of Clifford algebras.  Most of us no
longer attribute to the quaternions the cosmic significance that
Hamilton claimed for them, but they fit nicely into our understanding of
the scheme of things.
 
The octonions, on the other hand, do not.  Their relevance to geometry
was quite obscure until 1925, when \'Elie Cartan described `triality'
--- the symmetry between vectors and spinors in 8-dimensional Euclidean
space \cite{Cartan3}.  Their potential relevance to physics was noticed
in a 1934 paper by Jordan, von Neumann and Wigner on the foundations of
quantum mechanics \cite{JNW}.  However, attempts by Jordan and others to
apply octonionic quantum mechanics to nuclear and particle physics met
with little success.  Work along these lines continued quite slowly
until the 1980s, when it was realized that the octonions explain some
curious features of string theory \cite{KT}.  The Lagrangian for
the classical superstring involves a relationship between vectors and
spinors in Minkowski spacetime which holds only in 3, 4, 6, and 10
dimensions.  Note that these numbers are 2 more than the dimensions of
$\R,\C,\H$ and $\O$.  As we shall see, this is no coincidence: briefly, 
the isomorphisms 
\[ \begin{array}{lcl}
          \Sl(2,\R) &\iso& \so(2,1)   \\
          \Sl(2,\C) &\iso& \so(3,1)   \\   
          \Sl(2,\H) &\iso& \so(5,1)   \\
          \Sl(2,\O) &\iso& \so(9,1)   
\end{array}
\]
allow us to treat a spinor in one of these dimensions as a pair of   
elements of the corresponding division algebra.  It is fascinating   
that of these superstring Lagrangians, it is the 10-dimensional
octonionic one that gives the most promising candidate for a realistic
theory of fundamental physics!  However, there is still no {\it proof}
that the octonions are useful for understanding the real world.  We
can only hope that eventually this question will be settled one way or
another.

Besides their possible role in physics, the octonions are important   
because they tie together some algebraic structures that otherwise   
appear as isolated and inexplicable exceptions.  As we shall explain,   
the concept of an octonionic projective space $\OP^n$ only makes sense   
for $n \le 2$, due to the nonassociativity of $\O$.  This means that   
various structures associated to real, complex and quaternionic   
projective spaces have octonionic analogues only for $n \le 2$.    
   
Simple Lie algebras are a nice example of this phenomenon.   There are  
3 infinite families of `classical' simple Lie algebras, which come from
the isometry groups of the projective spaces $\RP^n$, $\CP^n$ and   
$\HP^n$.  There are also 5 `exceptional' simple Lie algebras.   These 
were discovered by Killing and Cartan in the late 1800s.  At the time,
the significance of these exceptions was shrouded in mystery: they did
not arise as symmetry groups of known structures.  Only later did their
connection to the octonions become clear.  It turns out that 4 of them
come from the isometry groups of the projective planes over $\O$, $\O
\tensor \C$, $\O \tensor \H$ and $\O  \tensor \O$.  The remaining one is
the automorphism group of the octonions!   
   
Another good example is the classification of simple formally real
Jordan algebras.  Besides several infinite families of these, there
is the `exceptional' Jordan algebra, which consists of $3 \times 3$
hermitian octonionic matrices.   Minimal projections in this Jordan
algebra correspond to points of $\OP^2$, and the automorphism group of
this algebra is the same as the isometry group of $\OP^2$.   
   
The octonions also have fascinating connections to topology.  In 1957,  
Raoul Bott computed the homotopy groups of the topological group
$\OO(\infty)$, which is the inductive limit of the orthogonal groups
$\OO(n)$ as $n \to \infty$.  He proved that they repeat with period
8:   
\[   \pi_{i+8}(\OO(\infty)) \iso \pi_i(\OO(\infty)).   \]   
This is known as `Bott periodicity'.  He also computed the first 8:   
\ban       
               \pi_0(\OO(\infty)) &\iso& \Z_2  \\    
               \pi_1(\OO(\infty)) &\iso& \Z_2  \\   
               \pi_2(\OO(\infty)) &\iso&  0    \\   
               \pi_3(\OO(\infty)) &\iso& \Z    \\   
               \pi_4(\OO(\infty)) &\iso&  0    \\   
               \pi_5(\OO(\infty)) &\iso&  0    \\   
               \pi_6(\OO(\infty)) &\iso&  0    \\   
               \pi_7(\OO(\infty)) &\iso& \Z       
\ean   
Note that the nonvanishing homotopy groups here occur in dimensions one   
less than the dimensions of $\R,\C,\H$, and $\O$.  This is no coincidence!   
In a normed division algebra, left multiplication by an element of norm   
one defines an orthogonal transformation of the algebra, and thus an   
element of $\OO(\infty)$.   This gives us maps from the spheres $S^0,    
S^1, S^3$ and $S^7$ to $\OO(\infty)$, and these maps generate the    
homotopy groups in those dimensions.     

Given this, one might naturally guess that the period-8 repetition in
the homotopy groups of $\OO(\infty)$ is in some sense `caused' by the
octonions.  As we shall see, this is true.  Conversely, Bott
periodicity is closely connected to the problem of how many pointwise
linearly independent smooth vector fields can be found on the
$n$-sphere \cite{Husemoller}.  There exist $n$ such vector fields only
when $n+1 = 1, 2, 4,$ or $8$, and this can be used to show that 
division algebras over the reals can only occur in these dimensions.

In what follows we shall try to explain the octonions and their role in 
algebra, geometry, and topology.  In Section \ref{constructing} we give 
four constructions of the octonions: first via their multiplication   
table, then using the Fano plane, then using the Cayley--Dickson   
construction and finally using Clifford algebras, spinors, and a  
generalized concept of `triality' advocated by Frank Adams \cite{Adams}.
Each approach has its own merits.  In Section \ref{proj} we discuss
the projective lines and planes over the normed division algebras --- 
especially $\O$ --- and describe their relation to Bott periodicity,  
the exceptional Jordan algebra, and the Lie algebra isomorphisms listed 
above.  Finally, in Section \ref{lie} we discuss octonionic  
constructions of the exceptional Lie groups, especially the `magic 
square'.   
   
\subsection{Preliminaries} \label{preliminaries}

Before our tour begins, let us settle on some definitions.  For us a
{\bf vector space} will always be a finite-dimensional module over the
field of real numbers.  An {\bf algebra} $A$ will be a vector space that
is equipped with a bilinear map $m \maps A \times A \to A$ called
`multiplication' and a nonzero element $1 \in A$ called the `unit' such
that $m(1,a) = m(a,1) = a$.  As usual, we abbreviate $m(a,b)$ as $ab$.
We do not assume our algebras are associative!  Given an algebra, we
will freely think of real numbers as elements of this algebra via the
map $\alpha \mapsto \alpha 1$.
   
An algebra $A$ is a {\bf division algebra} if given $a,b \in A$ with $ab
= 0$, then either $a = 0$ or $b = 0$.  Equivalently, $A$ is a division
algebra if the operations of left and right multiplication by any
nonzero element are invertible.  A {\bf normed division algebra} is an
algebra $A$ that is also a normed vector space with $\|ab\| = \|a\|
\|b\|$.  This implies that $A$ is a division algebra and that $\|1\| =
1$.

We should warn the reader of some subtleties.  We say an algebra $A$ has
{\bf multiplicative inverses} if for any nonzero $a \in A$ there is an
element $a^{-1} \in A$ with $aa^{-1} = a^{-1}a = 1$.  An associative
algebra has multiplicative inverses iff it is a division
algebra.  However, this fails for nonassociative algebras!  In Section
\ref{cayley-dickson} we shall construct algebras that have
multiplicative inverses, but are not division algebras.  On the other
hand, we can construct a division algebra without multiplicative
inverses by taking the quaternions and modifying the product slightly,
setting $i^2 = -1 + \epsilon j$ for some small nonzero real number
$\epsilon$ while leaving the rest of the multiplication table unchanged.
The element $i$ then has both right and left inverses, but they are not
equal.  (We thank David Rusin for this example.)
   
There are three levels of associativity.  An algebra is {\bf   
power-associative} if the subalgebra generated by any one element is   
associative.  It is {\bf alternative} if the subalgebra generated by any   
two elements is associative.  Finally, if the subalgebra generated by any    
three elements is associative, the algebra is associative.     
   
As we shall see, the octonions are not associative, but they are alternative.   
How does one check a thing like this?  By a theorem of Emil Artin    
\cite{Schafer}, an algebra $A$ is alternative iff for all $a,b \in A$ we have   
\ba  (aa)b = a(ab), \qquad (ab)a = a(ba), \qquad (ba)a = b(aa)    
\label{alternative}   \ea   
In fact, any two of these equations implies the remaining one, so people   
usually take the first and last as the definition of `alternative'.   
To see this fact, note that any algebra has a trilinear map    
\[  [\cdot,\cdot,\cdot] \maps A^3 \to A  \]   
called the {\bf associator}, given by   
\[               [a,b,c] = (ab)c - a(bc)   .\]   
The associator measures the failure of associativity just as the 
commutator $[a,b] = ab - ba$ measures the failure of commutativity.   
Now, the commutator is an alternating bilinear map, meaning that it 
switches sign whenever the two arguments are exchanged:   
\[         [a,b] = -[b,a]   \]   
or equivalently, that it vanishes when they are equal:   
\[         [a,a] = 0 .\]   
This raises the question of whether the associator is alternating too.     
In fact, this holds precisely when $A$ is alternative!  The reason is   
that each equation in (\ref{alternative}) says that the associator   
vanishes when a certain pair of arguments are equal, or equivalently,   
that it changes sign when that pair of arguments is switched.  Note,   
however, that if the associator changes sign when we switch the $i$th   
and $j$th arguments, and also when we switch the $j$th and $k$th   
arguments, it must change sign when we switch the $i$th and $k$th.      
Thus any two of equations (\ref{alternative}) imply the third.     
   
Now we can say what is so great about $\R,\C,\H,$ and $\O$:   
 
\begin{thm}  \et \label{hurwitz} 
$\R,\C,\H$, and $\O$ are the only normed division algebras.  
\end{thm}  
 
\begin{thm} \et \label{zorn}   
$\R,\C,\H$, and $\O$ are the only alternative division algebras.   
\end{thm}   
   
The first theorem goes back to an 1898 paper by Hurwitz \cite{Hurwitz}. 
It was subsequently generalized in many directions, for example, to  
algebras over other fields.   A version of the second theorem appears in
an 1930 paper by Zorn \cite{Zorn} --- the guy with the lemma.  For
modern proofs of both these theorems, see Schafer's excellent book on
nonassociative algebras \cite{Schafer}.  We sketch a couple proofs of
Hurwitz's theorem in Section \ref{clifford}.  

Note that we did {\it not} state that $\R,\C,\H$ and $\O$ are the only   
division algebras.  This is not true.  For example, we have already  
described a way to get 4-dimensional division algebras that do not have  
multiplicative inverses.  However, we do have this fact:   

\begin{thm} \et \label{bott-milnor}  All division algebras have dimension    
$1, 2, 4,$ or $8$.    
\end{thm}  

\noindent
This was independently proved by Kervaire \cite{Kervaire} and
Bott--Milnor \cite{BM} in 1958.  We will say a bit about the proof in
Section \ref{OP1}.  However, in what follows our main focus will not be
on general results about division algebras.  Instead, we concentrate on
special features of the octonions.  Let us begin by constructing them.
  
\section{Constructing the Octonions} \label{constructing}   
   
The most elementary way to construct the octonions is to give their   
multiplication table.  The octonions are an 8-dimensional algebra   
with basis $1, e_1,e_2,e_3,e_4,e_5,e_6,e_7$,    
and their multiplication is given in this table, which describes   
the result of multiplying the element in the $i$th row by the   
element in the $j$th column:   
   
\vskip 2em   
{\vbox{   
\begin{center}   
{\small   
\begin{tabular}{|c|c|c|c|c|c|c|c|c|}                    \hline   
      & $e_1$ & $e_2$ & $e_3$ & $e_4$  & $e_5$ & $e_6$ & $e_7$ \\ \hline   
  $e_1$ & $-1$  & $e_4$ & $e_7$ & $-e_2$ & $e_6$ & $-e_5$ & $-e_3$ \\ \hline   
$e_2$ & $-e_4$ & $-1$ & $e_5$ & $e_1$ & $-e_3$ & $e_7$ & $-e_6$     \\ \hline   
$e_3$ & $-e_7$ & $-e_5$ & $-1$ & $e_6$ & $e_2$ & $-e_4$ & $e_1$   \\ \hline   
$e_4$ & $e_2$ & $-e_1$ & $-e_6$ & $-1$ & $e_7$ & $e_3$ & $-e_5$   \\ \hline   
$e_5$ & $-e_6$ & $e_3$ & $-e_2$ & $-e_7$ & $-1$ & $e_1$ & $e_4$    \\ \hline   
$e_6$ & $e_5$ & $-e_7$ & $e_4$ & $-e_3$ & $-e_1$ & $-1$ & $e_2$     \\ \hline   
$e_7$ & $e_3$ & $e_6$ & $-e_1$ & $e_5$ & $-e_4$ & $-e_2$ & $-1$    \\  \hline   
\end{tabular}} 
\end{center}   
\vskip 1em 
\centerline{Table 1 --- Octonion Multiplication Table}  
}}   
\vskip 1em   
   
\noindent
Unfortunately, this table is almost completely unenlightening!  About the only 
interesting things one can easily learn from it are:  
\begin{itemize}   
\item $e_1,\dots,e_7$ are square roots of -1,    
\item $e_i$ and $e_j$ anticommute when $i \ne j$:    
\[     e_i e_j = -e_j e_i   \]   
\item the {\bf index cycling} identity holds:   
\[       e_i e_j = e_k \; \implies\; e_{i+1} e_{j+1} = e_{k+1}  \]   
where we think of the indices as living in $\Z_7$, and   
\item the {\bf index doubling} identity holds:   
\[       e_i e_j = e_k \; \implies \; e_{2i} e_{2j} = e_{2k} . \]   
\end{itemize}   
Together with a single nontrivial product like $e_1 e_2 = e_4$, these   
facts are enough to recover the whole multiplication table.  However, we   
really want a better way to remember the octonion product.  We should   
become as comfortable with multiplying octonions as we are with   
multiplying matrices!  And ultimately, we want a more conceptual    
approach to the octonions, which explains their special properties and   
how they fit in with other mathematical ideas.  In what follows, we give   
some more descriptions of octonion multiplication, starting with a nice   
mnemonic, and working up to some deeper, more conceptual ones.   
   
\subsection{The Fano plane}    \label{fano}   
   
The quaternions, $\H$, are a 4-dimensional algebra with basis $1,i,j,k$.   
To describe the product we could give a multiplication   
table, but it is easier to remember that:   
\begin{itemize}   
\item $1$ is the multiplicative identity,   
\item $i,j,$ and $k$ are square roots of -1,   
\item we have $ij = k$, $ji = -k$, and all identities obtained    
from these by cyclic permutations of $(i,j,k)$.     
\end{itemize}   
We can summarize the last rule in a picture:   
 
\centerline{\epsfysize=1.5in\epsfbox{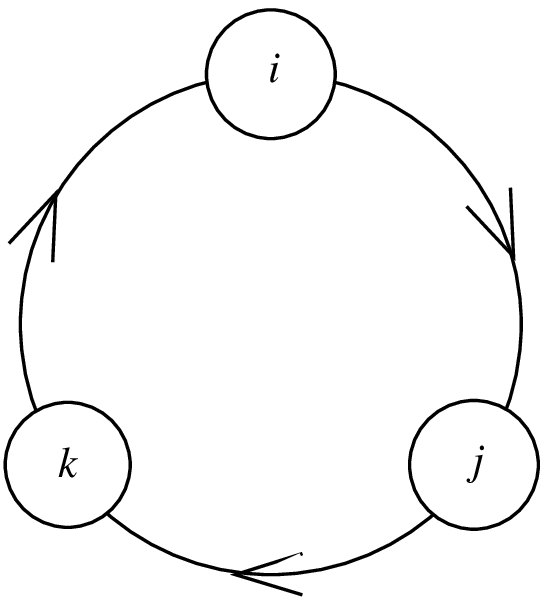}}   
\label{triangle}   
 
\noindent   
When we multiply two elements going clockwise around the circle we get   
the next one: for example, $ij = k$.  But when we multiply two   
going around counterclockwise, we get {\it minus} the next one:    
for example, $ji = -k$.     
   
We can use the same sort of picture to remember how to multiply   
octonions:   
\medskip

\centerline{\epsfysize=1.5in\epsfbox{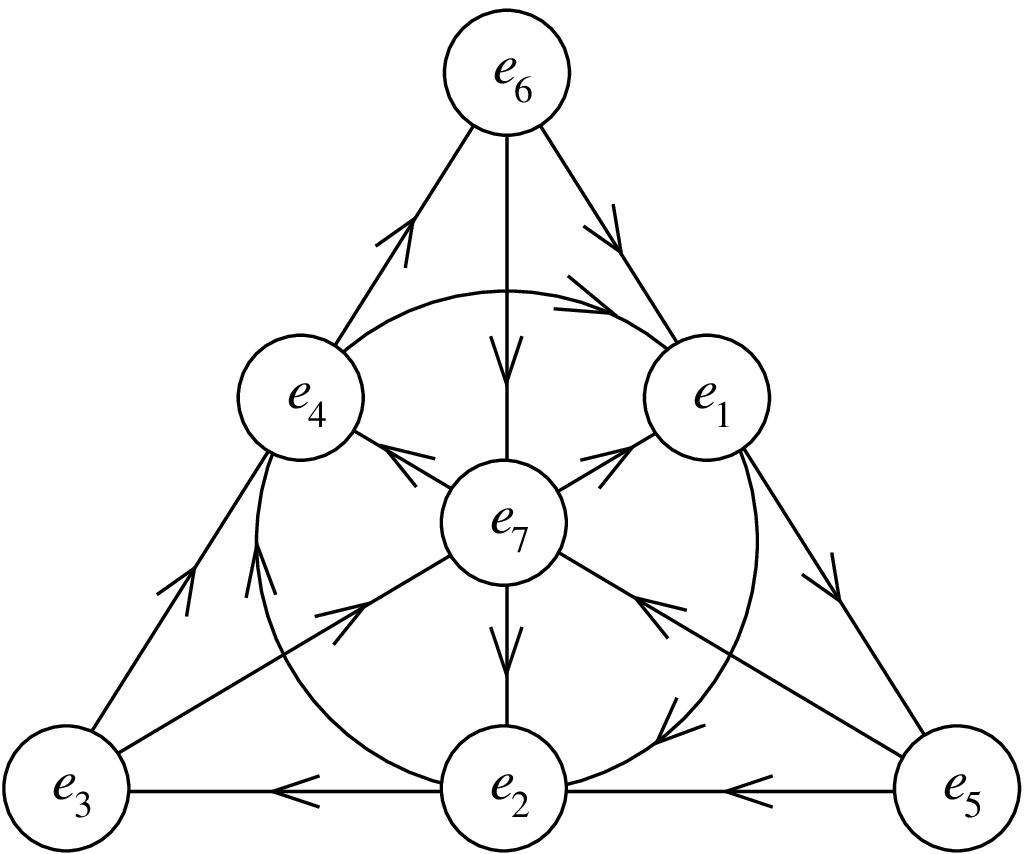}}   
\label{Fano}   
\medskip

\noindent   
This is the {\bf Fano plane}, a little gadget    
with 7 points and 7 lines.  The `lines' are the sides of the triangle,    
its altitudes, and the circle containing all the midpoints of the sides.   
Each pair of distinct points lies on a unique line.  Each line contains    
three points, and each of these triples has has a cyclic ordering    
shown by the arrows.  If $e_i, e_j,$ and $e_k$ are cyclically ordered    
in this way then    
\[            e_i e_j = e_k,  \qquad e_j e_i = -e_k  . \]   
Together with these rules:   
\begin{itemize}   
\item $1$ is the multiplicative identity,   
\item $e_1, \dots, e_7$ are square roots of -1,   
\end{itemize}   
the Fano plane completely describes the algebra structure of the
octonions.   Index-doubling corresponds to rotating the picture
a third of a turn.
   
This is certainly a neat mnemonic, but is there anything deeper lurking   
behind it?  Yes!  The Fano plane is the projective plane over the 2-element
field $\Z_2$.  In other words, it consists of lines through the origin 
in the vector space $\Z_2^3$.  Since every such line contains a single 
nonzero element, we can also think of the Fano plane as consisting of the 
seven nonzero elements of $\Z_2^3$.  If we think of the origin in $\Z_2^3$
as corresponding to $1 \in \O$, we get the following picture of the 
octonions:   

\medskip
\centerline{\epsfysize=1.5in\epsfbox{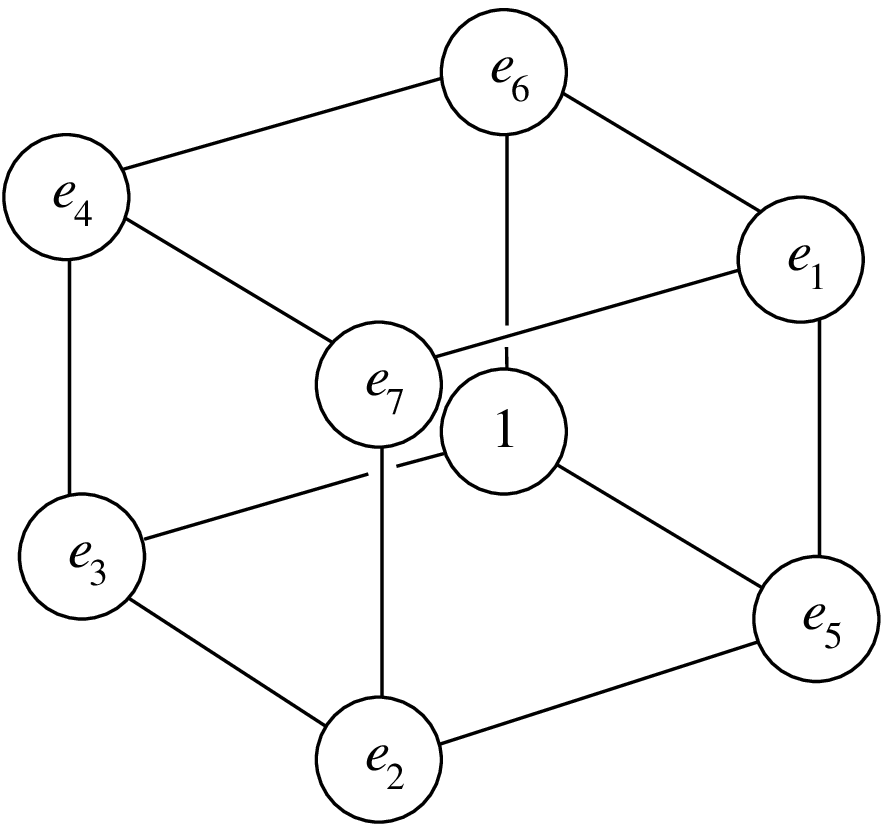}}   
\label{cube}   
\medskip

\noindent    
Note that planes through the origin of this 3-dimensional vector space    
give subalgebras of $\O$ isomorphic to the quaternions, lines through   
the origin give subalgebras isomorphic to the complex numbers, and   
the origin itself gives a subalgebra isomorphic to the real numbers.      
   
What we really have here is a description of the octonions as a
`twisted group algebra'.  Given any group $G$, the group algebra   
$\R[G]$ consists of all finite formal linear combinations of elements   
of $G$ with real coefficients.  This is an associative algebra with   
the product coming from that of $G$.  We can use any function    
\[       \alpha \maps G^2 \to \{ \pm 1 \}     \]   
to `twist' this product, defining a new product   
\[       \star \maps \R[G] \times \R[G] \to \R[G]  \]   
by:   
\[       g \star h = \alpha(g,h) \; gh,  \]   
where $g,h \in G \subset \R[G]$.  One can figure out an equation    
involving $\alpha$ that guarantees this new product will be associative.   
In this case we call $\alpha$ a `2-cocycle'.   If $\alpha$ satisfies  a   
certain extra equation, the product $\star$ will also be commutative,   
and we call $\alpha$ a `stable 2-cocycle'.  For example, the group   
algebra  $\R[\Z_2]$ is isomorphic to a product of 2 copies of $\R$,    
but we can twist it by a stable 2-cocyle to obtain the complex numbers.     
The group algebra $\R[\Z_2^2]$ is isomorphic to a product of 4 copies    
of $\R$, but we can twist it by a 2-cocycle to obtain the quaternions.      
Similarly, the group algebra $\R[\Z_2^3]$ is a product of 8 copies of $\R$,    
and what we have really done in this section is describe a function   
$\alpha$ that allows us to twist this group algebra to obtain the   
octonions.  Since the octonions are nonassociative, this function is   
not a 2-cocycle.  However, its coboundary is a `stable 3-cocycle', which  
allows one to define a new associator and braiding for the category of  
$\Z_2^3$-graded vector spaces, making it into a symmetric monoidal  
category \cite{AM}.  In this symmetric monoidal category, the octonions  
are a commutative monoid object.  In less technical terms: this category 
provides a context in which the octonions {\it are} commutative and 
associative!  So far this idea has just begun to be exploited. 
   
\subsection{The Cayley--Dickson construction}  \label{cayley-dickson}   
   
It would be nice to have a construction of the normed division algebras   
$\R,\C,\H,\O$ that explained why each one fits neatly inside the next.   
It would be nice if this construction made it clear why $\H$ is   
noncommutative and $\O$ is nonassociative.  It would be even better if   
this construction gave an infinite sequence of algebras, doubling in   
dimension each time, with the normed division algebras as the first   
four.  In fact, there is such a construction: it's called the   
Cayley--Dickson construction.   
   
As Hamilton noted, the complex number $a+bi$ can be thought of as a pair
$(a,b)$ of real numbers.  Addition is done component-wise, and
multiplication goes like this:
\[            (a,b)(c,d) = (ac - db,ad + cb)   \]   
We can also define the conjugate of a complex number by   
\[                (a,b)^* = (a,-b).   \]   
   
Now that we have the complex numbers, we can define the   
quaternions in a similar way.  A quaternion can be thought of    
as a pair of complex numbers.  Addition is done component-wise,   
and multiplication goes like this:   
\be          (a,b)(c,d) = (ac - db^*, a^* d + cb)      \label{cd1} \ee   
This is just like our formula for multiplication of complex numbers, but
with a couple of conjugates thrown in.  If we included them in 
the previous formula nothing would change, since the conjugate of a
real number is just itself.  We can also define the conjugate of a
quaternion by
\be          (a,b)^* = (a^*,-b).        \label{cd2} \ee   
   
The game continues!  Now we can define an octonion to be a pair of
quaternions.  As before, we add and multiply them using formulas
(\ref{cd1}) and (\ref{cd2}).  This trick for getting new algebras from
old is called the {\bf Cayley--Dickson construction}.
   
Why do the real numbers, complex numbers, quaternions   
and octonions have multiplicative inverses?  I take it as   
obvious for the real numbers.  For the complex numbers,    
one can check that   
\[    (a,b) (a,b)^* = (a,b)^* (a,b) = k (1,0)  \]   
where $k$ is a real number, the square of the norm of $(a,b)$.     
This means that whenever $(a,b)$ is nonzero, its multiplicative    
inverse is $(a,b)^*/k$.  One can check that the same holds for the    
quaternions and octonions.   
   
But this, of course, raises the question: why isn't there an {\it
infinite} sequence of division algebras, each one obtained from the
preceding one by the Cayley--Dickson construction?  The answer is that
each time we apply the construction, our algebra gets a bit worse.
First we lose the fact that every element is its own conjugate, then we
lose commutativity, then we lose associativity, and finally we lose the
division algebra property.
   
To see this clearly, it helps to be a bit more formal.  Define a {\bf   
$\ast$-algebra} to be an algebra $A$ equipped with a {\bf conjugation},   
that is, a real-linear map $\ast \maps A \to A$ with   
\[            a^{**} = a, \quad \quad (ab)^* = b^* a^*  \]   
for all $a,b \in A$.  We say a $\ast$-algebra is {\bf real} if  $a =   
a^*$ for every element $a$ of the algebra.  We say the $\ast$-algebra   
$A$ is {\bf nicely normed} if $a + a^* \in \R$ and $aa^* = a^* a > 0$ 
for all nonzero $a \in A$.  If $A$ is nicely normed we set   
\[             \Re(a) = (a + a^\ast)/2  \in \R, \qquad   
               \Im(a) = (a - a^\ast)/2 ,  \]   
and define a norm on $A$ by    
\[                 \|a\|^2 = aa^\ast .      \]   
If $A$ is nicely normed, it has multiplicative inverses given by   
\[               a^{-1} = a^\ast / \|a\|^2   .\]   
If $A$ is nicely normed and alternative, $A$ is a normed division   
algebra.  To see this, note that for any $a,b \in A$, all 4 elements   
$a,b,a^\ast,b^\ast$ lie in the associative algebra generated by $\Im(a)$   
and $\Im(b)$, so that    
\[     \|ab\|^2 = (ab)(ab)^\ast = ab(b^\ast a^\ast) =    
a(bb^\ast)a^\ast = \|a\|^2 \|b\|^2 . \]   
   
Starting from any $\ast$-algebra $A$, the Cayley--Dickson construction   
gives a new $\ast$-algebra $A'$.   Elements of $A'$ are pairs $(a,b) \in   
A^2$, and multiplication and conjugation are defined using equations   
(\ref{cd1}) and (\ref{cd2}).   The following propositions show the   
effect of repeatedly applying the Cayley--Dickson construction:   
   
\begin{prop} \et \label{CD1}   
$A'$ is never real.     
\end{prop}   
   
\begin{prop} \et \label{CD2}   
$A$ is real (and thus commutative) $\iff$ $A'$ is commutative.   
\end{prop}     
   
\begin{prop} \et \label{CD3}   
$A$ is commutative and associative $\iff$ $A'$ is associative.    
\end{prop}   
   
\begin{prop} \et \label{CD4}   
$A$ is associative and nicely normed $\iff$   
$A'$ is alternative and nicely normed.   
\end{prop}      
   
\begin{prop} \et\label{CD5}   
$A$ is nicely normed $\iff$ $A'$ is nicely normed.    
\end{prop}      
   
\noindent All of these follow from straightforward calculations; to   
prove them here would merely deprive the reader of the pleasure of doing   
so.  It follows from these propositions that:   
\begin{center}   
$\R$ is a real commutative associative nicely normed 
$\ast$-algebra $\implies$   
    
$\C$ is a commutative associative nicely normed $\ast$-algebra $\implies$   
   
$\H$ is an associative nicely normed $\ast$-algebra $\implies$   
   
$\O$ is an alternative nicely normed $\ast$-algebra    
\end{center}   
and therefore that $\R,\C,\H,$ and $\O$ are normed division algebras.     
It also follows that the octonions are neither real, nor commutative, nor   
associative.       
   
If we keep applying the Cayley--Dickson process to the octonions we get a   
sequence of $\ast$-algebras of dimension 16, 32, 64, and so on.   The   
first of these is called the {\bf sedenions}, presumably alluding to the   
fact that it is 16-dimensional \cite{LPS}.   It follows from the above   
results that all the $\ast$-algebras in this sequence are nicely normed   
but neither real, nor commutative, nor alternative.  They all have   
multiplicative inverses, since they are nicely normed.  But they are not   
division algebras, since an explicit calculation demonstrates that the   
sedenions, and thus all the rest, have zero divisors.   In fact    
\cite{Cohen,Moreno}, the zero divisors of norm one in the sedenions    
form a subspace that is homeomorphic to the exceptional Lie group $\G_2$.   
   
The Cayley--Dickson construction provides a nice way to obtain the
sequence $\R,\H,\C,\O$ and the basic properties of these algebras.  But
what is the meaning of this construction?  To answer this, it is better
to define $A'$ as the algebra formed by adjoining to $A$ an element $i$
satisfying $i^2 = -1$ together with the following relations:
\be  a(ib) = i(a^* b) ,   \qquad   
    (ai)b = (ab^*)i,      \qquad   
    (ia)(bi^{-1}) = (ab)^*  \label{cd3} \ee   
for all $a,b \in A$.  We make $A'$ into a $\ast$-algebra using the   
original conjugation on elements of $A$ and setting $i^* = -i$.  It is   
easy to check that every element of $A'$ can be uniquely written as $a +   
ib$ for some $a,b \in A$, and that this description of the   
Cayley--Dickson construction becomes equivalent to our previous one    
if we set $(a,b) = a + ib$.     
   
What is the significance of the relations in (\ref{cd3})?   Simply   
this: {\sl they express conjugation in terms of conjugation!}  This is a pun   
on the double meaning of the word `conjugation'.  What I really mean is   
that they express the $\ast$ operation in $A$ as conjugation by $i$.  In   
particular, we have   
\[       a^\ast = (ia)i^{-1} = i(ai^{-1})   \]   
for all $a \in A$.  Note that when $A'$ is associative, any one of the   
relations in (\ref{cd3}) implies the other two.  It is when $A'$ is   
nonassociative that we really need all three relations.   

This interpretation of the Cayley--Dickson construction makes it easier   
to see what happens as we repeatedly apply the construction starting with   
$\R$.  In $\R$ the $\ast$ operation does nothing, so when we do the   
Cayley--Dickson construction, conjugation by $i$ must have no effect on   
elements of $\R$.  Since $\R$ is commutative, this means that $\C = \R'$   
is commutative.  But $\C$ is no longer real, since $i^* = -i$.   
   
Next let us apply the Cayley--Dickson construction to $\C$.  Since $\C$   
is commutative, the $\ast$ operation in $\C$ is an automorphism.  Whenever   
we have an associative algebra $A$ equipped with an automorphism $\alpha$,    
we can always extend $A$ to a larger associative algebra by adjoining an    
invertible element $x$ with    
\[                \alpha(a) = xax^{-1}  \]    
for all $a \in A$.  Since $\C$ is associative, this means that $\C' =   
\H$ is associative.  But since $\C$ is not real, $\H$ cannot be   
commutative, since conjugation by the newly adjoined element $i$ must   
have a nontrivial effect.   
   
Finally, let us apply the Cayley--Dickson construction to $\H$.   Since   
$\H$ is noncommutative, the $\ast$ operation in $\H$ is not an   
automorphism; it is merely an antiautomorphism.   This means we cannot   
express it as conjugation by some element of a larger associative   
algebra.  Thus $\H' = \O$ must be nonassociative.     

\subsection{Clifford Algebras} \label{clifford}      

William Clifford invented his algebras in 1876 as an attempt to
generalize the quaternions to higher dimensions, and he published a
paper about them two years later \cite{Clifford}.  Given a real inner
product space $V$, the {\bf Clifford algebra} $\Cliff(V)$ is the
associative algebra freely generated by $V$ modulo the relations
\[  v^2 = -\|v\|^2 \] 
for all $v \in V$.  Equivalently, it is the associative algebra 
freely generated by $V$ modulo the relations 
\[              vw + wv = -2\langle v,w\rangle    \]  
for all $v, w \in V$.  If $V = \R^n$ with its usual inner product, we 
call this Clifford algebra $\Cliff(n)$.  Concretely, this is the 
associative algebra freely generated by $n$ anticommuting square roots 
of $-1$.  From this we easily see that  
\[  \Cliff(0) = \R, \qquad\qquad \Cliff(1) = \C, \qquad\qquad    
\Cliff(2) = \H .\]   
So far this sequence resembles the iterated Cayley-Dickson construction 
--- but the octonions are {\it not} a Clifford algebra, since they are 
nonassociative.   Nonetheless, there is a profound relation between  
Clifford algebras and normed division algebras.  This relationship gives
a nice way to prove that $\R, \C, \H$ and $\O$ are the only normed  
dvivision algebras.  It is also crucial for understanding the  
geometrical meaning of the octonions.    
  
To see this relation, first suppose $\K$ is a normed division algebra.    
Left multiplication by any element $a \in \K$ gives an operator  
\[   \begin{array}{ccccc}  
       L_a &\maps& \K & \to   &  \K  \\  
           &     & x &\mapsto&   ax .   
\end{array} \]  
If $\|a\| = 1$, the operator $L_a$ is norm-preserving, so it maps the  
unit sphere of $\K$ to itself.  Since $\K$ is a division algebra, we can  
find an operator of this form mapping any point on the unit sphere to  
any other point.  The only way the unit sphere in $\K$ can have this much  
symmetry is if the norm on $\K$ comes from an inner product.  Even better,  
this inner product is unique, since we can use the polarization identity  
\[   \langle x, y\rangle = {1\over 2}(\|x+y\|^2 - \|x\|^2 - \|y\|^2) \]   
to recover it from the norm.     
  
Using this inner product, we say an element $a \in \K$ is {\bf imaginary}  
if it is orthogonal to the element $1$, and we let $\Im(\K)$ be the space  
of imaginary elements of $\K$.  We can also think of $\Im(\K)$ as the  
tangent space of the unit sphere in $\K$ at the point $1$.  This has a  
nice consequence: since $a \mapsto L_a$ maps the unit sphere in $\K$ to  
the Lie group of orthogonal transformations of $\K$, it must send  
$\Im(\K)$ to the Lie algebra of skew-adjoint transformations of $\K$.  
In short, $L_a$ is skew-adjoint whenever $a$ is imaginary.  
  
The relation to Clifford algebras shows up when we compute the square of  
$L_a$ for $a \in \Im(\K)$.  We can do this most easily when $a$ has norm  
$1$.  Then $L_a$ is both orthogonal and skew-adjoint.  For any  
orthogonal transformation, we can find some orthonormal basis in which  
its matrix is block diagonal, built from $2 \times 2$ blocks that look  
like this:  
\[   \left( \begin{array}{cc} \cos \theta & \sin \theta \\  
                            -\sin \theta & \cos \theta   
            \end{array} \right)  \]  
and possibly a $1 \times 1$ block like this: $\left( 1 \right)$.  
Such a transformation can only be skew-adjoint if it consists solely of  
$2 \times 2$ blocks of this form:  
\[  \pm \left( \begin{array}{cc} 0 & 1 \\  
                             -1 & 0   
            \end{array} \right).  \]  
In this case, its square is $-1$.  We thus have $L_a^2 = -1$ when  
$a \in \Im(\K)$ has norm 1.  It follows that 
\[             L_a^2 = -\| a \|^2 \]  
for all $a \in \Im(\K)$.  We thus obtain a representation of the Clifford  
algebra $\Cliff(\Im(\K))$ on $\K$.   Any $n$-dimensional normed division  
algebra thus gives an $n$-dimensional representation of $\Cliff(n-1)$.    
As we shall see, this is very constraining.    
  
We have already described the Clifford algebras up to $\Cliff(2)$.  
Further calculations \cite{Harvey,Porteous} give the following table,  
where we use $A[n]$ to stand for $n\times n$ matrices with entries in  
the algebra $A$:  
   
\medskip  
{\vbox{   
\begin{center}   
{\small   
\begin{tabular}{|c|c|c|c|c|c|c|c|c|}                    \hline   
$n$ &    $\Cliff(n)$          \\ \hline   
$0$ &    $\R$                 \\ \hline   
$1$ &    $\C$                 \\ \hline   
$2$ &    $\H$                 \\ \hline   
$3$ &    $\H \oplus \H$       \\ \hline   
$4$ &    $\H[2]$              \\ \hline   
$5$ &    $\C[4]$              \\ \hline   
$6$ &    $\R[8]$              \\ \hline   
$7$ &    $\R[8] \oplus \R[8]$ \\ \hline   
\end{tabular}}  
\vskip 1em 
\centerline{Table 2 --- Clifford Algebras}  
\end{center}   
}}   
\medskip  
   
\noindent  
Starting at dimension 8, something marvelous happens: the table continues  
in the following fashion:   
\[    \Cliff(n+8) \iso \Cliff(n) \tensor \R[16]  .\]   
In other words, $\Cliff(n+8)$ consists of $16 \times 16$ matrices  
with entries in $\Cliff(n)$.  This `period-8' behavior was discovered  
by Cartan in 1908 \cite{Cartan2}, but we will take the liberty of 
calling it {\bf Bott periodicity}, since it has a far-ranging set of 
applications to topology, some of which were discovered by Bott.

Since Clifford algebras are built from matrix algebras over $\R,\C$ and
$\H$, it is easy to determine their representations.  Every
representation is a direct sum of irreducible ones.  The only
irreducible representation of $\R[n]$ is its obvious one via matrix
multiplication on $\R^n$.  Similarly, the only irreducible
representation of $\C[n]$ is the obvious one on $\C^n$, and the only
irreducible representation of $\H[n]$ is the obvious one on $\H^n$.
  
Glancing at the above table, we see that unless $n$ equals $3$ or $7$
modulo $8$, $\Cliff(n)$ is a real, complex or quaternionic matrix
algebra, so it has a unique irreducible representation.  For reasons to
be explained later, this irreducible representation is known as the
space of {\bf pinors} and denoted $P_n$.  When $n$ is $3$ or $7$ modulo
$8$, the algebra $\Cliff(n)$ is a direct sum of two real or quaternionic
matrix algebras, so it has two irreducible representations, which we
call the {\bf positive pinors} $P_n^+$ and {\bf negative pinors}
$P_n^-$.  We summarize these results in the following table:
  
\medskip  
{\vbox{   
\begin{center}   
{\small   
\begin{tabular}{|c|c|l|}                            \hline   
$n$ &    $\Cliff(n)$         & irreducible representations       \\ \hline
$0$ &    $\R$                & $P_0 = \R$                        \\ \hline   
$1$ &    $\C$                & $P_1 = \C$                        \\ \hline   
$2$ &    $\H$                & $P_2 = \H$                        \\ \hline   
$3$ &    $\H \oplus \H$      & $P^+_3 = \H,\, P^-_3 =\H$         \\ \hline   
$4$ &    $\H[2]$             & $P_4 = \H^2$                      \\ \hline   
$5$ &    $\C[4]$             & $P_5 = \C^4$                      \\ \hline   
$6$ &    $\R[8]$             & $P_6 =\R^8$                       \\ \hline   
$7$ &    $\R[8] \oplus \R[8]$& $P_7^+ = \R^8,\, P_7^- =\R^8$     \\ \hline   
\end{tabular}}  
\vskip 1em 
\centerline{Table 3 --- Pinor Representations}  
\end{center}   
}}   
\medskip  
  
\noindent 
Examining this table, we see that in the range of dimensions listed
there is an $n$-dimensional representation of $\Cliff(n-1)$ only for $n
= 1,2,4,$ and $8$.  What about higher dimensions?  By Bott periodicity,
the irreducible representations of $\Cliff(n+8)$ are obtained by
tensoring those of $\Cliff(n)$ by $\R^{16}$.  This multiplies the
dimension by 16, so one can easily check that for $n > 8$, the
irreducible representations of $\Cliff(n-1)$ always have dimension
greater than $n$.  

It follows that normed division algebras are only possible in dimensions  
$1,2,4,$ and $8$.  Having constructed $\R,\C,\H$ and $\O$, we also know  
that normed division algebras {\it exist} in these dimensions.  The only  
remaining question is whether they are {\it unique}.  For this it helps  
to investigate more deeply the relation between normed division algebras  
and the Cayley-Dickson construction.   In what follows, we outline an  
approach based on ideas in the book by Springer and Veldkamp \cite{SV}.   
  
First, suppose $\K$ is a normed division algebra.  Then there is a unique  
linear operator $\ast \maps \K \to \K$ such that $1^\ast = 1$ and $a^\ast  
= -a$ for $a \in \Im(\K)$.  With some calculation one can prove this  
makes $\K$ into a nicely normed $\ast$-algebra.    
  
Next, suppose that $\K_0$ is any subalgebra of the normed division algebra  
$\K$.  It is easy to check that $\K_0$ is a nicely normed $\ast$-algebra in  
its own right.  If $\K_0$ is not all of $\K$, we can find an element $i \in 
\K$ that is orthogonal to every element of $\K_0$.   Without loss of 
generality we shall assume this element has norm 1.  Since this element 
$i$ is orthogonal to $1 \in \K_0$, it is imaginary.  From the  definition 
of the $\ast$ operator it follows that $i^\ast = -i$, and from results 
earlier in this section we have $i^2 = -1$.  With  further calculation 
one can show that for all $a,a' \in \K_0$ we have   
\[  a(ia') = i(a^* a') ,  \qquad   
    (ai)a' = (aa'^*)i,      \qquad   
    (ia)(a'i^{-1}) = (aa')^*  \]  
A glance at equation (\ref{cd3}) reveals that these are exactly the  
relations defining the Cayley-Dickson construction!  With a little  
thought, it follows that the subalgebra of $\K$ generated by $\K_0$ and $i$  
is isomorphic as a $\ast$-algebra to $\K'_0$, the $\ast$-algebra obtained  
from $\K_0$ by the Cayley-Dickson construction.  
  
Thus, whenever we have a normed division algebra $\K$ we can find a 
chain of subalgebras $\R = \K_0 \subset \K_1 \subset \cdots \subset 
\K_n = \K$ such that $\K_{i+1} \iso \K_i'$.  To construct $\K_{i+1}$, we 
simply need to choose a norm-one element of $\K$ that is orthogonal to 
every element of $\K_i$.    It follows that the only normed division 
algebras of dimension 1, 2, 4 and 8 are $\R,\C,\H$ and $\O$.   This also 
gives an alternate proof that there are no normed division algebras of 
other dimensions: if there were any, there would have to be a 
16-dimensional one, namely $\O'$ --- the sedenions.  But as mentioned
in Section \ref{cayley-dickson}, one can check explicitly that the 
sedenions are not a division algebra.    

\subsection{Spinors and Trialities}  \label{triality}   
   
A nonassociative division algebra may seem like a strange thing to   
bother with, but the notion of triality makes it seem a bit more   
natural.  The concept of duality is important throughout linear algebra.   
The concept of triality is similar, but considerably subtler.  Given   
vector spaces $V_1$ and $V_2$, we may define a {\bf duality} to be a   
bilinear map    
\[    f \maps V_1 \times V_2 \to \R   \]   
that is nondegenerate, meaning that if we fix either argument   
to any nonzero value, the linear functional induced on the other vector   
space is nonzero.  Similarly, given vector spaces $V_1,V_2,$ and $V_3$,   
a {\bf triality} is a trilinear map    
\[    t \maps V_1 \times V_2 \times V_3 \to \R   \]   
that is nondegenerate in the sense that if we fix any two arguments to   
any nonzero values, the linear functional induced on the third vector   
space is nonzero.   
   
Dualities are easy to come by.  Trialities are much rarer.  For suppose   
we have a triality  
\[    t \maps V_1 \times V_2 \times V_3 \to \R  . \]   
By dualizing, we can turn this into a bilinear map   
\[     m \maps V_1 \times V_2 \to V_3^\ast     \]   
which we call `multiplication'.  By the nondegeneracy of our triality,   
left multiplication by any nonzero element of $V_1$ defines an   
isomorphism from $V_2$ to $V_3^\ast$.  Similarly, right multiplication   
by any nonzero element of $V_2$ defines an isomorphism from $V_1$ to   
$V_3^\ast$.  If we choose nonzero elements $e_1 \in V_1$ and $e_2 \in   
V_2$, we can thereby identify the spaces $V_1$, $V_2$ and $V_3^\ast$   
with a single vector space, say $V$.   Note that this identifies   
all three vectors $e_1 \in V_1$, $e_2 \in V_2$, and $e_1e_2 \in V_3^\ast$   
with the same vector $e \in V$.   We thus obtain a product   
\[      m \maps V \times V \to V        \]   
for which $e$ is the left and right unit.  Since left or right   
multiplication by any nonzero element is an isomorphism, $V$ is   
actually a division algebra!   Conversely, any division algebra   
gives a triality.     

It follows from Theorem \ref{bott-milnor}  that trialities only occur in
dimensions 1, 2, 4, or 8.  This theorem is quite deep.  By comparison,
Hurwitz's classification of {\it normed} division algebras is easy to
prove.  Not surprisingly, these correspond to a special sort of
triality, which we call a `normed' triality.  
  
To be precise, a {\bf normed triality} consists of inner product  
spaces $V_1, V_2, V_3$ equipped with a trilinear map   
$t \maps V_1 \times V_2 \times V_3 \to \R$ with  
\[      |t(v_1, v_2, v_3)| \le \|v_1\| \, \|v_2\| \, \|v_3 \|,   \]    
and such that for all $v_1, v_2$ there exists $v_3 \ne 0$ for which this   
bound is attained --- and similarly for cyclic permutations of $1,2,3$.    
Given a normed triality, picking unit vectors in any two of the spaces 
$V_i$ allows us to identify all three spaces and get a normed division 
algebra.  Conversely, any normed division algebra gives a normed triality.    
  
But where do normed trialities come from?  They come from the theory of  
spinors!  From Section \ref{clifford}, we already know that any  
$n$-dimensional normed division algebra is a representation of  
$\Cliff(n-1)$, so it makes sense to look for normed trialities here. 
In fact, representations of $\Cliff(n-1)$ give certain representations  
of $\Spin(n)$, the double cover of the rotation group in $n$ dimensions.  
These are called `spinors'.  As we shall see, the relation between   
spinors and vectors gives a nice way to construct normed trialities in  
dimensions 1, 2, 4 and 8.  
  
To see how this works, first let $\Pin(n)$ be the group sitting inside  
$\Cliff(n)$ that consists of all products of unit vectors in $\R^n$.  
This group is a double cover of the orthogonal group $\OO(n)$, where  
given any unit vector $v \in \R^n$, we map both $\pm v \in \Pin(n)$ to  
the element of $\OO(n)$ that reflects across the hyperplane perpendicular to  
$v$.  Since every element of $\OO(n)$ is a product of reflections, this  
homomorphism is indeed onto.      
  
Next, let $\Spin(n) \subset \Pin(n)$ be the subgroup consisting of all  
elements that are a product of an even number of unit vectors in  
$\R^n$.  An element of $\OO(n)$ has determinant 1 iff it is the product  
of an even number of reflections, so just as $\Pin(n)$ is a double cover  
of $\OO(n)$, $\Spin(n)$ is a double cover of $\SO(n)$.  Together with a  
French dirty joke which we shall not explain, this analogy is the origin  
of the terms `$\Pin$' and `pinor'.    
  
Since $\Pin(n)$ sits inside $\Cliff(n)$, the irreducible representations
of $\Cliff(n)$ restrict to representations of $\Pin(n)$, which turn out
to be still irreducible.  These are again called {\bf pinors}, and we
know what they are from Table 3.  Similarly, $\Spin(n)$ sits inside the
subalgebra
\[            \Cliff_0(n) \subseteq \Cliff(n) \]  
consisting of all linear combinations of products of an even number of
vectors in $\R^n$.  Thus the irreducible representations of
$\Cliff_0(n)$ restrict to representations of $\Spin(n)$, which turn out
to be still irreducible.  These are called {\bf spinors} --- but we warn
the reader that this term is also used for many slight variations on
this concept.
  
In fact, there is an isomorphism   
\[      \phi \maps \Cliff(n-1) \to \Cliff_0(n) \]  
given as follows:  
\[       \phi(e_i) = e_i e_n  ,   \qquad \qquad 1 \le i \le n-1 ,\]  
where $\{e_i\}$ is an orthonormal basis for $\R^n$.  Thus spinors in $n$  
dimensions are the same as pinors in $n-1$ dimensions!   Table 3  
therefore yields the following table, where we use similar notation but  
with `$S$' instead of `$P$':  
  
\medskip  
{\vbox{   
\begin{center}   
{\small   
\begin{tabular}{|c|c|l|}                           \hline   
$n$ &  $\Cliff_0(n)$    & irreducible representations            \\ \hline  
$1$ &  $\R$             & $S_1 = \R$                             \\ \hline   
$2$ &  $\C$             & $S_2 = \C$                             \\ \hline   
$3$ &  $\H$             & $S_3 = \H$                             \\ \hline   
$4$ &  $\H \oplus \H$   & $S_4^+ = \H, \, S_4^- = \H$            \\ \hline   
$5$ &  $\H[2]$          & $S_5 = \H^2$                           \\ \hline   
$6$ &  $\C[4]$          & $S_6 = \C^4$                           \\ \hline   
$7$ &  $\R[8]$          & $S_7 = \R^8$                           \\ \hline   
$8$ &  $\R[8] \oplus \R[8]$ & $S_8^+ = \R^8,\, S_8^- = \R^8$     \\ \hline   
\end{tabular}}  
\vskip 1em 
\centerline{Table 4 --- Spinor Representations}  
\end{center}   
}}   
\medskip   
  
\noindent   
We call $S_n^+$ and $S_n^-$ the {\bf right-handed} and {\bf left-handed}  
spinor representations.  For $n > 8$ we can work out the spinor  
representations using Bott periodicity:  
\[      S_{n+8} \iso S_n \tensor \R^{16}   \]  
and similarly for right-handed and left-handed spinors.

Now, besides its pinor representation(s), the group $\Pin(n)$ also has
an irreducible representation where we first apply the 2--1 homomorphism
$\Pin(n) \to \OO(n)$ and then use the obvious representation of $\OO(n)$
on $\R^n$.  We call this the {\bf vector} representation, $V_n$.  As a
vector space $V_n$ is just $\R^n$, and $\Cliff(n)$ is generated by
$\R^n$, so we have an inclusion
\[                  V_n \hookrightarrow \Cliff(n)  .\]  
Using this, we can restrict the action of the Clifford algebra on pinors  
to a map  
\[    
\begin{array}{lcc}       
m_n \maps& V_n \times P_n^\pm \to P_n^\pm &n \equiv 3,7 \bmod 8  \\  
m_n \maps& V_n \times P_n \to P_n         & {\rm otherwise.}   
\end{array}   
\]  
This map is actually an intertwining operator between 
representations of $\Pin(n)$.  If we restrict the vector representation 
to the subgroup $\Spin(n)$, it remains irreducible.  This is not always 
true for the pinor representations, but we can always decompose them as 
a direct sum of spinor representations.  Applying this decomposition to 
the map $m_n$, we get a map 
\[    
\begin{array}{lc}        
m_n \maps V_n \times S_n^\pm \to S_n^\mp   &n \equiv 0,4 \bmod 8  \\   
m_n \maps V_n \times S_n \to S_n         & {\rm otherwise.}    
\end{array}    
\]   
All the spinor representations appearing here are self-dual, so we can  
dualize the above maps and reinterpret them as trilinear maps   
\[  
\begin{array}{lc}        
t_n \maps V_n \times S_n^+ \times S_n^- \to \R  & n \equiv 0,4 \bmod 8  \\   
t_n \maps V_n \times S_n \times S_n \to \R      & {\rm otherwise.}    
\end{array}   
\]   
 
These trilinear maps are candidates for trialities!  However, they can  
only be trialities when the dimension of the vector representation  
matches that of the relevant spinor representations.    In the range of  
the above table this happens only for $n = 1,2,4,8$.  In these cases we 
actually do get normed trialities, which in turn give normed division algebras: 
\[ 
\begin{array}{ll} 
  t_1 \maps V_1 \times S_1 \times S_1 \to \R  & 
{\rm \; gives \;} \R.   \\ 
   t_2 \maps V_2 \times S_2 \times S_2 \to \R & 
{\rm \; gives \;} \C.    \\ 
   t_4 \maps V_4 \times S_4^+ \times S_4^- \to \R & 
{\rm \; gives \;} \H.    \\ 
   t_8 \maps V_8 \times S_8^+ \times S_8^- \to \R & 
{\rm \; gives \;} \O .  
\end{array} 
\] 
In higher dimensions, the spinor representations become bigger than the 
vector representation, so we get no more trialities this way --- and of 
course, none exist.  
 
Of the four normed trialities, the one that gives the octonions 
has an interesting property that the rest lack.  To see this property, 
one must pay careful attention to the difference between a normed triality 
and a normed division algebra.  To construct a normed division $\K$ 
algebra from the normed triality $t \maps V_1 \times V_2 \times V_3 \to \R$, 
we must arbitrarily choose unit vectors in two of the three spaces, so 
the symmetry group of $\K$ is smaller than that of $t$.  More precisely,  
let us define a {\bf automorphism} of the normed triality $t \maps V_1 \times 
V_2 \times V_3 \to \R$ to be a triple of norm--preserving maps  
$f_i \maps V_i \to V_i$ such that  
\[        t(f_1(v_1), f_2(v_2), f_3(v_3)) = t(v_1,v_2,v_3)  \] 
for all $v_i \in V_i$.   These automorphisms form a group we call  
$\Aut(t)$.  If we construct a normed division algebra $\K$ from $t$ 
by choosing unit vectors $e_1 \in V_1, e_2 \in V_2$, we have  
\[   
\Aut(\K) \iso \{(f_1,f_2,f_3) \in \Aut(t)\; \colon \; f_1(e_1) = e_1,  
\; f_2(e_2) = e_2 \} . 
\] 
 
In particular, it turns out that: 
\be
\begin{array}{lclclcl} 
    1 &\iso& \Aut(\R) & \subseteq& \Aut(t_1) & \iso &
\{ (g_1,g_2,g_3) \in \OO(1)^3 \colon \; g_1g_2g_3 = 1 \}  \\        
    \Z_2 &\iso& \Aut(\C) & \subseteq& \Aut(t_2)& \iso &
\{ (g_1,g_2,g_3) \in \U(1)^3 \colon \; g_1g_2g_3 = 1 \} \times \Z_2  \\        
    \SO(3) &\iso& \Aut(\H) &\subseteq &\Aut(t_4)& \iso  &
\Sp(1)^3 / \{ \pm(1,1,1) \} \\ 
     \G_2 &\iso& \Aut(\O) &\subseteq& \Aut(t_8)& \iso& \Spin(8)  
\end{array} 
\label{Aut(t)}
\ee
where 
\[ \OO(1) \iso \Z_2, \qquad \U(1) \iso \SO(2) , \qquad \Sp(1) \iso \SU(2) \] 
are the unit spheres in $\R$, $\C$ and $\H$, respectively ---
the only spheres that are Lie groups.
$\G_2$ is just another name for the automorphism group of 
the octonions; we shall study this group in Section \ref{G2}.   
The bigger group $\Spin(8)$ acts as automorphisms of the triality  
that gives the octonions, and it does so in an interesting way.
Given any element $g \in \Spin(8)$, there exist unique elements
$g_\pm \in \Spin(8)$ such that 
\[      t(g(v_1), g_+(v_2), g_-(v_3)) = t(v_1,v_2,v_3)  \] 
for all $v_1 \in V_8, v_2 \in S^+_8,$ and $v_3 \in S^-_8$.
Moreover, the maps 
\[   \alpha_\pm \maps g \to g_\pm  \] 
are outer automorphisms of $\Spin(8)$.  In fact ${\rm Out(\Spin(8))}$
is the permutation group on 3 letters, and there exist outer
automorphisms that have the effect of permuting the vector, left-handed
spinor, and right-handed spinor representations any way one likes;
$\alpha_+$ and $\alpha_-$ are among these.   

In general, outer automorphisms of simple Lie groups come from
symmetries of their Dynkin diagrams.  Of all the simple Lie groups,
$\Spin(8)$ has the most symmetrical Dynkin diagram!  It looks like this:

\medskip 
\centerline{\epsfysize=1.0in\epsfbox{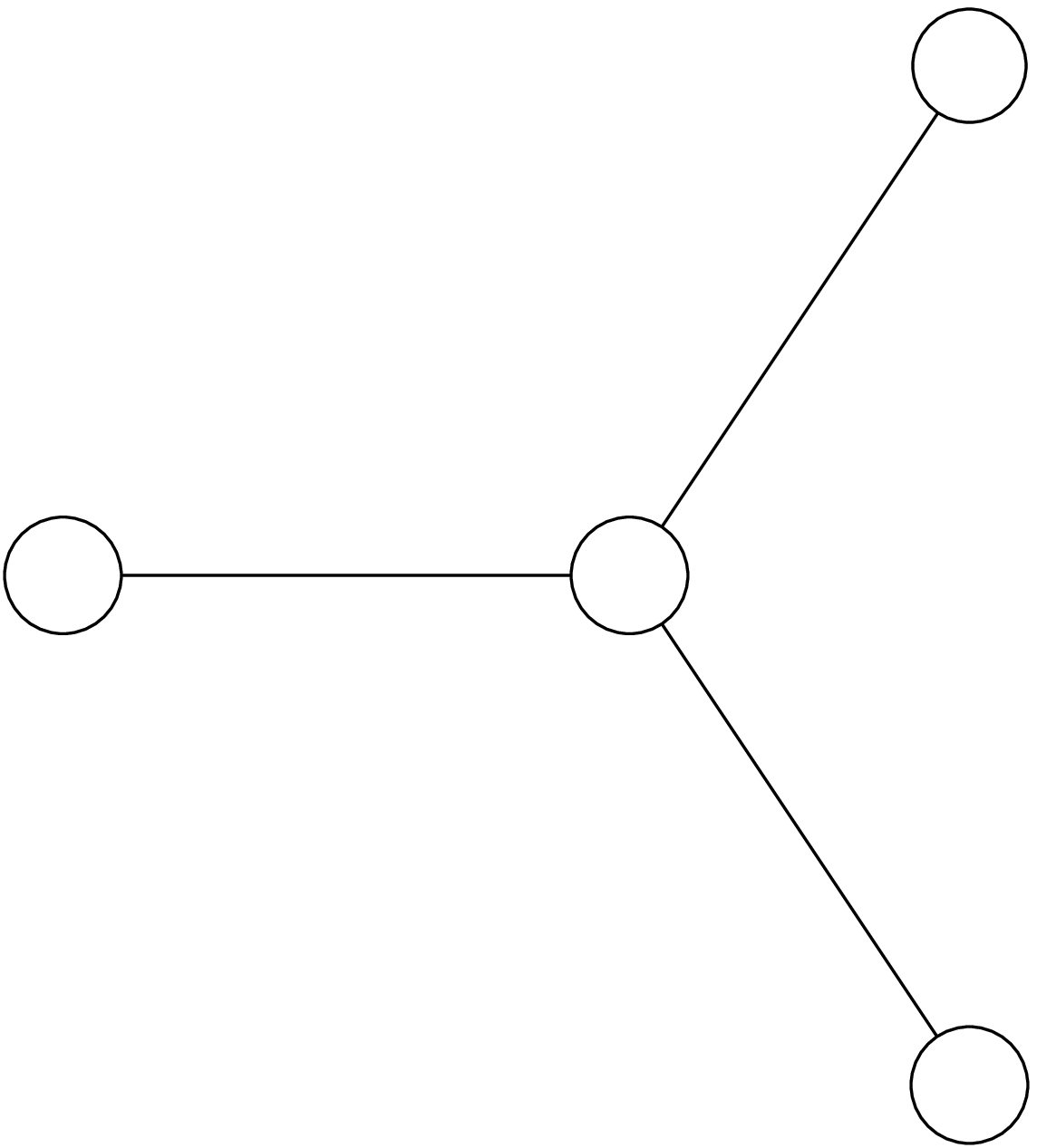}}   
\label{triality.figure}   
\medskip
 
\noindent   
Here the three outer nodes correspond to the vector, left-handed spinor
and right-handed spinor representations of $\Spin(8)$, while the central
node corresponds to the adjoint representation --- that is, the
representation of $\Spin(8)$ on its own Lie algebra, better known as
$\so(8)$.  The outer automorphisms corresponding to the symmetries of
this diagram were discovered in 1925 by Cartan \cite{Cartan3}, who
called these symmetries {\bf triality}.  The more general notion of
`triality' we have been discussing here came later, and is apparently
due to Adams \cite{Adams}.
 
The construction of division algebras from trialities has tantalizing 
links to physics.   In the Standard Model of particle physics, all 
particles other than the Higgs boson transform either as vectors or 
spinors.  The vector particles are also called `gauge bosons', and they
serve to carry the {\it forces} in the Standard Model.  The spinor
particles  are also called `fermions', and they correspond to the basic
forms of {\it matter}: quarks and leptons.   The interaction between
matter and the forces is described by a trilinear map involving two
spinors and one vector.  This map is often drawn as a Feynman diagram:  

\medskip
\centerline{\epsfysize=1.0in\epsfbox{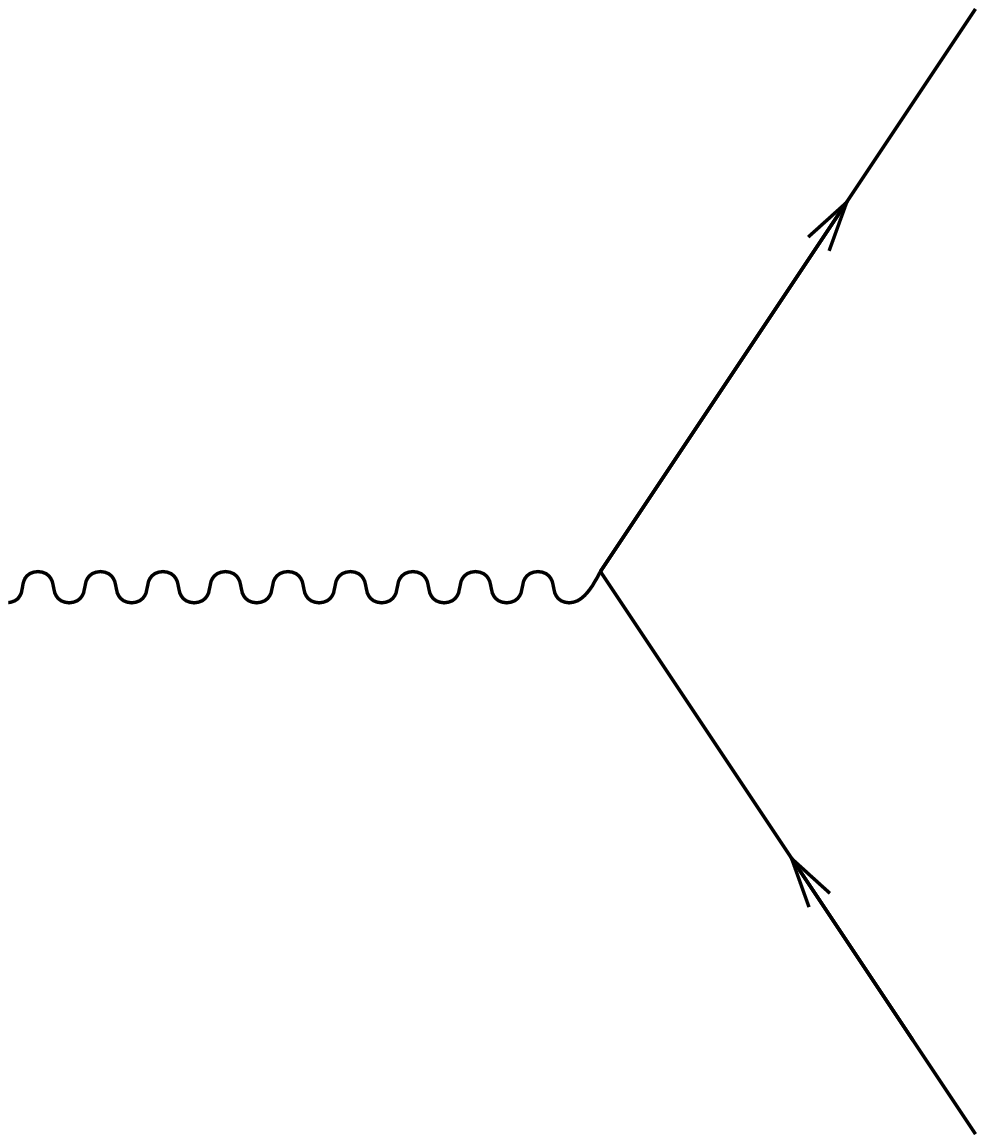}}   
\label{feynman}   
\medskip 
\noindent  
where the straight lines denote spinors and the wiggly one denotes a 
vector.  The most familiar example is the process whereby an electron 
emits or absorbs a photon.   
 
It is fascinating that the same sort of mathematics can be used both to
construct the normed division algebras and to describe the interaction
between matter and forces.  Could this be important for physics?  One
{\it prima facie} problem with this speculation is that physics uses
spinors associated to Lorentz groups rather than rotation groups, due to
the fact that spacetime has a Lorentzian rather than Euclidean metric.
However, in Section \ref{lorentz} we describe a way around this problem.
Just as octonions give the spinor representations of $\Spin(8)$, pairs
of octonions give the spinor representations of $\Spin(9,1)$.  This is
one reason so many theories of physics work best when spacetime is
10-dimensional!  Examples include superstring theory \cite{Deligne,GSW},
supersymmetric gauge theories \cite{Evans,KT,Schray}, and Geoffrey
Dixon's extension of the Standard Model based on the algebra $\C \tensor
\H \tensor \O$, in which the 3 forces arise naturally from the three
factors in this tensor product \cite{Dixon}.
 
\section{Octonionic Projective Geometry}   \label{proj} 
 
Projective geometry is a venerable subject that has its origins in the 
study of perspective by Renaissance painters.  As seen by the eye, 
parallel lines --- e.g., train tracks --- appear to meet at a `point at 
infinity'.  When one changes ones viewpoint, distances and angles appear 
to change, but points remain points and lines remain lines.  These facts 
suggest a modification of Euclidean plane geometry, based on a set of 
points, a set of lines, and relation whereby a point `lies on' a line, 
satisfying the following axioms: 
\begin{itemize} 
\item For any two distinct points, there is a unique line on which they 
both lie. 
\item For any two distinct lines, there is a unique point which lies on 
both of them. 
\item There exist four points, no three of which lie on the same line. 
\item There exist four lines, no three of which have the same point lying 
on them.   
\end{itemize} 
A structure satisfying these axioms is called a {\bf projective plane}. 
Part of the charm of this definition is that it is `self-dual': if we 
switch the words `point' and `line' and switch who lies on whom, it 
stays the same.   
 
We have already met one example of a projective plane in Section  
\ref{fano}: the smallest one of all, the Fano plane.  The example 
relevant to perspective is the real projective plane, $\RP^2$.  Here the 
points are lines through the origin in $\R^3$, the lines are planes 
through the origin in $\R^3$, and the relation of `lying on' is taken to 
be inclusion.  Each point $(x,y) \in \R^2$ determines a point in 
$\RP^2$, namely the line in $\R^3$ containing the origin and the point 
$(x,y,-1)$:

\centerline{\epsfysize=2in\epsfbox{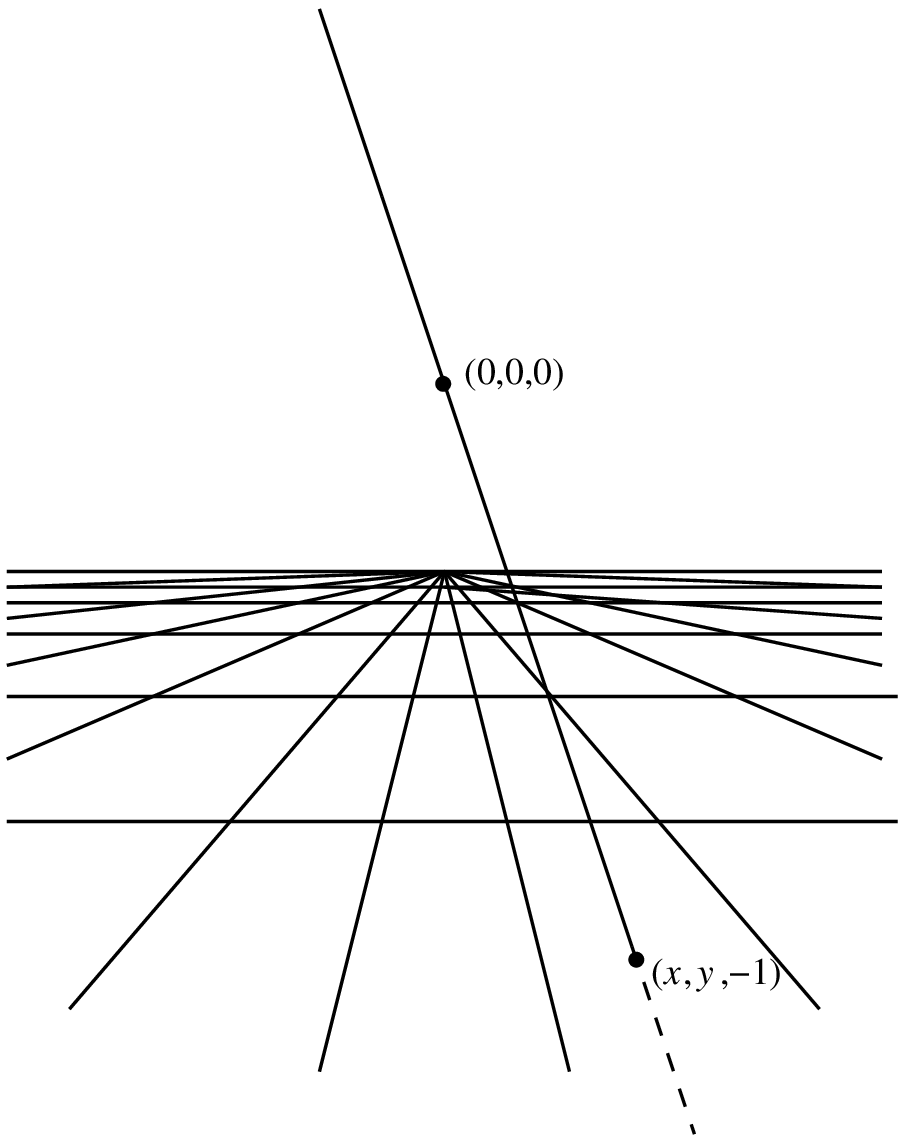}}   
\label{plane}   

\noindent
There are also other points in $\RP^2$, the `points at infinity',
corresponding to lines through the origin in $\R^3$ that do not
intersect the plane $\{z = -1\}$.   For example, any point on the 
horizon in the above picture determines a point at infinity. 
 
Projective geometry is also interesting in higher dimensions. 
One can define a {\bf projective space} by the following axioms: 
\begin{itemize} 
\item For any two distinct points $p,q$, there is a unique line 
$pq$ on which they both lie. 
\item For any line, there are at least three points lying on this line. 
\item If $a,b,c,d$ are distinct points and there is a point lying 
on both $ab$ and $cd$, then there is a point lying on both $ac$ and $bd$. 
\end{itemize} 
Given a projective space and a set $S$ of points in this space, we 
define the {\bf span} of $S$ to be the smallest set $T$ of points
containing $S$ such that if $a$ and $b$ lie in $T$, so do all points
on the line $ab$.  The {\bf dimension}  
of a projective space is defined to be one less than the minimal 
cardinality of a set that spans the whole space.  The reader may enjoy 
showing that a 2-dimensional projective space is the same thing as a 
projective plane \cite{Garner}. 
 
If $\K$ is any field, there is an $n$-dimensional projective space  
called $\KP^n$ where the points are lines through the origin in 
$\K^{n+1}$, the lines are planes through the origin in $\K^{n+1}$, and 
the relation of `lying on' is inclusion.  In fact, this construction 
works even when $\K$ is a mere {\bf skew field}: a ring such that every 
nonzero element has a left and right multiplicative inverse.   We just 
need to be a bit careful about defining lines and planes through the 
origin in $\K^{n+1}$.  To do this, we use the fact that $\K^{n+1}$ is a 
$\K$-bimodule in an obvious way.  We take a line through the origin to 
be any set  
\[   L = \{ \alpha x \; \colon\; \alpha \in \K \}  \] 
where $x \in \K^{n+1}$ is nonzero, and take a plane through the 
origin to be any set 
\[   P = \{ \alpha x + \beta y \; \colon \; \alpha,\beta \in \K \} \] 
where $x,y \in \K^{n+1}$ are elements such that  
$\alpha x + \beta y = 0$ implies $\alpha,\beta = 0$.  
 
Given this example, the question naturally arises whether {\it every} 
projective $n$-space is of the form $\KP^n$ for some skew field $\K$. 
The answer is quite surprising: yes, but only if $n > 2$.  Projective 
planes are more subtle \cite{Stevenson}.  A projective plane comes 
from a skew field if and only if it satisfies an extra axiom, the 
`axiom of Desargues', which goes as follows.  Define a {\bf triangle}  
to be a triple of points that don't all lie on the same line.  Now,  
suppose we have two triangles $xyz$ and $x'y'z'$.  The sides of each  
triangle determine three lines, say $LMN$ and $L'M'N'$.  Sometimes  
the line through $x$ and $x'$, the line through $y$ and $y'$, and  
the line through $z$ and $z'$ will all intersect at the same point: 

\vskip 1em
{\hskip 15em}{\epsfysize=2in\epsfbox{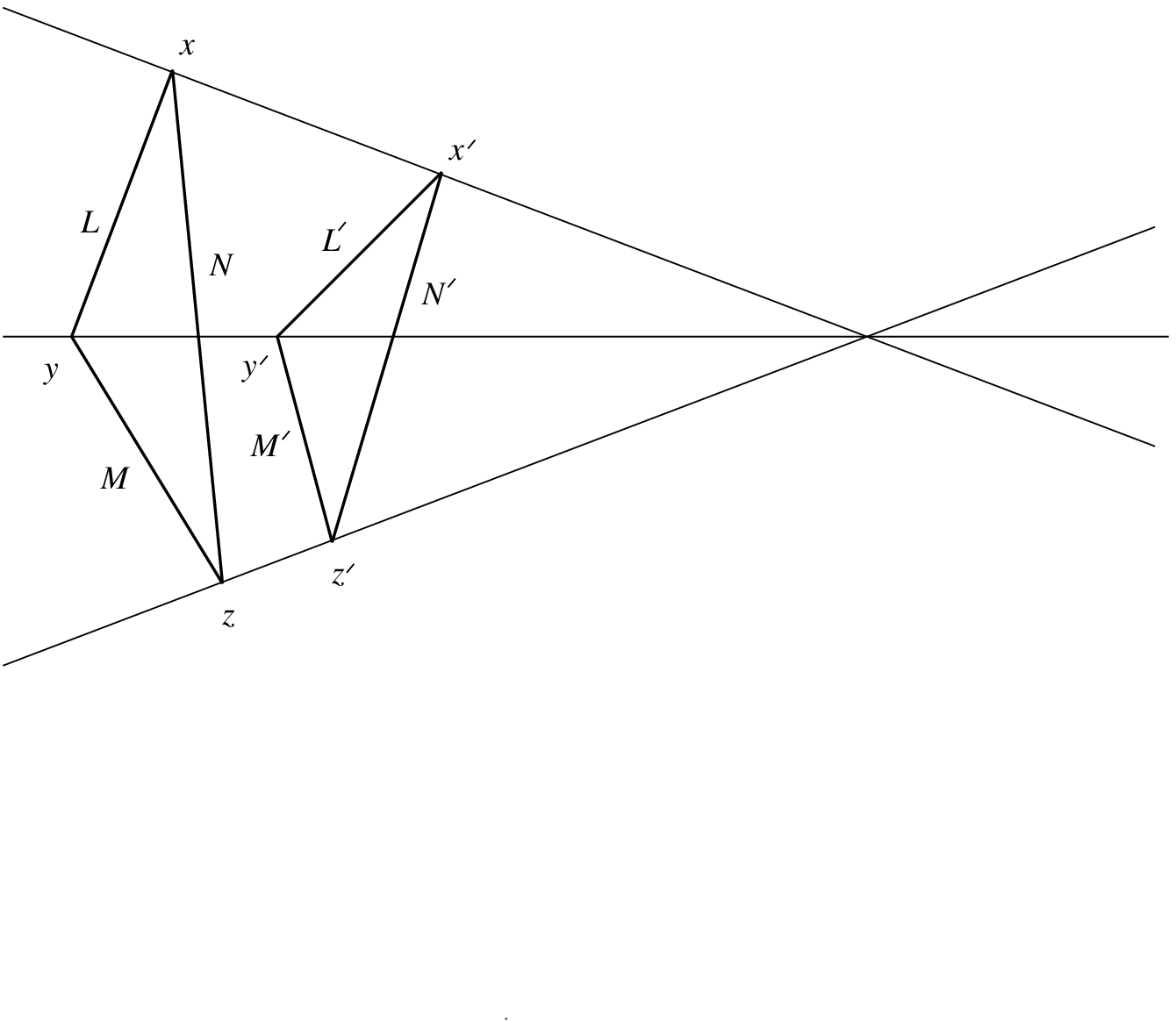}}   
\label{desargues1}   
\vskip -3em 

\noindent 
The {\bf axiom of Desargues} says that whenever this happens, something 
else happens: the intersection of $L$ and $L'$, the intersection of $M$ 
and $M'$, and the intersection of $N$ and $N'$ all lie on the same line: 
 
{\hskip 11em}{\epsfysize=2in\epsfbox{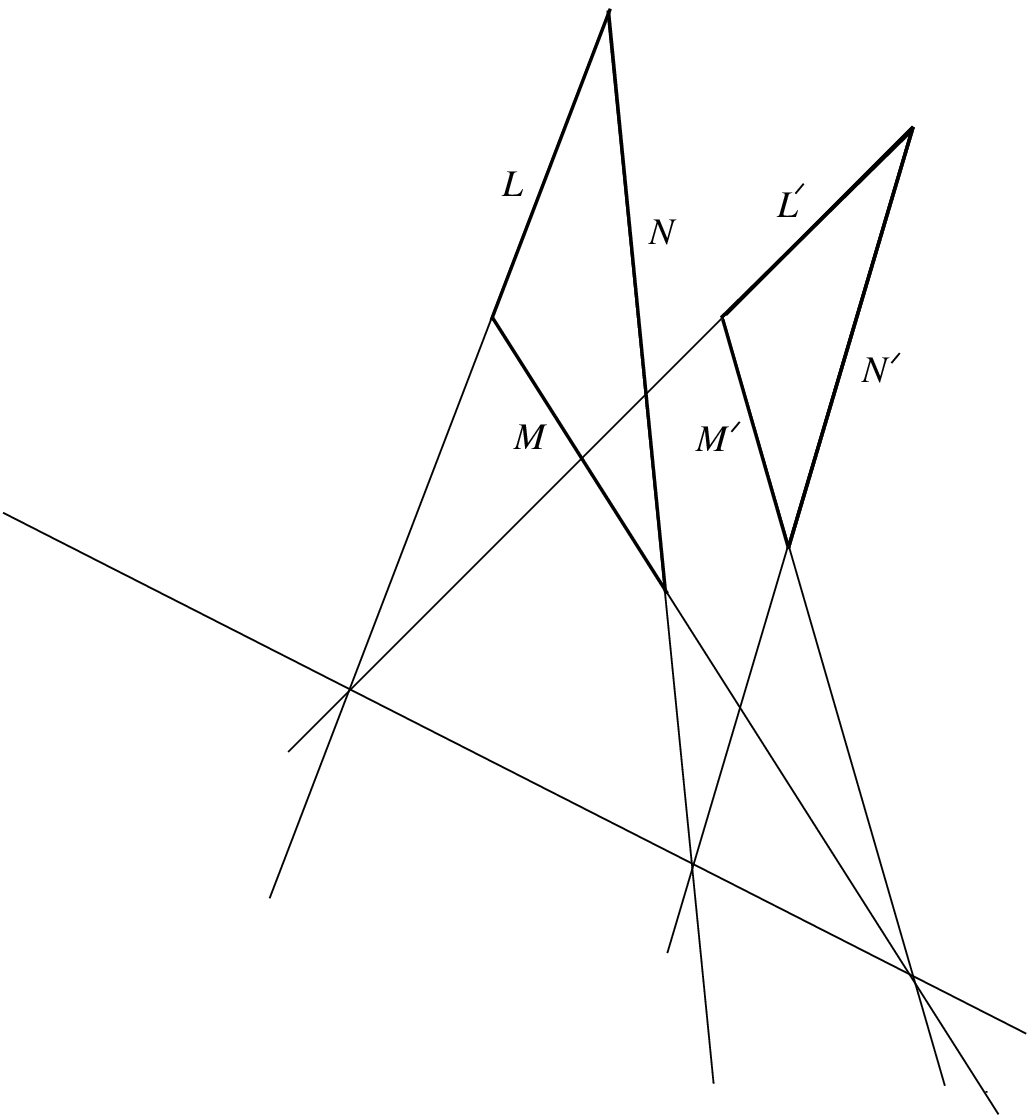}}   
\label{desargues2}

\noindent  
This axiom holds automatically for projective spaces of dimension 3 
or more, but not for projective planes.  A projective plane satisfying  
this axiom is called {\bf Desarguesian}. 
 
The axiom of Desargues is pretty, but what is its connection to skew 
fields?  Suppose we start with a projective plane $P$ and try to 
reconstruct a skew field from it.  We can choose any line $L$, choose 
three distinct points on this line, call them $0, 1$, and $\infty$, and 
set $\K = L - \{\infty\}$.   Copying geometric constructions that work 
when $P = \RP^2$, we can define addition and multiplication of points in 
$\K$.   In general the resulting structure $(\K,+,0,\cdot,1)$ will not 
be a skew field.  Even worse, it will depend in a nontrivial way on the 
choices made.  However, if we assume the axiom of Desargues, these 
problems go away.  We thus obtain a one-to-one correspondence between 
isomorphism classes of skew fields and isomorphism classes of 
Desarguesian projective planes. 
 
Projective geometry was very fashionable in the 1800s, with such  
worthies as Poncelet, Brianchon, Steiner and von Staudt making important 
contributions.  Later it was overshadowed by other forms of geometry.   
However, work on the subject continued, and in 1933 Ruth Moufang 
constructed a remarkable example of a non-Desarguesian projective plane 
using the octonions \cite{Moufang}.  As we shall see, this projective 
plane deserves the name $\OP^2$.   
 
The 1930s also saw the rise of another reason for interest in projective 
geometry: quantum mechanics!  Quantum theory is distressingly different 
from the classical Newtonian physics we have learnt to love.  In 
classical mechanics, observables are described by real-valued functions. 
In quantum mechanics, they are often described by hermitian $n \times 
n$ complex matrices.   In both cases, observables are closed under 
addition and multiplication by real scalars.  However, in quantum 
mechanics, observables do not form an associative algebra.   Still,  
one can raise an observable to a power, and from squaring one 
can construct a commutative but nonassociative product: 
\[     a \circ b = {1\over 2}((a+b)^2 - a^2 - b^2) 
                  = {1\over 2}(ab + ba) . \] 
In 1932, Pascual Jordan attempted to understand this situation better by 
isolating  the bare minimum axioms that an `algebra of observables' 
should satisfy \cite{Jordan}.  He invented the definition of what is now 
called a {\bf formally real Jordan algebra}: a commutative and 
power-associative algebra satisfying 
\[  a_1^2 + \cdots + a_n^2 = 0 \quad \implies \quad a_1 = \cdots = a_n = 0  \] 
for all $n$.   The last condition gives the algebra a partial 
ordering: if we write $a \le b$ when the element $b - a$ is a sum of 
squares, it says that $a \le b$ and $b \le a$ imply $a = b$.  Though it 
is not obvious, any formally real Jordan algebra satisfies the identity  
\[         a \circ (b \circ a^2) = (a \circ b) \circ a^2 \] 
for all elements $a$ and $b$.   Any commutative algebra satisfying 
this identity is called a {\bf Jordan algebra}.  Jordan algebras are 
automatically power-associative.   
 
In 1934, Jordan published a paper with von Neumann and Wigner  
classifying all formally real Jordan algebras \cite{JNW}.  The 
classification is nice and succinct.  An {\bf ideal} in the Jordan algebra 
$A$ is a subspace $B \subseteq A$ such that $b \in B$ implies $a \circ b 
\in B$ for all $a \in A$.  A Jordan algebra $A$ is {\bf simple} if its 
only ideals are $\{0\}$ and $A$ itself.  Every formally real Jordan 
algebra is a direct sum of simple ones.  The simple formally real Jordan 
algebras consist of 4 infinite families and one exception.   
\begin{enumerate} 
\item The algebra $\h_n(\R)$  
with the product $a \circ b = {1\over 2}(ab + ba)$.   
\item The algebra $\h_n(\C)$  
with the product $a \circ b = {1\over 2}(ab + ba)$.   
\item The algebra $\h_n(\H)$  
with the product $a \circ b = {1\over 2}(ab + ba)$.   
\item The algebra $\R^n \oplus \R$ with the product  
\[  (v,\alpha) \circ (w, \beta) =  
(\alpha w + \beta v, \langle v,w\rangle + \alpha \beta).  \] 
\item The algebra $\h_3(\O)$  
with the product $a \circ b = {1\over 2}(ab + ba)$.   
\end{enumerate} 
Here we say a square matrix with entries in the $\ast$-algebra $A$ is 
{\bf hermitian} if it equals its conjugate transpose, and we let 
$\h_n(A)$ stand for the hermitian $n \times n$ matrices with entries in 
$A$.  Jordan algebras in the fourth family are called {\bf spin 
factors}, while $\h_3(\O)$ is called the {\bf exceptional Jordan 
algebra}.  This classification raises some obvious questions.  Why does 
nature prefer the Jordan algebras $\h_n(\C)$ over all the rest?  Or does 
it?  Could the other Jordan algebras --- even the exceptional one --- 
have some role to play in quantum physics?  Despite much research, these 
questions remain unanswered to this day. 
 
The paper by Jordan, von Neumann and Wigner appears to have been 
uninfluenced by Moufang's discovery of $\OP^2$, but in fact they are 
related.  A {\bf projection} in a formally real Jordan algebra is 
defined to be an element $p$ with $p^2 = p$.  In the familiar case of 
$\h_n(\C)$, these correspond to hermitian matrices with eigenvalues $0$ 
and $1$, so they are used to describe observables that assume only two 
values --- e.g., `true' and `false'.  This suggests treating projections 
in a formally real Jordan algebra as propositions in a kind of `quantum 
logic'.  The partial order helps us do this: given projections $p$ and 
$q$, we say that $p$ `implies' $q$ if $p \le q$.    
 
The relation between Jordan algebras and quantum logic is already 
interesting \cite{Emch}, but the real fun starts when we note 
that projections in $\h_n(\C)$ correspond to subspaces of $\C^n$.   This 
sets up a relationship to projective geometry \cite{Varadarajan}, since 
the projections onto 1-dimensional subspaces correspond to points in 
$\CP^n$, while the projections onto 2-dimensional subspaces correspond 
to lines.  Even better, we can work out the dimension of a subspace $V 
\subseteq \C^n$ from the corresponding projection $p \maps \C^n \to V$ 
using only the partial order on projections: $V$ has dimension $d$ iff 
the longest chain of distinct projections  
\[   0 = p_0 < \cdots < p_i = p  \] 
has length $i = d$.  In fact, we can use this to define the {\bf rank} 
of a projection in any formally real Jordan algebra.  We can then try to 
construct a projective space whose points are the rank-1 projections and 
whose lines are the rank-2 projections, with the relation of `lying on' 
given by the partial order $\le$.   
 
If we try this starting with $\h_n(\R)$, $\h_n(\C)$ or $\h_n(\H)$, we 
succeed when $n \ge 2$, and we obtain the projective spaces $\RP^n$, 
$\CP^n$ and $\HP^n$, respectively.  If we try this starting with the 
spin factor $\R^n \oplus \R$ we succeed when $n \ge 2$, and obtain a 
series of 1-dimensional projective spaces related to Lorentzian 
geometry.  Finally, in 1949 Jordan \cite{Jordan2} discovered that if we 
try this construction starting with the exceptional Jordan algebra, we 
get the projective plane discovered by Moufang: $\OP^2$.  
 
In what follows we describe the octonionic projective plane 
and exceptional Jordan algebra in more detail.  But first let us
consider the octonionic projective line, and the Jordan algebra
$\h_2(\O)$.
 
\subsection{Projective Lines}  \label{OP1} 
 
A one-dimensional projective space is called a {\bf projective line}. 
Projective lines are not very interesting from the viewpoint of 
axiomatic projective geometry, since they have only one line on which 
all the points lie.  Nonetheless, they can be geometrically and 
topologically interesting.  This is especially true of the octonionic 
projective line.  As we shall see, this space has a deep connection to 
Bott periodicity, and also to the Lorentzian geometry of 10-dimensional 
spacetime. 
 
Suppose $\K$ is a normed division algebra.  We have already defined 
$\KP^1$ when $\K$ is associative, but this definition does not work well 
for the octonions: it is wiser to take a detour through Jordan 
algebras.  Let $\h_2(\K)$ be the space of $2 \times 2$ hermitian 
matrices with entries in $\K$.  It is easy to check that this becomes a 
Jordan algebra with the product $a \circ b = {1\over 2}(ab + ba)$.  We 
can try to build a projective space from this Jordan algebra using the 
construction in the previous section.  To see if this 
succeeds, we need to ponder the projections in $\h_2(\K)$.  A little 
calculation shows that besides the trivial projections 0 and 1, they 
are all of the form  
\[   
\left( \begin{array}{c}  x^* \\ y^* \end{array} \right)  
\left( \begin{array}{cc} \! x  &  y \! \end{array} \right)  
=  
\left( \begin{array}{cc}  
                         x^* x  & x^* y    \\  
                         y^* x  & y^* y 
\end{array} \right)  
\] 
where $(x,y) \in \K^2$ has 
\[         \|x\|^2 + \|y\|^2 = 1.   \] 
These nontrivial projections all have rank 1, so they are the points of 
our would--be projective space.  Our would--be projective space has just 
one line, corresponding to the projection 1, and all the points lie on 
this line.  It is easy to check that the axioms for a projective space 
hold.  Since this projective space is 1-dimensional, we have succeeded 
in creating the {\bf projective line over} $\K$.  We call the set of 
points of this projective line $\KP^1$. 
 
Given any nonzero element $(x,y) \in \K^2$, we can normalize it and then 
use the above formula to get a point in $\KP^1$, which we call 
$[(x,y)]$.   This allows us to describe $\KP^1$ in terms 
of lines through the origin, as follows.  Define an equivalence relation 
on nonzero elements of $\K^2$ by 
\[        (x,y) \sim (x',y') \; \iff \; [(x,y)] = [(x',y')]  .\] 
We call an equivalence class for this relation a {\bf line through the 
origin} in $\K^2$.  We can then identify  points in $\KP^1$ with lines 
through the origin in $\K^2$.  
 
Be careful: when $\K$ is the octonions, the line through the 
origin containing $(x,y)$ is not always equal to 
\[    \{(\alpha x, \alpha y)\; \colon \; \alpha \in \K\}.   \] 
This is only true when $\K$ is associative, or when $x$ or $y$ is 
$1$.   Luckily, we have $(x,y) \sim (y^{-1}x,1)$ when $y \ne 0$ and 
$(x,y) \sim (1,x^{-1}y)$ when $x \ne 0$.  Thus in either case we get a 
concrete description of the line through the origin containing $(x,y)$: 
when $x \ne 0$ it equals 
\[    \{(\alpha(y^{-1}x), \alpha)\; \colon \; \alpha \in \K\} ,  \] 
and when $y \ne 0$ it equals 
\[    \{(\alpha, \alpha(x^{-1}y)\; \colon \; \alpha \in \K\} .  \] 
In particular, the line through the origin containing $(x,y)$ is  
always a real vector space isomorphic to $\K$. 
 
We can make $\KP^1$ into a manifold as follows.  By the above 
observations, we can cover it with two coordinate charts: one containing 
all points of the form $[(x,1)]$, the other containing all points of the 
form $[(1,y)]$.   It is easy to check that $[(x,1)] = [(1,y)]$ iff $y = 
x^{-1}$, so the transition function from the first chart to the second 
is the map $x \mapsto x^{-1}$.  Since this transition function and its 
inverse are smooth on the intersection of the two charts, $\KP^1$ 
becomes a smooth manifold.  

When pondering the geometry of projective lines it is handy to
visualize the complex case, since $\CP^1$ is just the familiar
`Riemann sphere'.  In this case, the map   
\[    x \mapsto [(x,1)]    \]  
is given by stereographic projection: 
  
\begin{figure}[h]   
\centerline{\epsfysize=1.5in\epsfbox{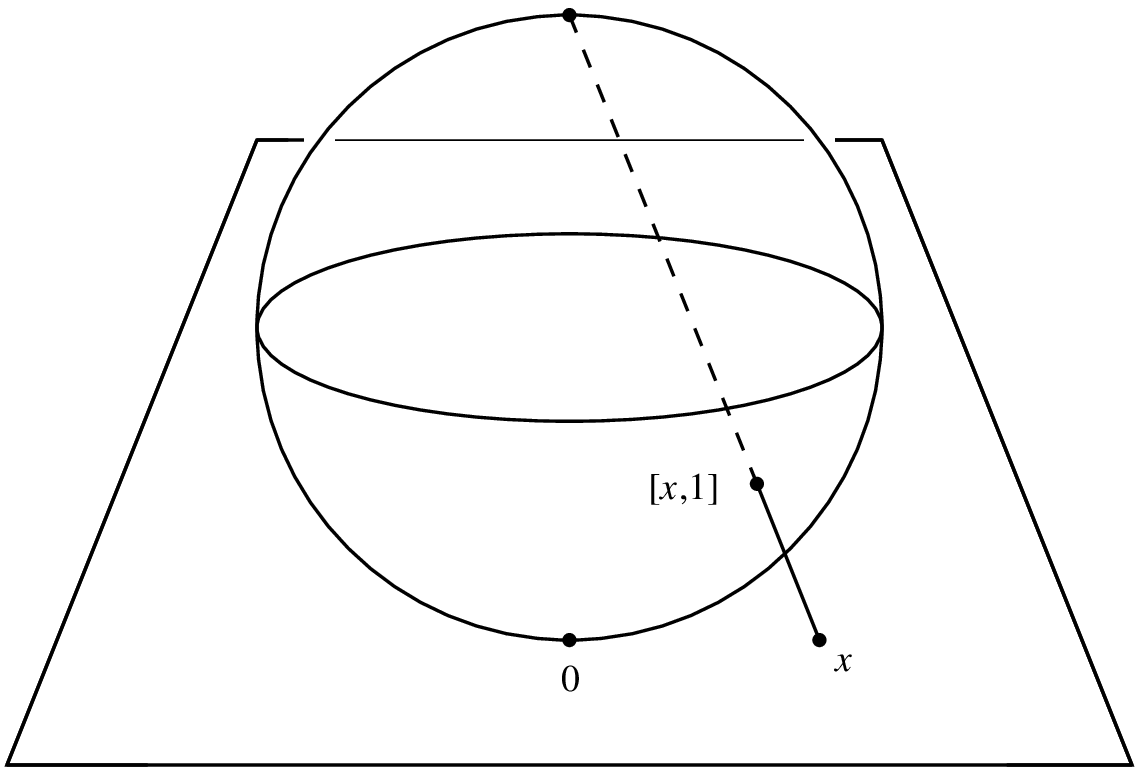}}   
\label{stereo}   
\end{figure}   
  
\noindent  
where we choose the sphere to have diameter 1.   This map from $\C$ to 
$\CP^1$ is one-to-one and almost onto, missing only the point at 
infinity, or `north pole'.  Similarly, the map  
\[    y \mapsto [(1,y)]    \]  
misses only the south pole. Composing the first map with the inverse of 
the second, we get the map  $x \mapsto x^{-1}$, which goes by the name 
of `conformal inversion'.  The southern hemisphere of the Riemann 
sphere consists of all points $[(x,1)]$ with $\|x\| \le 1$, while the 
northern hemisphere consists of all $[(1,y)]$ with $\|y\| \le 1$.  Unit 
complex numbers $x$ give points $[(x,1)] = [(1,x^{-1})]$ on the equator.    
  
All these ideas painlessly generalize to $\KP^1$ for any normed division  
algebra $\K$.  First of all, as a smooth manifold $\KP^1$ is just a   
sphere with dimension equal to that of $\K$:  
\[  
\begin{array}{ccl}       
          \RP^1 &\iso& S^1   \\  
          \CP^1 &\iso& S^2   \\  
          \HP^1 &\iso& S^4   \\  
          \OP^1 &\iso& S^8.  
\end{array}   
\]  
We can think of it as the one-point compactification of $\K$.   The  
`southern hemisphere', `northern hemisphere', and `equator' of $\K$ have 
descriptions exactly like those given above for the complex case.  Also, 
as in the complex case, the maps $x \mapsto [(x,1)]$ and $y \mapsto [(1,y)]$ 
are angle-preserving with respect to the usual Euclidean metric on $\K$ 
and the round metric on the sphere.   

One of the nice things about $\KP^1$ is that it comes equipped with a 
vector bundle whose fiber over the point $[(x,y)]$ is the line 
through the origin corresponding to this point.  This bundle is called 
the {\bf canonical line bundle}, $L_\K$.  Of course, when we are working 
with a particular division algebra, `line' means a copy of this division 
algebra, so if we think of them as real vector bundles, $L_\R, L_\C, 
L_\H$ and $L_\O$ have dimensions 1,2,4, and 8, respectively.   

These bundles play an important role in topology, so it is good to
understand them in a number of ways.  In general, any $k$-dimensional
real vector bundle over $S^n$ can be formed by taking trivial bundles
over the northern and southern hemispheres and gluing  them together
along the equator via a map $f \maps S^{n-1} \to \OO(k)$.  We must
therefore be able to build the canonical line bundles $L_\R, L_\C,
L_\H$ and $L_\O$ using maps  
\[  
\begin{array}{cccl}       
          f_\R \maps &S^0& \to & \OO(1)   \\  
          f_\C \maps &S^1& \to & \OO(2)   \\  
          f_\H \maps &S^3& \to & \OO(4)   \\  
          f_\O \maps &S^7& \to & \OO(8).        
\end{array}   
\]  
What are these maps?  We can describe them all simultaneously.  Suppose  
$\K$ is a normed division algebra of dimension $n$.  In the southern  
hemisphere of $\KP^1$, we can identify any fiber of $L_\K$  with $\K$ by 
mapping the point $(\alpha x, \alpha)$ in the line $[(x,1)]$ to the 
element $\alpha \in \K$.  This trivializes the canonical line bundle 
over the southern hemisphere.  Similarly, we can trivialize this bundle 
over the northern hemisphere by mapping the point $(\beta,\beta y)$ in 
the line $[(1,y)]$ to the element $\beta \in \K$.  If $x \in \K$ has norm 
one, $[(x,1)] = [(1,x^{-1})]$ is a point on the equator, so we get two 
trivializations of the fiber over this point.  These are related as 
follows: if $(\alpha x, \alpha) =  (\beta, \beta x^{-1})$ then $\beta = 
\alpha x$.  The map $\alpha \mapsto \beta$ is thus right multiplication 
by $x$.  In short,   
\[       f_\K \maps S^{n-1} \to \OO(n)  \]  
is just the map sending any norm-one element $x \in \K$ to the operation 
of right multiplication by $x$.   
 
The importance of the map $f_\K$ becomes clearest if we form the 
inductive limit of the groups $\OO(n)$ using the obvious inclusions 
$\OO(n) \hookrightarrow \OO(n+1)$, obtaining a topological group 
called $\OO(\infty)$.  Since $\OO(n)$ is included in $\OO(\infty)$, 
we can think of $f_\K$ as a map from $S^{n-1}$ to $\OO(\infty)$.
Its homotopy class $[f_\K]$ has the following marvelous property,
mentioned in the Introduction:

\vbox{
\begin{itemize}   
\item $[f_\R]$ generates $\pi_0(\OO(\infty)) \iso \Z_2$.     
\item $[f_\C]$ generates $\pi_1(\OO(\infty)) \iso \Z$.   
\item $[f_\H]$ generates $\pi_3(\OO(\infty)) \iso \Z$.   
\item $[f_\O]$ generates $\pi_7(\OO(\infty)) \iso \Z$. 
\end{itemize}  
}

Another nice perspective on the canonical line bundles $L_\K$ comes from 
looking at their unit sphere bundles.  Any fiber of $L_\K$ is naturally 
an inner product space, since it is a line through the origin in $\K^2$. 
If we take the unit sphere in each fiber, we get a bundle of 
$(n-1)$-spheres over $\KP^1$ called the {\bf Hopf bundle}: 
\[               p_\K \maps E_\K \to \KP^1   \] 
The projection $p_\K$ is called the {\bf Hopf map}.  The total space 
$E_\K$ consists of all the unit vectors in $\K^2$, so it is a sphere  
of dimension $2n-1$.  In short, the Hopf bundles look like this: 
\[     
\begin{array}{crccl}       
\K = \R: \qquad &S^0 \hookrightarrow &S^1&    \to &S^1   \\ 
\K = \C: \qquad &S^1 \hookrightarrow &S^3&    \to &S^2   \\ 
\K = \H: \qquad &S^3 \hookrightarrow &S^7&    \to &S^4   \\ 
\K = \O: \qquad &S^7 \hookrightarrow &S^{15}& \to &S^8    
\end{array} 
\] 

We can understand the Hopf maps better by thinking about inverse images 
of points.   The inverse image $p_\K^{-1}(x)$ of any point $x \in S^n$ 
is a $(n-1)$-sphere in $S^{2n-1}$, and the inverse image of any pair of 
distinct points is a pair of linked spheres of this sort.  When $\K = 
\C$  we get linked circles in $S^3$, which form the famous {\bf Hopf link}: 
 
\medskip 
\centerline{\epsfysize=1.0in\epsfbox{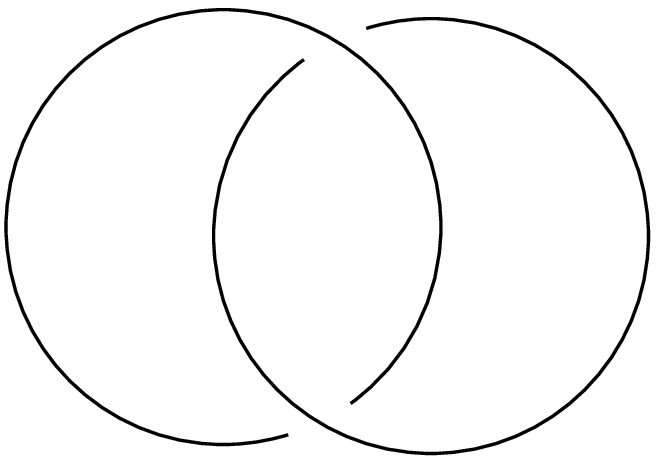}}   
\label{hopf}   
\medskip 
 
\noindent  
When $\K = \O$, we get a pair of linked 7-spheres in $S^{15}$.   
 
To quantify this notion of linking, we can use the `Hopf invariant'.   
Suppose for a moment that $n$ is any natural number greater than one,  
and let $f \maps S^{2n-1} \to S^n$ be any smooth map.  If $\omega$ is 
the normalized volume form on $S^n$, then $f^* \omega$ is a closed 
$n$-form on $S^{2n-1}$.  Since the $n$th cohomology of $S^{2n-1}$ 
vanishes, $f^\ast \omega = d \alpha$ for some $(n-1)$-form $\alpha$. 
We define the {\bf Hopf invariant} of $f$ to be the number 
\[       H(f) = \int_{S^{2n-1}} \alpha \wedge d\alpha  .\] 
This is easily seen to be invariant under smooth homotopies of the map $f$.   

To see how the Hopf invariant is related to linking, we can compute it 
using homology rather than cohomology.  If we take any two regular 
values of $f$, say $x$ and $y$, the inverse images of these points are 
compact oriented $(n-1)$-dimensional submanifolds of $S^{2n-1}$.   We 
can always find an oriented $n$-dimensional submanifold $X \subset 
S^{2n-1}$ that has boundary equal to $f^{-1}(x)$ and that intersects 
$f^{-1}(y)$ transversely.  The dimensions of $X$ and $f^{-1}(y)$ add up 
to $2n-1$, so their intersection number is well-defined.  By the  
duality between homology and cohomology, this number equals the Hopf 
invariant $H(f)$.   This shows that the Hopf invariant is an integer.  
Moreover, it shows that when the Hopf invariant is nonzero, the inverse 
images of $x$ and $y$ are linked. 

Using either of these approaches we can compute the Hopf invariant of 
$p_\C$, $p_\H$ and $p_\O$.  They all turn out to have Hopf invariant 1.  
This implies, for example, that the inverse images of distinct points 
under $p_\O$ are nontrivially linked 7-spheres in $S^{15}$.  It 
also implies that $p_\C$, $p_\H$ and $p_\O$ give nontrivial elements of 
$\pi_{2n-1}(S^n)$ for $n = 2, 4$, and $8$.  In fact, these elements
generate the torsion-free part of $\pi_{2n-1}(S^n)$.

A deep study of the Hopf invariant is one way to prove that any division
algebra must have dimension 1, 2, 4 or 8.  One can show that if there
exists an $n$-dimensional division algebra, then $S^{n-1}$ must be {\bf
parallelizable}: it must admit $n - 1$ pointwise linearly independent
smooth vector fields.  One can also show that for $n > 1$, $S^{n-1}$ is
parallelizable iff there exists a map $f \maps S^{2n-1} \to S^n$ with
$H(f) = 1$ \cite{AH,BM,Kervaire}.  The hard part is showing that a map from
$S^{2n-1}$ to $S^n$ can have Hopf invariant 1 only if $n = 2, 4$, or
$8$.  This was proved by Adams sometime about 1958 \cite{Adams0}.

\subsection{$\OP^1$ and Bott Periodicity} \label{bott}  
 
We already touched upon Bott periodicity when we mentioned that the 
Clifford algebra $\Cliff_{n+8}$ is isomorphic to the algebra of $16  
\times 16$ matrices with entries lying in $\Cliff_n$.  This is but one 
of many related `period-8' phenomena that go by the name of Bott  
periodicity.  The appearance of the number 8 here is no coincidence: all 
these phenomena are related to the octonions!  Since this marvelous   
fact is somewhat under-appreciated, it seems worthwhile to say a bit   
about it.   Here we shall focus on those aspects that are related to   
$\OP^1$ and the canonical octonionic line bundle over this space.    
  
Let us start with K-theory.   This is a way of gaining information about 
a topological space by studying the vector bundles over it.  If the  
space has holes in it, there will be nontrivial vector bundles that  
have `twists' as we go around these holes.  The simplest example is   
the `M\"obius strip' bundle over $S^1$, a 1-dimensional real vector 
bundle which has a $180^\circ$ twist as we go around the circle.  In  
fact, this is just the canonical line bundle $L_\R$.   The canonical 
line bundles $L_\C, L_\H$ and $L_\O$ provide higher-dimensional 
analogues of this example. 
 
K-theory tells us to study the vector bundles over a topological space
$X$ by constructing an abelian group as follows.  First, take the set
consisting of all isomorphism classes of real vector bundles over $X$.
Our ability to take direct sums of vector bundles gives this set an
`addition' operation making it into a commutative monoid.  Next, adjoin
formal `additive inverses' for all the elements of this set, obtaining
an abelian group.  This group is called $KO(X)$, the {\bf real K-theory}
of $X$.  Alternatively we could start with complex vector bundles and 
get a group called $K(X)$, but here we will be interested in real vector
bundles.
 
Any real vector bundle $E$ over $X$ gives an element $[E] \in KO(X)$, and 
these elements generate this group.  If we pick a point in $X$, there is an
obvious homomorphism $\dim \maps KO(X) \to \Z$ sending $[E]$ to the
dimension of the fiber of $E$ at this point.  Since the dimension is a
rather obvious and boring invariant of vector bundles, it is nice to
work with the kernel of this homomorphism, which is called the {\bf
reduced} real K-theory of $X$ and denoted $\widetilde{KO}(X)$.  This is
an invariant of pointed spaces, i.e.\ spaces equipped with a designated
point or {\bf basepoint}.
 
Any sphere becomes a pointed space if we take the north pole as
basepoint.  The reduced real K-theory of the first eight spheres 
looks like this:
\ban       
               \widetilde{KO}(S^1) &\iso& \Z_2  \\    
               \widetilde{KO}(S^2) &\iso& \Z_2  \\    
               \widetilde{KO}(S^3) &\iso& 0  \\    
               \widetilde{KO}(S^4) &\iso& \Z  \\    
               \widetilde{KO}(S^5) &\iso& 0  \\    
               \widetilde{KO}(S^6) &\iso& 0  \\    
               \widetilde{KO}(S^7) &\iso& 0  \\    
               \widetilde{KO}(S^8) &\iso& \Z  
\ean  
where, as one might guess, 
\begin{itemize}    
\item $[L_\R]$ generates $\widetilde{KO}(S^1)$.
\item $[L_\C]$ generates $\widetilde{KO}(S^2)$.
\item $[L_\H]$ generates $\widetilde{KO}(S^4)$.
\item $[L_\O]$ generates $\widetilde{KO}(S^8)$.
\end{itemize}
As mentioned in the previous section, one can build any $k$-dimensional
real vector bundle over $S^n$ using a map $f \maps S^{n-1} \to \OO(k)$.
In fact, isomorphism classes of such bundles are in one-to-one correspondence
with homotopy classes of such maps.  Moreover, two such bundles determine 
the same element of $\widetilde{KO}(X)$ if and only if the corresponding
maps become homotopy equivalent after we compose them with the 
inclusion $\OO(k) \hookrightarrow \OO(\infty)$, where 
$\OO(\infty)$ is the direct limit of the groups $\OO(k)$.   It follows that
\[        \widetilde{KO}(S^n) \iso \pi_{n-1}(\OO(\infty)) .\]
This fact gives us the list of homotopy groups of $\OO(\infty)$ which
appears in the Introduction.   It also means that to prove
Bott periodicity for these homotopy groups:
\[   \pi_{i+8}(\OO(\infty)) \iso \pi_i(\OO(\infty)),   \]   
it suffices to prove Bott periodicity for real K-theory:
\[     \widetilde{KO}(S^{n+8}) \iso \widetilde{KO}(S^n) . \]

Why do we have Bott periodicity in real K-theory?  It turns out
that there is a graded ring $KO$ with 
\[              KO_n = \widetilde{KO}(S^n)  .\]
The product in this ring comes from our ability to take
`smash products' of spheres and also of real vector bundles over these
spheres.  Multiplying by $[L_\O]$ gives an isomorphism
\[
\begin{array}{ccc}
       \widetilde{KO}(S^n) &\to& \widetilde{KO}(S^{n+8})  \\
                       x   &\mapsto& [L_\O] \,x 
\end{array}
\]
In other words, the canonical octonionic line bundle over $\OP^1$
generates Bott periodicity!

There is much more to say about this fact and how it relates to Bott
periodicity for Clifford algebras, but alas, this would take us too far
afield.  We recommend that the interested reader turn to some
introductory texts on K-theory, for example the one by Dale Husemoller
\cite{Husemoller}.  Unfortunately, all the books I know downplay the
role of the octonions.  To spot it, one must bear in mind the relation 
between the octonions and Clifford algebras, discussed in Section
\ref{clifford} above.
 
\subsection{$\OP^1$ and Lorentzian Geometry}  \label{lorentz} 
 
In Section \ref{OP1} we sketched a systematic approach to projective 
lines over the normed division algebras.  The most famous example is the 
Riemann sphere, $\CP^1$.  As emphasized by Penrose \cite{PR}, this space 
has a fascinating connection to Lorentzian geometry --- or in other words, 
special relativity.  All conformal transformations of the Riemann 
sphere come from fractional linear transformations   
\[        z \mapsto {az + b\over cz + d}, \qquad \qquad a,b,c,d \in \C.  \] 
It is easy to see that the group of such transformations is isomorphic 
to $\PSL(2,\C)$: $2 \times 2$ complex matrices with determinant 1,  
modulo scalar multiples of the identity.  Less obviously, it is also 
isomorphic to the Lorentz group $\SO_0(3,1)$: the identity component of 
the group of linear transformations of $\R^4$ that preserve the 
Minkowski metric 
\[     x \cdot y = x_1 y_1 + x_2 y_2 + x_3 y_3 - x_4 y_4 .\] 
This fact has a nice explanation in terms of the `heavenly sphere'.  
Mathematically, this is the 2-sphere consisting of all lines of the form 
$\{\alpha x\}$ where $x \in \R^4$ has $x \cdot x = 0$.  In special 
relativity such lines represent light rays, so the heavenly sphere is 
the sphere on which the stars appear to lie when you look at the night 
sky.  This sphere inherits a conformal structure from the Minkowski 
metric on $\R^4$.  This allows us to identify the heavenly sphere with 
$\CP^1$, and it implies that the Lorentz group acts as conformal 
transformations of $\CP^1$.  In concrete terms, what this means is that 
if you shoot past the earth at nearly the speed of light, the 
constellations in the sky will appear distorted, but all {\it angles} 
will be preserved. 
 
In fact, these results are not special to the complex case: the same
ideas work for the other normed division algebras as well!  The algebras
$\R, \C, \H$ and $\O$ are related to Lorentzian geometry in 3, 4, 6, and
10 dimensions, respectively \cite{MD,MS,MS2,Schray,Sudbery}.  Even
better, a full explanation of this fact brings out new relationships
between the normed division algebras and spinors.  In what follows we
explain how this works for all 4 normed division algebras, with special
attention to the peculiarities of the octonionic case.
 
To set the stage, we first recall the most mysterious of the four 
infinite series of Jordan algebras listed at the beginning of Section 
\ref{proj}: the spin factors.  We described these quite concretely, but 
a more abstract approach brings out their kinship to Clifford algebras.  
Given an $n$-dimensional real inner product space $V$, let the {\bf spin 
factor} $\J(V)$ be the Jordan algebra freely generated by $V$ modulo
relations  
\[     v^2 = \|v\|^2  .\]  
Polarizing and applying the commutative law, we obtain   
\[   v\circ w = \langle v, w \rangle, \]  
so $\J(V)$ is isomorphic to $V \oplus \R$ with the product  
\[  (v,\alpha) \circ (w, \beta) =   
(\alpha w + \beta v, \langle v,w\rangle + \alpha \beta). \] 
 
Though Jordan algebras were invented to study quantum mechanics, the 
spin factors are also deeply related to special relativity.  We can 
think of $\J(V) \iso V \oplus \R$ as {\bf Minkowksi spacetime}, with 
$V$ as space and $\R$ as time.  The reason is that $\J(V)$ is naturally
equipped with a symmetric bilinear form of signature $(n,1)$, the {\bf
Minkowski metric}: 
\[    (v,\alpha)\cdot (w,\beta) = \langle v,w\rangle - \alpha \beta. \] 
The group of linear transformations preserving the Minkowski metric is 
called $\OO(n,1)$, and the identity component of this is called the {\bf 
Lorentz group}, $\SO_0(n,1)$.   We define the {\bf lightcone} ${\rm 
C}(V)$ to consist of all nonzero $x \in \J(V)$ with $x \cdot x = 0$.  A 
1-dimensional subspace of $\J(V)$ spanned by an element of the 
lightcone is called a {\bf light ray}, and the space of all light rays is 
called the {\bf heavenly sphere} ${\rm S}(V)$.  We can identify the 
heavenly sphere with the unit sphere in $V$, since every light ray is 
spanned by an element of the form $(v,1)$ where $v \in V$ has norm one. 
Here is a picture of the lightcone and the heavenly sphere when $V$ is 
2-dimensional: 
 
\medskip 
\centerline{\epsfysize=1.5in\epsfbox{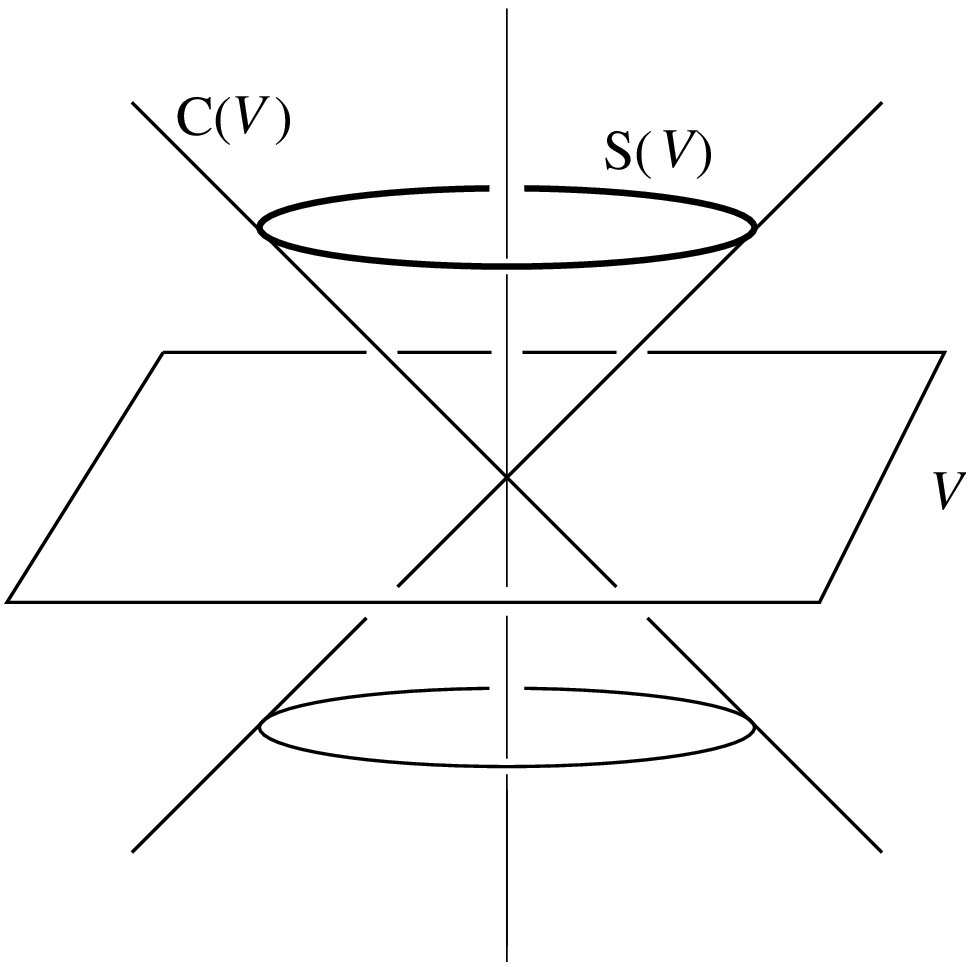}}   
\label{heavenly}   
\medskip 
 
When $V$ is at least 2-dimensional, we can build a projective space from 
the Jordan algebra $\J(V)$.  The result is none other than the heavenly 
sphere!   To see this, note that aside from the elements 0 and 1, all 
projections in $\J(V)$ are of the form $p = \textstyle{1\over 2}(v, 1)$ 
where $v \in V$ has norm one.  These are the points of our projective 
space, but as we have seen, they also correspond to points of the 
heavenly sphere.  Our projective space has just one line, corresponding 
to the projection $1 \in \J(V)$.  We can visualize this line as the 
heavenly sphere itself.  
 
What does all this have to do with normed division algebras?   To answer 
this, let $\K$ be a normed division algebra of dimension $n$.  Then  
the Jordan algebra $h_2(\K)$ is secretly a spin factor!  There is an  
isomorphism  
\[   \phi \maps \h_2(\K) \to J(\K \oplus \R) \iso \K \oplus \R \oplus \R \] 
given by 
\be  \phi \left( \begin{array}{cc}     \alpha + \beta & x     \\  
                                   x^\ast & \alpha - \beta \\  
\end{array} \right) = (x, \beta, \alpha) , 
\qquad \qquad x \in \K, \; \alpha, \beta \in \R .  \label{h2} \ee
Furthermore, the determinant of matrices in $\h_2(\K)$ is well-defined  
even when $\K$ is noncommutative or nonassociative: 
\[  \det \left( \begin{array}{cc}  \alpha + \beta & x     \\  
                                   x^\ast & \alpha - \beta \\  
\end{array} \right) = \alpha^2 - \beta^2 - \|x\|^2 ,  \] 
and clearly we have  
\[      \det(a) = -\phi(a) \cdot \phi(a) \] 
for all $a \in \h_2(\K)$.    
 
These facts have a number of nice consequences.  First of all, since the 
Jordan algebras $\J(\K \oplus \R)$ and $\h_2(\K)$ are isomorphic, so are 
their associated projective spaces.   We have seen that the former space 
is the heavenly sphere ${\rm S}(\K \oplus \R)$, and that the latter is 
$\KP^1$.  It follows that 
\[    \KP^1 \iso {\rm S}(\K \oplus \R) . \] 
This gives another proof of something we already saw in Section 
\ref{OP1}: $\KP^1$ is an $n$-sphere.  But it shows more.   The Lorentz 
group $\SO_0(n+1,1)$ has an obvious action on the heavenly sphere, and 
the usual conformal structure on the sphere is invariant under this 
action.  Using the above isomorphism we can transfer this group action 
and invariant conformal structure to $\KP^1$ in a natural way. 
 
Secondly, it follows that the determinant-preserving linear 
transformations of $\h_2(\K)$ form a group isomorphic to $\OO(n+1,1)$.  
How can we find some transformations of this sort?  If $\K = \R$, this 
is easy: when $g \in \SL(2,\R)$ and $x \in \h_2(\R)$, we again have 
$gxg^* \in \h_2(\R)$, and   
\[     \det(gxg^*) = \det(x). \]  
This gives a homomorphism from $\SL(2,\R)$ to ${\rm O}(2,1)$.  This 
homomorphism is two--to--one, since both $g = 1$ and $g = -1$ act 
trivially, and it maps $\SL(2,\R)$ onto the identity component of ${\rm 
O}(2,1)$.  It follows that $\SL(2,\R)$ is a double cover of  
$\SO_0(2,1)$.  The exact same construction works for $\K = \C$, so 
$\SL(2,\C)$ is a double cover of $\SO_0(3,1)$.    
 
For the other two normed division algebras the above calculation
involving determinants breaks down, and it even becomes tricky to define
the group $\SL(2,\K)$, so we start by working at the Lie algebra level.
We say a $m \times m$ matrix with entries in the normed division algebra
$\K$ is {\bf traceless} if the sum of its diagonal entries is zero.  Any
such traceless matrix acts as a real--linear operator on $\K^m$.  When
$\K$ is commutative and associative, the space of operators coming from
$m \times m$ traceless matrices with entries in $\K$ is closed under
commutators, but otherwise it is not, so we define $\Sl(m,\K)$ to be the
Lie algebra of operators on $\K^m$ {\it generated} by operators of this
form.  This Lie algebra in turn generates a Lie group of real-linear
operators on $\K^m$, which we call $\SL(m,\K)$.  Note that
multiplication in this group is given by composition of real-linear
operators, which is associative even for $\K = \O$.
 
The Lie algebra $\Sl(m,\K)$ comes born with a representation:  
its {\bf fundamental representation} as real-linear operators on $\K^m$,  
given by 
\[   a \maps x \mapsto ax ,           \qquad \qquad x \in \K^m \] 
whenever $a \in \Sl(m,\K)$ actually corresponds to a traceless $m \times 
m$ matrix with entries in $\K$.  Tensoring the fundamental representation  
with its dual, we get a representation of $\Sl(m,\K)$ on the space 
of matrices $\K[m]$, given by 
\[   a \maps x \mapsto ax + xa^*,       \qquad \qquad x \in \K[m]   \] 
whenever $a$ is a traceless matrix with entries in $\K$.  Since $ax + 
xa^*$ is hermitian whenever $x$ is, this representation restricts to a 
representation of $\Sl(m,\K)$ on $\h_m(\K)$.   This in turn can be 
exponentiated to obtain a representation of the group $\SL(m,\K)$ on 
$\h_m(\K)$. 
 
Now let us return to the case $m = 2$.  One can prove that the 
representation of $\SL(2,\K)$ on $\h_2(\K)$ is determinant-preserving 
simply by checking that 
\[    {d \over dt} \det(x + t(ax + xa^*))\, \Bigr|_{t = 0} = 0 \] 
when $x$ lies in $\h_2(\K)$ and $a \in \K[2]$ is traceless.   Here the 
crucial thing is to make sure that the calculation is not spoiled by 
noncommutativity or nonassociativity.  It follows that we have a  
homomorphism 
\[     \alpha_\K \maps \SL(2,\K) \to \SO_0(n+1,1)  \] 
One can check that this is onto, and that its kernel consists of the 
matrices $\pm 1$.  Thus if we define 
\[   \PSL(2,\K) = \SL(2,\K) / \{\pm 1\} ,  \] 
we get isomorphisms 
\[     
\begin{array}{ccl}       
          \PSL(2,\R)& \iso &\SO_0(2,1)   \\ 
          \PSL(2,\C)& \iso &\SO_0(3,1)   \\   
          \PSL(2,\H)& \iso &\SO_0(6,1)   \\    
          \PSL(2,\O)& \iso &\SO_0(9,1) .   
\end{array} 
\] 
Putting this together with our earlier observations, it follows that 
$\PSL(2,\K)$ acts as conformal transformations of $\KP^1$. 
 
We conclude with some words about how all this relates to spinors.  The 
machinery of Clifford algebras and spinors extends effortlessly from the 
case of inner product spaces to vector spaces equipped with an 
indefinite metric.  In particular, the Lorentz group $\SO_0(n+1,1)$ has 
a double cover called $\Spin(n+1,1)$, and this group has certain 
representations called spinor representations.  When $n = 1,2,4$ or $8$, 
we actually have 
\[         \Spin(n+1,1) \iso \SL(2,\K)  \] 
where $\K$ is the normed division algebra of dimension $n$.   The
fundamental representation of $\SL(2,\K)$ on $\K^2$ is the left-handed
spinor representation of $\Spin(n+1,1)$.  Its dual is the right-handed
spinor representation.  Moreover, the interaction between vectors and
spinors that serves as the basis of supersymmetric theories of physics
in spacetimes of dimension 3, 4, 6 and 10 is just the action of
$\h_2(\K)$ on $\K^2$ by matrix multiplication.  In a 
Feynman diagram, this is represented as follows:

\medskip
\hskip 25em \raise  1ex \hbox{$\K^2$}

\centerline{\raise7.5ex \hbox{$\h_2(\K)$} \kern -.5em
\epsfysize=1.0in\epsfbox{feynman.eps}}   

\hskip 25em $\K^2$
\label{feynman2}   

In the case $\K = \C$, Penrose \cite{PR} has described a nice trick for
getting points on the heavenly sphere from spinors.
In fact, it also works for other normed division algebras: 
if $(x,y) \in \K^2$ is nonzero, the hermitian matrix 
\[   
\left( \begin{array}{c}  x \\ y \end{array} \right)  
\left( \begin{array}{cc} \! x^\ast  &  y^\ast \! \end{array} \right)  
=  
\left( \begin{array}{cc}  
                         x x^\ast  & x y^\ast    \\  
                         y x^\ast  & y y^\ast    
\end{array} \right)   
\] 
is nonzero but has determinant zero, so it defines a point on the
heavenly sphere.  If we restrict to spinors of norm one, this trick
reduces to the Hopf map.  Moreover, it clarifies the curious double role
of $\KP^1$ as both the heavenly sphere in special relativity and a space
of propositions in the quantum logic associated to the Jordan algebra
$\h_2(\K)$: any point on the heavenly sphere corresponds to a
proposition specifying the state of a spinor!

\subsection{$\OP^2$ and the Exceptional Jordan Algebra}  \label{OP2} 
 
The octonions are fascinating in themselves, but the magic really starts
when we use them to construct the exceptional Jordan algebra $\h_3(\O)$
and its associated projective space, the octonionic projective plane. 
The symmetry groups of these structures turn out to be exceptional Lie 
groups, and triality gains an eerie pervasive influence over the
proceedings, since an element of $\h_3(\O)$ consists of 3 octonions and
3 real numbers.  Using the relation between normed division algebras and
trialities, we get an isomorphism 
\be 
\begin{array} {ccc} 
  \h_3(\O) & \iso & \R^3 \oplus V_8 \oplus S_8^+ \oplus S_8^-   \\ 
{} & {} & {}  \\
\left( \begin{array}{ccc}  
                         \alpha  &  z^*  & y^*    \\  
                         z       & \beta & x      \\ 
                         y       & x^*   & \gamma   
\end{array} \right) & \mapsto & ((\alpha,\beta,\gamma),x,y,z)   
\end{array} 
\label{jordan.triality} \ee  
where $\alpha,\beta,\gamma \in \R$ and $x,y,z \in \O$.   Examining the
Jordan product in $\h_3(\O)$ then reveals a wonderful fact: while
superficially this product is defined using the $\ast$-algebra structure
of $\O$, it can actually be defined using only the natural maps 
\[      V_8 \times S_8^+ \to S_8^- , \qquad
        V_8 \times S_8^- \to S_8^+ , \qquad
        S_8^+ \times S_8^- \to V_8  
\]  
together with the inner products on these 3 spaces.   All this
information is contained in the normed triality 
\[     t_8 \maps V_8 \times S_8^+ \times S_8^- \to \R ,\]
so any automorphism of this triality gives a automorphism of $\h_3(\O)$.
In Section \ref{triality} we saw that $\Aut(t_8) \iso \Spin(8)$.  With
a little thought, it follows that 
\[    \Spin(8) \subseteq \Aut(\h_3(\O))   .\] 

However, this picture of $\h_3(\O)$ in terms of 8-dimensional Euclidean
geometry is just part of a bigger picture --- a picture set in
10-dimensional Minkowski spacetime!  If we regard $\h_2(\O)$ as sitting
in the lower right-hand corner of $\h_3(\O)$, we get an isomorphism
\be 
\begin{array} {ccc} 
  \h_3(\O) & \iso & \R \oplus \h_2(\O) \oplus \O^2   \\ 
{} & {} & {}  \\
\left( \begin{array}{cc}  
                         \alpha &  \psi^*    \\  
                         \psi   &  a      \\ 
\end{array} \right) & \mapsto & (\alpha,a,\psi)   
\end{array} 
\label{jordan.10d} \ee  
We saw in Section \ref{lorentz} that $a \in \h_2(\O)$ and $\psi \in
\O^2$ can be identified with a vector and a spinor in 10-dimensional
Minkowski spacetime, respectively.   Similarly, $\alpha$ is a scalar.  

This picture gives a representation of $\Spin(9,1)$ as linear
transformations of $\h_3(\O)$.  Unfortunately, most of these
transformations do not preserve the Jordan product on $\h_3(\O)$.  As we
shall see, they only preserve a lesser structure on $\h_3(\O)$: the {\it
determinant}.  However, the transformations coming from the subgroup
$\Spin(9) \subset \Spin(9,1)$ do preserve the Jordan product.  We can
see this as follows.  As a representation of $\Spin(9)$, $\h_2(\O)$
splits into `space' and `time':
\[       \h_2(\O) \iso V_9 \oplus \R  \]
with the two pieces corresponding to the traceless elements of 
$\h_2(\O)$ and the real multiples of the identity, respectively. 
On the other hand, the spinor representation of $\so(9)$ splits
as $S_8^+ \oplus S_8^-$ when we restrict it to $\so(8)$, so we 
have 
\[       \O^2 \iso S_9  .\]
We thus obtain an isomorphism
\be 
\begin{array} {ccc} 
  \h_3(\O) & \iso & \R^2 \oplus V_9 \oplus S_9   \\ 
{} & {} & {}  \\
\left( \begin{array}{cc}  
                         \alpha &  \psi^*    \\  
                         \psi   &  a + \beta     \\ 
\end{array} \right) & \mapsto & ((\alpha,\beta),a,\psi)   
\end{array} 
\label{jordan.9d} \ee  
where $a \in \h_2(\O)$ has vanishing trace and $\beta$ is a real
multiple of the identity.  In these terms, one can easily check that the
Jordan product in $\h_3(\O)$ is built from invariant operations on
scalars, vectors and spinors in 9 dimensions.  It follows that
\[    \Spin(9) \subseteq \Aut(\h_3(\O))   .\] 
For more details on this, see Harvey's book \cite{Harvey}.

This does not exhaust all the symmetries of $\h_3(\O)$, since there are 
other automorphisms coming from the permutation group on 3 letters, 
which acts on $(\alpha,\beta,\gamma) \in \R^3$ and $(x,y,z) \in \O^3$ in 
an obvious way.  Also, any matrix $g \in \OO(3)$ acts by conjugation as 
an automorphism of $\h_3(\O)$; since the entries of $g$ are real, there 
is no problem with nonassociativity here.   The group $\Spin(9)$ is 
36-dimensional, but the full automorphism group $\h_3(\O)$ is much  
bigger: it is 52-dimensional.  As we explain in Section \ref{F4}, it 
goes by the name of $\F_4$. 

However, we can already do something interesting with the automorphisms 
we have: we can use them to diagonalize any element of $\h_3(\O)$.   To 
see this, first note that the rotation group, and thus $\Spin(9)$, acts 
transitively on the unit sphere in $V_9$.  This means we can use an  
automorphism in our $\Spin(9)$ subgroup to bring any element of $\h_3(\O)$ 
to the form 
\[ 
\left( \begin{array}{ccc}  
                         \alpha  &  z^*  & y^*         \\  
                         z       & \beta & x           \\ 
                         y       & x^* & \gamma   \end{array} \right)  
\]  
where $x$ is real.   The next step is to apply an automorphism 
that makes $y$ and $z$ real while leaving $x$ alone.  To do this, note 
that the subgroup of $\Spin(9)$ fixing any nonzero vector in $V_9$ is 
isomorphic to $\Spin(8)$.  When we restrict the representation $S_9$ to 
this subgroup, it splits as $S_8^+ \oplus S_8^-$, and with some work 
\cite{Harvey} one can show that $\Spin(8)$ acts on $S_8^+ \oplus S_8^- 
\iso \O^2$ in such a way that any element $(y,z) \in \O^2$ can be 
carried to an element with both components real.  The final step is to 
take our element of $\h_3(\O)$ with all real entries and use an 
automorphism to diagonalize it.  We can do this by conjugating it with a 
suitable matrix in $\OO(3)$.   
 
To understand $\OP^2$, we need to understand projections in $\h_3(\O)$.  
Here is where our ability to diagonalize matrices in $\h_3(\O)$ via 
automorphisms comes in handy.  Up to automorphism, every projection in 
$\h_3(\O)$ looks like one of these four: 
\[  
p_0 = 
\left( \begin{array}{ccc}  
                         0  &  0  & 0         \\  
                         0  &  0  & 0           \\ 
                         0  &  0  & 0   \end{array} \right) , \; \;
p_1 =  
\left( \begin{array}{ccc}  
                         1  &  0  & 0         \\  
                         0  &  0  & 0           \\ 
                         0  &  0  & 0   \end{array} \right) , \; \;
p_2 =  
\left( \begin{array}{ccc}  
                         1  &  0  & 0         \\  
                         0  &  1  & 0           \\ 
                         0  &  0  & 0   \end{array} \right) , \; \;
p_3 =  
\left( \begin{array}{ccc}  
                         1  &  0  & 0         \\  
                         0  &  1  & 0           \\ 
                         0  &  0  & 1   \end{array} \right) .  \] 
Now, the trace of a matrix in $\h_3(\O)$ is invariant under 
automorphisms, because we can define it using only the Jordan algebra 
structure: 
\[        \tr(a) = {1\over 9} \tr(L_a)    , \qquad \qquad a \in \h_3(\O) \] 
where $L_a$ is left multiplication by $a$.  It follows that the trace of
any projection in $\h_3(\O)$ is 0,1,2, or 3.  Furthermore, the rank of
any projection $p \in \h_3(\O)$ equals its trace.  To see this, first
note that $\tr(p) \ge \rank(p)$, since $p < q$ implies $\tr(p) <
\tr(q)$, and the trace goes up by integer steps.  Thus we only need
show $\tr(p) \le \rank(p)$.  For this it suffices to consider the four
projections shown above, as both trace and rank are invariant under
automorphisms.  Since $p_0 < p_1 < p_2 < p_3$, it is clear that for
these projections we indeed have $\tr(p) \le \rank(p)$.

It follows that the points of the octonionic projective plane are 
projections with trace 1 in $\h_3(\O)$, while the lines are projections 
with trace 2.  A calculation \cite{Harvey} shows that any projection 
with trace 1 has the form 
\be   p =  
\left( \begin{array}{c}  x \\ y \\ z \end{array} \right)  
\left( \begin{array}{ccc} \! x^\ast  &  y^\ast & z^\ast \! \end{array} \right)  
=  
\left( \begin{array}{ccc}  
                         x x^*  & x y^*   & x z^* \\ 
                         y x^*  & y y^*   & y z^* \\ 
                         z x^*  & z y^*   & z z^*  
\end{array} \right)   
\label{projection} 
\ee 
where $(x,y,z) \in \O^3$ has 
\[     (xy)z = x(yz), \qquad \qquad  \|x\|^2 + \|y\|^2 + \|z\|^2 = 1 .\] 
On the other hand, any projection with trace 2 is of the form $1 - p$ 
where $p$ has trace 1.  This sets up a one-to-one correspondence between 
points and lines in the octonionic projective plane.  If we use this 
correspondence to think of both as trace-1 projections, the point $p$ 
lies on the line $p'$ if and only if $p < 1 - p'$.  Of course, $p < 1 - p'$ 
iff $p' < 1 - p$.   The symmetry of this relation means the octonionic 
projective plane is self-dual!  This is also true of the real, complex 
and quaternionic projective planes.  In all cases, the operation that 
switches points and lines corresponds in quantum logic to the `negation'
of propositions \cite{Varadarajan}.

We use $\OP^2$ to stand for the set of points in the octonionic 
projective plane.   Given any nonzero element $(x,y,z) \in \O^3$ with 
$(xy)z = x(yz)$, we can normalize it and then use equation 
(\ref{projection}) to obtain a point $[(x,y,z)] \in \OP^2$.  Copying the 
strategy that worked for $\OP^1$, we can make $\OP^2$ into a smooth manifold  
by covering it with three coordinate charts:
\begin{itemize}
\item one chart containing all points of the form $[(x,y,1)]$, 
\item one chart containing all points of the form $[(x,1,z)]$, 
\item one chart containing all points of the form $[(1,y,z)]$. 
\end{itemize}
Checking that this works is a simple calculation.  The only interesting
part is to make sure that whenever the associative law might appear
necessary, we can either use the alternativity of the octonions or the
fact that only triples with $(xy)z = x(yz)$ give points $[(x,y,z)] \in
\OP^2$.
 
We thus obtain the following picture of the octonionic projective plane.
As a manifold, $\OP^2$ is 16-dimensional.  The lines in $\OP^2$ are
copies of $\OP^1$, and thus 8-spheres.  For any two distinct points in
$\OP^2$, there is a unique line on which they both lie.  For any two
distinct lines, there is a unique point lying on both of them.  There is
a `duality' transformation that maps points to lines and vice versa
while preserving this incidence relation.  In particular, since the
space of all points lying on any given line is a copy of $\OP^1$, so 
is the space of all lines containing a given point!

To dig more deeply into the geometry of $\OP^2$ one needs another
important structure on the exceptional Jordan algebra: the determinant. 
We saw in Section \ref{lorentz} that despite noncommutativity and
nonassociativity, the determinant of a matrix in $\h_2(\O)$ is a
well-defined and useful concept.  The same holds for $\h_3(\O)$!   We
can define the {\bf determinant} of a matrix in $\h_3(\O)$ by
\[ 
\det \left( \begin{array}{ccc}  
                         \alpha  &  z^*  & y^*         \\  
                         z       & \beta & x           \\ 
                         y       & x^* & \gamma   \end{array} \right) = 
\alpha \beta \gamma - (\alpha \|x\|^2 + \beta \|y\|^2 + \gamma \|z\|^2)
+ 2 \Re(xyz) .
\]  
We can express this in terms of the trace and product via
\[   \det(a) = {1\over 3} \tr(a^3) - {1\over 2} \tr(a^2) \tr(a) + 
               {1\over 6} {\tr(a)}^3  .\]
This shows that the determinant is invariant under all automorphisms of
$\h_3(\O)$.  However, the determinant is invariant under an even bigger
group of linear transformations.  As we shall see in Section \ref{E6},
this group is 78-dimensional: it is a noncompact real form of the
exceptional Lie group $\E_6$.  This extra symmetry makes it worth
seeing how much geometry we can do starting with just the determinant
and the vector space structure of $\h_3(\O)$.  

The determinant is a cubic form on $\h_3(\O)$, so there is a unique
symmetric trilinear form 
\[  (\cdot,\cdot,\cdot) \maps \h_3(\O) \times \h_3(\O) \times \h_3(\O)
\to \R \]
such that
\[      (a,a,a) = \det(a) .\]
By dualizing this, we obtain the so-called {\bf cross product}
\[  \times \maps \h_3(\O) \times \h_3(\O) \to \h_3(\O)^*. \]
Explicitly, this is given by
\[     (a \times b)(c) = (a,b,c) .\]
Despite its name, this product is commutative.  

We have already seen that points of $\OP^2$ correspond to trace-1
projections in $\h_3(\O)$.  Freudenthal \cite{Freudenthal} noticed that
these are the same as elements $p \in \h_3(\O)$ with $\tr(p) = 1$ and $p
\times p = 0$.    Even better, we can drop the equation $\tr(p) = 1$ as
long as we promise to work with {\it equivalence classes} of nonzero 
elements satisfying $p \times p = 0$, where two such elements are
equivalent when one is a nonzero real multiple of the other.  Each
such equivalence class $[p]$ corresponds to a unique point of $\OP^2$,
and we get all the points this way.

Given two points $[p]$ and $[q]$, their cross product $p \times q$ is
well-defined up to a nonzero real multiple.  This suggests that we define
a `line' to be an equivalence class of elements $p \times q \in
\h_3(\O)^*$, where again two such elements are deemed equivalent if one
is a nonzero real multiple of the other.   Freudenthal showed that we
get a projective plane isomorphic to $\OP^2$ if we take these as our
definitions of points and lines and decree that the point $[p]$ lies on
the line $[L]$ if and only if $L(p) = 0$.   Note that this equation
makes sense even though $L$ and $p$ are only well-defined up to 
nonzero real multiples.

One consequence of all this is that one can recover the structure of
$\OP^2$ as a projective plane starting from just the determinant on
$\h_3(\O)$: we did not need the Jordan algebra structure!  However, to
get a `duality' map switching points and lines while preserving the
incidence relation, we need a bit more: we need the nondegenerate
pairing 
\[        \langle a,b \rangle =  \tr(ab)    \]
on $\h_3(\O)$.  This sets up an isomorphism 
\[        \h_3(\O) \iso \h_3(\O)^* . \]
This isomorphism turns out to map points to lines, and in fact, it sets
up a one-to-one correspondence between points and lines.  We can use
this correspondence to think of both points and lines in $\OP^2$ as
equivalence classes of elements of $\h_3(\O)$.  In these terms, the
point $p$ lies on the line $\ell$ iff $\langle \ell,p \rangle = 0$.  This
relationship is symmetrical!  It follows that if we switch points and
lines using this correspondence, the incidence relation is preserved.

We thus obtain a very pretty setup for working with $\OP^2$.  If we
use the isomorphism between $\h_3(\O)$ and its dual to reinterpret
the cross product as a map 
\[  \times \maps \h_3(\O) \times \h_3(\O) \to \h_3(\O) ,\]
then not only is the line through distinct points $[p]$ and $[q]$ given
by $[p \times q]$, but also the point in which two distinct lines
$[\ell]$ and $[m]$ meet is given by $[\ell \times m]$.  A 
triple of points $[p], [q]$ and $[r]$ is collinear iff $(p,q,r) = 0$,
and a triple of lines $[\ell]$, $[m]$, $[n]$ meets at a point iff
$(\ell,m,n) = 0$.  In addition, there is a delightful bunch of
identities relating the Jordan product, the determinant, the cross
product and the inner product in $\h_3(\O)$.  

For more on octonionic geometry, the reader is urged to consult the
original papers by Freudenthal
\cite{Freudenthal4,Freudenthal,Freudenthal2,Freudenthal3}, Jacques Tits
\cite{Tits,Tits2} and Tonny Springer
\cite{Springer,Springer2,Springer3}.  The book by Helmut Salzmann {\it
et al} is also good \cite{Salzmann}.  Unfortunately, we must now bid
goodbye to this subject and begin our trip through the exceptional
groups.  However, we shall return to study the symmetries of $\OP^2$ and
the exceptional Jordan algebra in Sections \ref{F4} and \ref{E6}.

\section{Exceptional Lie Algebras}         \label{lie}   
 
On October 18th, 1887, Wilhelm Killing wrote a letter to Friedrich Engel
saying that he had classified the simple Lie algebras.  In the next
three years, this revolutionary work was published in a series of papers
\cite{Killing}.  Besides what we now call the `classical' simple Lie
algebras, he claimed to have found 6 `exceptional' ones --- new
mathematical objects whose existence had never before been suspected. In
fact he gave a rigorous construction of only the smallest of these. In
his 1894 thesis, Cartan \cite{Cartan0} constructed all of them and
noticed that the two 52-dimensional exceptional Lie algebras discovered
by Killing were isomorphic, so that there are really only 5.

The Killing-Cartan classification of simple Lie algebras introduced much
of the technology that is now covered in any introductory course on the
subject, such as roots and weights.  In what follows we shall avoid this
technology, since we wish instead to see the exceptional Lie algebras as
octonionic relatives of the classical ones --- slightly eccentric relatives,
but still having a close connection to {\it geometry}, in particular the
Riemannian geometry of projective planes.  It is also for this reason
that we shall focus on the compact real forms of the simple Lie algebras.

The classical simple Lie algebras can be organized in three infinite
families:
\[
\begin{array}{lcl}
 \so(n) &=&  \{ x \in \R[n] \colon \; x^* = -x, \; \tr(x) = 0\},   \\
 \su(n) &=&  \{ x \in \C[n] \colon \; x^* = -x, \; \tr(x) = 0\},   \\
 \symp(n) &=&  \{ x \in \H[n] \colon \; x^* = -x \}.   
\end{array}
\]
The corresponding Lie groups are
\[
\begin{array}{lcl}
 \SO(n) &=&  \{ x \in \R[n] \colon \; xx^* = 1, \; \det(x) = 1\},   \\
 \SU(n) &=&  \{ x \in \C[n] \colon \; xx^* = 1, \; \det(x) = 1\},   \\
 \Sp(n) &=&  \{ x \in \H[n] \colon \; xx^* = 1 \}.   
\end{array}
\]
These arise naturally as symmetry groups of projective spaces over
$\R$, $\C$, and $\H$, respectively.   More precisely, they arise as
groups of {\bf isometries}: transformations that preserve a specified
Riemannian metric.  Let us sketch how this works, as a warmup for the 
exceptional groups.  

First consider the projective space $\RP^n$.  We can think of this as
the unit sphere in $\R^{n+1}$ with antipodal points $x$ and $-x$
identified.  It thus inherits a Riemannian metric from the sphere,
and the obvious action of the rotation group $\OO(n+1)$ as isometries
of the sphere yields an action of this group as isometries of $\RP^n$
with this metric.  In fact, with this metric, the group of all isometries 
of $\RP^n$ is just 
\[       \Isom(\RP^n) \iso \OO(n+1)/\OO(1)  \]
where $\OO(1) = \{\pm 1\}$ is the subgroup of $\OO(n+1)$ that acts trivially
on $\RP^n$.  The Lie algebra of this isometry group is 
\[       \isom(\RP^n) \iso \so(n+1) .\]

The case of $\CP^n$ is very similar.  We can think of this as the unit
sphere in  $\C^{n+1}$ with points $x$ and $\alpha x$ identified whenever
$\alpha$ is a unit complex number.  It thus inherits a Riemannian metric
from this sphere, and the unitary group $\U(n+1)$ acts as isometries.  
If we consider only the connected component of the isometry group and 
ignore the orientation-reversing isometries coming from complex 
conjugation, we have
\[       \Isom_0(\CP^n) \iso \U(n+1)/\U(1)  \]
where $\U(1)$ is the subgroup that acts trivially on $\CP^n$.  
The Lie algebra of this isometry group is
\[       \isom(\CP^n) \iso \su(n+1)  .\]

The case of $\HP^n$ is subtler, since we must take the noncommutativity
of the quaternions into account.   We can think of $\HP^n$ as the unit
sphere in $\H^{n+1}$ with points $x$ and $\alpha x$ identified whenever
$\alpha$ is a unit quaternion, and as before, $\HP^n$ inherits a
Riemannian metric.  The group $\Sp(n+1)$ acts as isometries of $\HP^n$,
but this action comes from {\it right} multiplication, so 
\[       \Isom(\HP^n) \iso \Sp(n+1)/\{ \pm 1 \},  \]
since not $\Sp(1)$ but only its center $\{\pm 1\}$ acts trivially
on $\HP^n$ by right multiplication.  At the Lie algebra level, this 
gives
\[       \isom(\HP^n) \iso \symp(n+1) . \]

For lovers of the octonions, it is tempting to try a similar
construction starting with $\OP^2$.  While nonassociativity makes things
a bit tricky, we show in Section \ref{F4} that it can in fact be done. 
It turns out that $\Isom(\OP^2)$ is one of the exceptional Lie groups,
namely $\F_4$.  Similarly, the exceptional Lie groups $\E_6$, $\E_7$
and $\E_8$ are in a certain subtle sense the isometry groups of
projective planes over the algebras $\C \tensor \O$, $\H \tensor \O$ and
$\O \tensor \O$.  Together with $\F_4$, these groups can all be 
defined by the so-called `magic square' construction, which makes use of
much of the algebra we have described so far.  We explain three versions
of this construction in Section \ref{magic}.  We then treat the groups
$\E_6, \E_7$ and $\E_8$ individually in the following sections.  But
first, we must introduce $\G_2$: the smallest of the exceptional Lie
groups, and none other than the automorphism group of the octonions.

\subsection{$\G_2$}       \label{G2}   
 
In 1914, \'Elie Cartan noted that the smallest of the exceptional Lie 
groups, $\G_2$, is the automorphism group of the octonions 
\cite{Cartan}.  Its Lie algebra $\g_2$ is therefore $\Der(\O)$, the 
derivations of the octonions.  Let us take these facts as definitions of 
$\G_2$ and its Lie algebra, and work out some of the consequences.   
   
What are automorphisms of the octonions like?  One way to analyze this
involves subalgebras of the octonions.  Any octonion $e_1$ whose square
is $-1$ generates a subalgebra of $\O$ isomorphic to $\C$. If we then
pick any octonion $e_2$ with square equal to $-1$ that anticommutes with
$e_1$, the elements $e_1,e_2$ generate a subalgebra isomorphic to $\H$.
Finally, if we pick any octonion $e_3$ with square equal to $-1$ that
anticommutes with $e_1,e_2,$ and $e_1e_2$, the elements $e_1,e_2,e_3$
generate all of $\O$.  We call such a triple of octonions a {\bf basic
triple}.  Given any basic triple, there exists a unique way to define
$e_4, \dots , e_7$ so that the whole multiplication table in Section
\ref{constructing} holds.  In fact, this follows from the remarks on the
Cayley--Dickson construction at the end of Section \ref{clifford}.
   
It follows that given any two basic triples, there exists a unique   
automorphism of $\O$ mapping the first to the second.  Conversely, it is
obvious that any automorphism maps basic triples to basic triples. 
This gives a nice description of the group $\G_2$, as follows.   
    
Fix a basic triple $e_1,e_2,e_3$.   There is a unique automorphism   
of the octonions mapping this to any other basic triple, say   
$e'_1,e'_2,e'_3$.  Now our description of basic triples so far has   
been purely algebraic, but we can also view them more geometrically   
as follows: a basic triple is any triple of unit imaginary   
octonions (i.e.\ imaginary octonions of norm one) such that each is   
orthogonal to the algebra generated by the other two.  This means that   
our automorphism can map $e_1$ to any point $e'_1$ on the 6-sphere of   
unit imaginary octonions, then map $e_2$ to any point $e'_2$ on the   
5-sphere of unit imaginary octonions that are orthogonal to $e'_1$, and   
then map $e_3$ to any point $e'_3$ on the 3-sphere of unit imaginary   
octonions that are orthogonal to $e'_1,e'_2$ and $e'_1 e'_2$.  It follows   
that      
\[    \dim \G_2 = \dim S^6 + \dim S^5 + \dim S^3 = 14  .\]    
   
The triality description of the octonions in Section \ref{triality}   
gives another picture of $\G_2$.  First, recall that $\Spin(8)$ is the   
automorphism group of the triality $t_8 \maps V_8 \times S^+_8 \times  
S^-_8 \to \R$.  To construct the octonions from this triality we need to  
pick unit vectors in any two of these spaces, so we can think of $\G_2$  
as the subgroup of $\Spin(8)$ fixing unit vectors in $V_8$ and $S^+_8$.  
The subgroup of $\Spin(8)$ fixing a unit vector in $V_8$ is just   
$\Spin(7)$, and when we restrict the representation $S^+_8$ to   
$\Spin(7)$, we get the spinor representation $S_7$.  Thus $\G_2$ is the   
subgroup of $\Spin(7)$ fixing a unit vector in $S_7$.  Since $\Spin(7)$   
acts transitively on the unit sphere $S^7$ in this spinor representation
\cite{Adams}, we have 
\[          \Spin(7)/\G_2 = S^7 .  \]   
It follows that   
\[   \dim \G_2 = \dim (\Spin(7)) - \dim S^7 = 21 - 7 = 14  .\]   
   
This picture becomes a bit more vivid if we remember that after choosing
unit vectors in $V_8$ and $S^+_8$, we can identify both these
representations with the octonions, with both unit vectors corresponding
to $1 \in \O$.  Thus what we are really saying is this: the subgroup of
$\Spin(8)$ that fixes $1$ in the vector representation on $\O$ is
$\Spin(7)$; the subgroup that fixes $1$ in both the vector and
right-handed spinor representations is $\G_2$.  This subgroup also fixes
the element $1$ in the left-handed spinor representation of $\Spin(8)$
on $\O$.
   
Now, using the vector representation of $\Spin(8)$ on $\O$, we    
get homomorphisms   
\[       \G_2 \hookrightarrow \Spin(8) \to \SO(\O)  \]   
where $\SO(\O) \iso \SO(8)$ is the rotation group of the octonions,   
viewed as a real vector space with the inner product $\langle x,y\rangle   
= \Re(x^\ast y)$.   The map from $\Spin(8)$ to $\SO(\O)$ is two-to-one,   
but when we restrict it to $\G_2$ we get a one-to-one map   
\[       \G_2 \hookrightarrow \SO(\O) . \]   
   
At the Lie algebra level, this construction gives an inclusion   
\[      \g_2 \hookrightarrow \so(\O)  \]   
where $\so(\O) \iso \so(8)$ is the Lie algebra of skew-adjoint   
real-linear transformations of the octonions.   Since $\g_2$ is   
14-dimensional and $\so(\O)$ is 28-dimensional, it is nice to see   
exactly where the extra 14 dimensions come from.  In fact, they come   
from two copies of $\Im(\O)$, the 7-dimensional space consisting of all   
imaginary octonions.     
   
More precisely, we have:   
\be  \so(\O) = \g_2 \oplus L_{\Im(\O)} \oplus R_{\Im(\O)} \label{so(O)} \ee   
(a direct sum of vector spaces, not Lie algebras), where $L_{\Im(\O)}$    
is the space of linear transformations of $\O$ given by left multiplication   
by imaginary octonions and $R_{\Im(\O)}$ is the space of linear   
transformations of $\O$ given by right multiplication by imaginary   
octonions \cite{Schafer}.   To see this, we first check that left   
multiplication by an imaginary octonion is skew-adjoint.  Using   
polarization, it suffices to note that   
\[           \langle x,ax \rangle = \Re(x^*(ax)) = \Re((x^*a)x) =   
\Re((a^*x)^*x) = -\Re((ax)^* x) = -\langle ax,x \rangle \]   
for all $a \in \Im(\O)$ and $x \in \O$.  Note that this calculation only   
uses the alternative law, not the associative law, since $x, x^\ast$   
and $a$ all lie in the algebra generated by the two elements $a$ and    
$\Im(x)$.  A similar argument shows that right multiplication by   
an imaginary octonion is skew-adjoint.  It follows that    
$\g_2$, $L_{\Im(\O)}$ and $R_{\Im(\O)}$ all naturally lie in $\so(8)$.  
Next, with some easy calculations we can check that 
\[     L_{\Im(\O)} \cap R_{\Im(\O)} = \{0\} \]
and 
\[    \g_2 \cap (L_{\Im(\O)} + R_{\Im(\O)}) = \{0\} .\]
Using the fact that the dimensions of the 3 pieces adds to 28,
equation (\ref{so(O)}) follows.
   
We have seen that $\G_2$ sits inside $\SO(8)$, but we can do better: it   
actually sits inside $\SO(7)$.  After all, every automorphism of the   
octonions fixes the identity, and thus preserves the space of octonions   
orthogonal to the identity.  This space is just $\Im(\O)$, so   
we have an inclusion   
\[       \G_2 \hookrightarrow \SO(\Im(\O))  \]   
where $\SO(\Im(\O)) \iso \SO(7)$ is the rotation group of the imaginary   
octonions.  At the Lie algebra level this gives an inclusion   
\[        \g_2 \hookrightarrow \so(\Im(\O))  . \]   
   
Since $\g_2$ is 14-dimensional and $\so(\Im(\O))$ is 21-dimensional, it   
is nice to see where the 7 extra dimensions come from.  Examining   
equation (\ref{so(O)}), it is clear that these extra dimensions must   
come from the transformations in $ L_{\Im(\O)} \oplus R_{\Im(\O)}$ that   
annihilate the identity $1 \in \O$.  The transformations that do this are   
precisely those of the form   
\[           \ad_a = L_a - R_a \]   
for $a \in \Im(\O)$.   We thus have    
\be  \so(\Im(\O)) \iso \Der(\O) \oplus \ad_{\Im(\O)}  \label{so(Im(O))}  \ee   
where $\ad_{\Im(\O)}$ is the 7-dimensional space of such transformations.   
   
We may summarize some of the above results as follows:   
   
\begin{thm} \et \label{g2-description} The compact real form   
of the Lie algebra $\g_2$ is given by   
\[     \g_2 = \Der(\O) \subset \so(\Im(\O)) \subset \so(\O)  \]   
and we have    
\ban      \so(\Im(\O)) &=& \Der(\O) \oplus \ad_{\Im(\O)}    \\   
          \so(\O) &=& \Der(\O) \oplus L_{\Im(\O)} \oplus R_{\Im(\O)}    
\ean   
where the Lie brackets in $\so(\Im(\O))$ and $\so(\O)$ are built   
from natural bilinear operations on the summands.     
\end{thm}   
   
As we have seen, $\G_2$ has a 7-dimensional representation $\Im(\O)$.   
In fact, this is the smallest nontrivial representation of $\G_2$,
so it is worth understanding in as many ways as possible.  The space
$\Im(\O)$ has at least three natural structures that are preserved by
the transformations in $\G_2$.  These give more descriptions of  
$\G_2$ as a symmetry group, and they also shed some new light on the   
octonions.  The first two of the structures we describe are analogous to
more familiar ones that exist on the 3-dimensional space of imaginary 
quaternions, $\Im(\H)$.  The third makes explicit use of the   
nonassociativity of the octonions.   
   
First, both $\Im(\H)$ and $\Im(\O)$ are closed under the commutator. In    
the case of $\Im(\H)$, the commutator divided by 2 is the familiar {\bf   
cross product} in 3 dimensions:   
\[            a \times b = {1\over 2}[a,b]  .\]   
We can make the same definition for $\Im(\O)$, obtaining a 7-dimensional   
analog of the cross product.   For both $\Im(\H)$ and $\Im(\O)$ the   
cross product is bilinear and anticommutative.   The cross product   
makes $\Im(\H)$ into a Lie algebra, but not $\Im(\O)$.  For both   
$\Im(\H)$ and $\Im(\O)$, the cross product has two nice    
geometrical properties.  On the one hand, its norm is determined by the    
formula   
\[    \|a \times b\|^2 + \langle a,b\rangle^2 = \|a\|^2 \, \|b\|^2 , \]   
or equivalently,    
\[   \|a \times b\| = | {\sin \theta}| \, \|a\| \, \|b\| ,  \]   
where $\theta$ is the angle between $a$ and $b$.  On the other hand, $a   
\times b$ is orthogonal to $a$ and $b$.  Both these properties follow   
from easy calculations.  For $\Im(\H)$, these two properties are enough   
to determine $x \times y$ up to a sign.  For $\Im(\O)$ they are not ---   
but they become so if we also use the fact that $x \times y$ lies inside   
a copy of $\Im(\H)$ that contains $x$ and $y$.    
   
It is clear that the group of all real-linear transformations of   
$\Im(\H)$ preserving the cross product is just $\SO(3)$, which is also   
the automorphism group of the quaternions.  One can similarly show that   
the group of real-linear transformations of $\Im(\O)$ preserving the   
cross product is exactly $\G_2$.  To see this, start by noting that any   
element of $\G_2$ preserves the cross product on $\Im(\O)$, since the   
cross product is defined using octonion multiplication.  To show that   
conversely any transformation preserving the cross product lies in   
$\G_2$, it suffices to express the multiplication of imaginary octonions    
in terms of their cross product.  Using this identity:   
\[  a \times b = ab + \langle a, b\rangle   , \]   
it actually suffices to express the inner product on $\Im(\O)$ in terms   
of the cross product.  Here the following identity does the job:    
\be  \langle a,b\rangle =   
-\textstyle{1\over 6}\tr(a \times (b \times \cdot\,))    
\label{inner} \ee   
where the right-hand side refers to the trace of the map   
\[       a \times (b \times \cdot\,) \maps \Im(\O) \to \Im(\O)  .\]   
   
Second, both $\Im(\H)$ and $\Im(\O)$ are equipped with a natural 3-form,   
or in other words, an alternating trilinear functional.  This is given   
by   
\[        \phi(x,y,z) = \langle x,yz \rangle  .\]   
In the case of $\Im(\H)$ this is just the usual volume form, and the   
group of real-linear transformations preserving it is $\SL(3, \R)$.   In   
the case of $\Im(\O)$, the real-linear transformations preserving $\phi$   
are exactly those in the group $\G_2$.  A proof of this by Robert  
Bryant can be found in Reese Harvey's book \cite{Harvey}.   The 3-form   
$\phi$ is important in the theory of `Joyce manifolds' \cite{Joyce},  
which are 7-dimensional Riemannian manifolds with holonomy group equal  
to $\G_2$.  
   
Third, both $\Im(\H)$ and $\Im(\O)$ are closed under the associator.    
For $\Im(\H)$ this is boring, since the associator vanishes.  On  the  
other hand, for $\Im(\O)$ the associator is interesting.  In fact, it  
follows from results of Harvey \cite{Harvey} that a real-linear  
transformation $T \maps \Im(\O) \to \Im(\O)$ preserves the associator  
if and only if $\pm T$ lies in $\G_2$.  Thus the symmetry group of the
associator is slightly bigger than $\G_2$: it is $\G_2 \times \Z_2$.  
   
Now we must make an embarrassing admission: these three structures on  
$\Im(\O)$ are all almost the same thing!   Starting with the cross
product   
\[    \times \maps \Im(\O) \times \Im(\O) \to \Im(\O)  \]   
we can recover the usual inner product on $\Im(\O)$ by equation   
(\ref{inner}).  This inner product allows us to dualize the cross
product and obtain a trilinear functional, which is, up to a constant,   
just the 3-form    
\[     \phi \maps \Im(\O) \times \Im(\O) \times \Im(\O) \to \R .\]   
The cross product also determines an orientation on $\Im(\O)$ (we leave   
this as an exercise for the reader).  This allows us to take the Hodge   
dual of $\phi$, obtaining a 4-form $\psi$, i.e.\ an alternating   
tetralinear functional    
\[     \psi \maps \Im(\O) \times \Im(\O) \times \Im(\O) \times \Im(\O)   
\to \R .\]   
Dualizing yet again, this gives a ternary operation which, up   
to a constant multiple, is the associator:   
\[  [\cdot, \cdot, \cdot] \maps \Im(\O) \times \Im(\O) \times \Im(\O)    
\to \Im(\O) .\]   
   
We conclude this section with a handy explicit formula for all the   
derivations of the octonions.    In an associative algebra $A$, any   
element $x$ defines an {\bf inner derivation} $\ad_x \maps A \to A$ by   
\[      \ad_x (a) = [x,a]  \]   
where the bracket denotes the commutator $xa - ax$.   
In a nonassociative algebra, this formula usually does not define   
a derivation.  However, if $A$ is alternative, any pair of elements   
$x,y \in A$ define a derivation $D_{x,y} \maps A \to A$ by   
\be       D_{x,y} a = [[x,y],a] - 3[x,y,a]          \label{D}  \ee   
where $[a,b,x]$ denotes the associator $(ab)x - a(bx)$.   Moreover, 
when $A$ is a normed division algebra, every derivation is a linear 
combination of derivations of this form.  Unfortunately, proving these 
facts seems to require some brutal calculations \cite{Schafer}.  
   
\subsection{$\F_4$}   \label{F4}   

The second smallest of the exceptional Lie groups is the 52-dimensional
group $\F_4$.  The geometric meaning of this group became clear in a
number of nearly simultaneous papers by various mathematicians.  In
1949, Jordan constructed the octonionic projective plane using
projections in $\h_3(\O)$.   One year later, Armand Borel \cite{Borel}
noted that $\F_4$ is the isometry group of a 16-dimensional projective
plane.  In fact, this plane is none other than than $\OP^2$.   Also
in 1950, Claude Chevalley and Richard Schafer \cite{CS} showed that
$\F_4$ is the automorphism group of $\h_3(\O)$.   In 1951, Freudenthal
\cite{Freudenthal4} embarked upon a long series of papers in which he
described not only $\F_4$ but also the other exceptional Lie groups
using octonionic projective geometry.  To survey these developments, one
still cannot do better than to read his classic 1964 paper on Lie groups
and the foundations of geometry \cite{Freudenthal3}.

Let us take Chevalley and Schafer's result as the definition of $\F_4$:
\[     \F_4 = \Aut(\h_3(\O))  .\]
Its Lie algebra is thus
\[   \f_4 = \Der(\h_3(\O)).     \label{f4.1}  \]
As we saw in Section \ref{OP2}, points of $\OP^2$ correspond to trace-1
projections in the exceptional Jordan algebra.  It follows that $\F_4$
acts as transformations of $\OP^2$.  In fact, we can equip $\OP^2$ with
a Riemannian metric for which $\F_4$ is the isometry group. To get a
sense of how this works, let us describe $\OP^2$ as a quotient space of
$\F_4$.   

In Section \ref{OP2} we saw that the exceptional Jordan algebra can
be built using natural operations on the scalar, vector and spinor
representations of $\Spin(9)$.  This implies that $\Spin(9)$ is a
subgroup of $\F_4$.  Equation (\ref{jordan.9d}) makes it clear that
$\Spin(9)$ is precisely the subgroup fixing the element
 \[   \left( \begin{array}{ccc} 1 & 0 & 0 \\  
                                0 & 0 & 0 \\  
                                0 & 0 & 0 \\  
\end{array} \right).  \]  
Since this element is a trace-one projection, it corresponds to a point
of $\OP^2$.  We have already seen that $\F_4$ acts transitively on 
$\OP^2$.  It follows that
\be     \OP^2 \iso \F_4 /\Spin(9)  . \label{F4/Spin(9)} \ee
  
This fact has various nice spinoffs.  First, it gives an easy way to  
compute the dimension of $\F_4$:  
\[    \dim(\F_4) = \dim(\Spin(9)) + \dim(\OP^2) =  
36 + 16 = 52.\]  
Second, since $\F_4$ is compact, we can take any Riemannian metric on $\OP^2$ 
and average it with respect to the action of this group.  The isometry
group of the resulting metric will automatically include $\F_4$ as a
subgroup.  With more work \cite{Besse}, one can show that actually
\[     \F_4 = \Isom(\OP^2) \]
and thus
\[   \f_4 = \isom(\OP^2).     \label{f4.2}  \]

Equation (\ref{F4/Spin(9)}) also implies that the tangent space of our
chosen point in $\OP^2$ is isomorphic to $\f_4/\so(9)$.  But we already
know that this tangent space is just $\O^2$, or in other words, the
spinor representation of $\so(9)$.  We thus have
\be   \f_4 \iso \so(9) \oplus S_9    \label{f4.3}   \ee   
as vector spaces, where $\so(9)$ is a Lie subalgebra.  The bracket in  
$\f_4$ is built from the bracket in $\so(9)$, the action $\so(9) \tensor  
S_9 \to S_9$, and the map $S_9 \tensor S_9 \to \so(9)$  obtained by  
dualizing this action.  We can also rewrite this description
of $\f_4$ in terms of the octonions, as follows:   
\[   \f_4 \iso \so(\O \oplus \R) \oplus \O^2    \label{f4.4}   \]
   
This last formula suggests that we decompose $\f_4$ further using the   
splitting of $\O \oplus \R$ into $\O$ and $\R$.   
It is easily seen by looking at matrices that for all $n,m$ we have    
\be     
\so(n+m) \iso \so(n) \oplus \so(m) \, \oplus\, V_n \tensor V_m.    
\label{so(n+m)} \ee   
Moreover, when we restrict the representation  
$S_9$ to $\so(8)$, it splits as a direct sum $S_8^+ \oplus S_8^-$.  
Using these facts and equation (\ref{f4.3}), we see  
\be      
\f_4 \iso \so(8) \oplus V_8 \oplus S_8^+ \oplus S_8^-     
\label{f4.5}   
\ee   
This formula emphasizes the close relation between $\f_4$ and triality:   
the Lie bracket in $\f_4$ is completely built out of maps involving   
$\so(8)$ and its three 8-dimensional irreducible representations!     
We can rewrite this in a way that brings out the role of the octonions:   
\[
\f_4 \iso \so(\O) \oplus \O^3    
\label{f4.6}    
\]
  
While elegant, none of these descriptions of $\f_4$ gives a  convenient
picture of all the derivations of the exceptional Jordan algebra.   In
fact, there is a nice picture of this sort for $\h_3(\K)$ whenever $\K$
is a normed division algebra.    One way to get a derivation of the
Jordan algebra $\h_3(\K)$ is to take a derivation of $\K$  and let it
act on each entry of the matrices in $\h_3(\K)$.  Another way uses
elements of 
\[    \sa_3(\K) = \{  x \in \K[3] \colon \; x^* = -x,\; \tr(x) = 0 \}  .\]
Given $x \in \sa_3(\K)$, there is a derivation $\ad_x$ of $\h_3(\K)$ given 
by
\[          \ad_x (a) = [x,a]   .\]   
In fact \cite{BS}, every derivation of $\h_3(\K)$ can be uniquely 
expressed as a linear combination of derivations of these two sorts,
so we have 
\be   \Der(\h_3(\K)) \iso \Der(\K) \oplus \sa_3(\K)   \label{der(h3K)} \ee
as vector spaces.   In the case of the octonions, this decomposition 
says that
\[        \f_4 \iso \g_2 \oplus \sa_3(\O) .     \label{f4.7} \]

In equation (\ref{der(h3K)}), the subspace $\Der(\K)$ is always a Lie
subalgebra, but $\sa_3(\K)$ is not unless $\K$ is commutative and
associative --- in which case $\Der(\K)$ vanishes.  Nonetheless, there
is a formula for the brackets in $\Der(\h_3(\K))$ which applies in every
case \cite{OV}.  Given $D,D' \in \Der(\K)$ and $x,y \in \sa_3(\K)$, we
have
\be
\begin{array}{lcl}
            [D,D'] &=& DD' - D'D   \cr
         [D,\ad_x] &=& \ad_{Dx}    \cr
     [\ad_x,\ad_y] &=& \ad_{[x,y]_0} + 
\displaystyle{{1\over 3} \sum_{i,j = 1}^3   D_{x_{ij},y_{ij}}  }
\end{array} 
\label{der(h3K)-bracket}
\ee
where $D$ acts on $x$ componentwise, $[x,y]_0$ is the trace-free   
part of the commutator $[x,y]$, and $D_{x_{ij},y_{ij}}$ is the
derivation of $\K$ defined using equation (\ref{D}).  
   
Summarizing these different descriptions of $\f_4$, we have:   
\begin{thm} \et \label{f4-description}  The compact real form of    
$\f_4$ is given by    
\ban   \f_4 &\iso& \isom(\OP^2)  \\   
            &\iso& \Der(\h_3(\O)) \\   
            &\iso& \Der(\O) \oplus \sa_3(\O) \\   
            &\iso& \so(\O \oplus \R) \oplus \O^2 \\   
            &\iso& \so(\O) \oplus \O^3   
\ean   
where in each case the Lie bracket is built from    
natural bilinear operations on the summands.     
\end{thm}   

\subsection{The Magic Square}  \label{magic}   
   
Around 1956, Boris Rosenfeld \cite{Rosenfeld1} had the remarkable idea
that just as $\F_4$ is the isometry group of the projective plane over
the octonions, the exceptional Lie groups $\E_6$, $\E_7$ and $\E_8$ are
the isometry groups of projective planes over the following three
algebras, respectively:
\begin{itemize}
\item the {\bf bioctonions}, $\C \tensor \O$,
\item the {\bf quateroctonions}, $\H \tensor \O$,
\item the {\bf octooctonions}, $\O \tensor \O$.
\end{itemize}
There is definitely something right about this idea, because one would
expect these projective planes to have dimensions 32, 64, and 128, and
there indeed do exist compact Riemannian manifolds with these
dimensions having $\E_6$, $\E_7$ and $\E_8$ as their isometry groups. 
The problem is that the bioctonions, quateroctonions and and
octooctonions are not division algebras, so it is a nontrivial matter to
define projective planes over them!

The situation is not so bad for the bioctonions: $\h_3(\C \tensor \O)$
is a simple Jordan algebra, though not a formally real one, and one can
use this to define $(\C \tensor \O)\P^2$ in a manner modeled after one
of the constructions of $\OP^2$.  Rosenfeld claimed that a similar
construction worked for the quateroctonions and octooctonions, but this
appears to be false.  Among other problems, $\h_3(\H \tensor \O)$ and
$\h_3(\O \tensor \O)$  do not become Jordan algebras under the product
$a \circ b =  {1\over 2}(ab + ba)$.  Scattered throughout the literature
\cite{Besse,Freudenthal3,Freudenthal5} one can find frustrated comments
about the lack of a really nice construction of $(\H \tensor \O)\P^2$
and $(\O \tensor \O)\P^2$.  One problem is that these spaces do {\it
not} satisfy the usual axioms for a projective plane.   Tits addressed
this problem in his theory of `buildings', which allows one to construct
a geometry having any desired algebraic group as symmetries
\cite{Tits4}.  But alas, it still seems that the quickest way to get our
hands on the quateroctonionic and octooctonionic `projective planes' is
by {\it starting} with the Lie groups $\E_7$ and $\E_8$ and then taking
quotients by suitable subgroups.

In short, more work must be done before we can claim to fully understand
the geometrical meaning of the Lie groups $\E_6, \E_7$ and $\E_8$. 
Luckily, Rosenfeld's ideas can be used to motivate a nice construction
of their Lie algebras.  This goes by the name of the `magic square'.  
Tits \cite{Tits3} and Freudenthal \cite{Freudenthal2} found two very
different versions of this construction in about 1958, but we shall
start by presenting a simplified version published by E.\ B.\ Vinberg
\cite{Vinberg} in 1966.

First consider the projective plane $\KP^2$ when $\K$ is a normed
division algebra $\K$.  The points of this plane are the rank-1
projections in the Jordan algebra $\h_3(\K)$, and this plane admits a
Riemannian metric such that
\[   \isom(\KP^2) \iso \Der(\h_3(\K)). \] 
Moreover, we have seen in equation (\ref{der(h3K)}) that 
\[   \Der(\h_3(\K)) \iso \Der(\K) \oplus \sa_3(\K) . \]
Combined with Rosenfeld's observations, these facts might lead one to 
hope that whenever we have a pair of normed division algebras $\K$ and 
$\K'$, there is a Riemannian manifold $(\K \tensor \K')\P^2$ with
\[   \isom((\K\tensor \K')\P^2) \iso 
\Der(\K) \oplus \Der(\K') \oplus \sa_3(\K \tensor \K')  \]
where for any $\ast$-algebra $A$ we define
\[
\begin{array}{lcl}
    \sh_n(A) & =& \{ x \in A[n] \colon \; x^* = x, \; \tr(x) = 0\}   \\
    \sa_n(A) & =& \{ x \in A[n] \colon \; x^* = -x, \; \tr(x) = 0\}.  
\end{array}
\]

This motivated Vinberg's definition of the {\bf magic square} Lie
algebras:
\[
    \M(\K,\K') = \Der(\K) \oplus \Der(\K') \oplus \sa_3(\K \tensor \K').
\]
Now, when $\K \tensor \K'$ is commutative and associative, $\sa_3(\K
\tensor \K')$ is a Lie algebra with the commutator as its Lie bracket,
but in the really interesting cases it is not.  Thus to make 
$\M(\K,\K')$ into a Lie algebra we must give it a rather subtle bracket.
We have already seen the special case $\K' = \R$ in equation
(\ref{der(h3K)-bracket}).  In general, the Lie bracket in $\M(\K,\K')$
is given as follows:
\begin{enumerate}
\item $\Der(\K)$ and $\Der(\K')$ are commuting Lie subalgebras of $\M(\K,\K')$.
\item The bracket of $D \in \Der(K) \oplus \Der(\K')$ 
with $x \in \sa_3(\K \tensor \K')$ is given by applying 
$D$ to every entry of the matrix $x$, using the natural action of
$\Der(K) \oplus \Der(\K')$ as derivations of $\K \tensor \K'$.
\item Given $X,Y \in \sa_3(\K \tensor \K')$, 
\[      [X,Y] = [X,Y]_0 + 
{1\over 3}\displaystyle{\sum_{i,j = 1}^3   D_{X_{ij},Y_{ij}}  } .\]
Here $[X,Y]_0$ is the traceless part of the $3 \times 3$ matrix $[X,Y]$,
and given $x,y \in \K \tensor \K'$ we define $D_{x,y} \in 
\Der(\K) \oplus \Der(\K')$ in the following way: $D_{x,y}$ is
real-bilinear in $x$ and $y$, and 
\[     D_{a \tensor a',b \tensor b'} = \langle a',b' \rangle D_{a,b} +
                                       \langle a,b \rangle D_{a',b'} \]
where $a,b \in \K$, $a',b' \in \K'$, and $D_{a,b},D_{a',b'}$ are defined
as in equation (\ref{D}).
\end{enumerate}
With this construction we magically obtain the following square of Lie 
algebras:
\vskip 1em   
{\vbox{   
\begin{center}   
\renewcommand{\baselinestretch}{1.3}  
{\small
\begin{tabular}{|c|c|c|c|c|}                    \hline   
           & $\K' = \R$ & $\K' = \C$ & $\K' = \H$ & $\K' = \O$     \\  \hline
$\K = \R$ &  $\so(3)$ & $\su(3)$  & $\symp(3)$  & $\f_4$           \\  \hline
$\K = \C$ &  $\su(3)$ & $\su(3) \oplus \su(3)$ & $\su(6)$ & $\e_6$ \\  \hline
$\K = \H$ &  $\symp(3)$ & $\su(6)$  & $\so(12)$  & $\e_7$          \\  \hline
$\K = \O$ &  $\f_4$   & $\e_6$    & $\e_7$     & $\e_8$            \\  \hline
\end{tabular}} \vskip 1em   
Table 5 --- Magic Square Lie Algebras $\M(\K,\K')$ 
\end{center}   
}   }
\vskip 0.5em   
\noindent
We will mainly be interested in the last row (or column), which is the
one involving the octonions.  In this case we can take the magic square
construction as {\it defining} the Lie algebras $\f_4$, $\e_6$, $\e_7$
and $\e_8$.  This definition turns out to be consistent with our earlier
definition of $\f_4$.

Starting from Vinberg's definition of the magic square Lie algebras, we
can easily recover Tits' original definition.  To do so, we need two
facts.  First, 
\[  
\sa_3(\K \otimes \K') 
\iso \sa_3(\K') \, \oplus \, (\Im(\K) \! \tensor \! \sh_3(\K')).
\]
This is easily seen by direct examination of the relevant matrices.
Second, 
\[    \Der(\h_3(\K)) \iso \Der(\K) \oplus \sa_3(\K)  \]
as vector spaces.  This is just equation (\ref{der(h3K)}).  Starting with
Vinberg's definition and applying these two facts, we obtain
\[
\begin{array}{lcl}
\M(\K,\K') &=& \Der(\K) \oplus \Der(\K') \oplus \sa_3(\K \tensor \K') \\
&\iso& \Der(\K) \oplus \Der(\K') 
\oplus \sa_3(\K') \, \oplus \, (\Im(\K) \tensor \sh_3(\K')) \\
&\iso& 
\Der(\K) \oplus \Der(\h_3(\K')) \, \oplus \, (\Im(\K) \tensor \sh_3(\K')) .
\end{array}
\]
The last line is Tits' definition of the magic square Lie algebras.
Unlike Vinberg's, it is not manifestly symmetrical in $\K$ and $\K'$.
This unhappy feature is somewhat made up for by the fact that $\Der(\K)
\oplus \Der(\h_3(\K'))$ is a nice big Lie subalgebra.  This subalgebra
acts on $\Im(\K) \tensor \sh_3(\K')$ in an obvious way, using the fact
that any derivation of $\K$ maps $\Im(\K)$ to itself, and any derivation
of $\h_3(\K')$ maps $\sh_3(\K')$ to itself.  However, the bracket of two
elements of $(\Im(\K) \tensor \sh_3(\K'))$ is a bit of a mess.

Yet another description of the magic square was recently given by
Barton and Sudbery \cite{BS}.  This one emphasizes the role of 
trialities.  Let $\Tri(\K)$ be the Lie algebra of the group $\Aut(t)$,
where $t$ is the normed triality giving the normed division algebra
$\K$.  From equation (\ref{Aut(t)}) we have
\be
\begin{array}{lcl} 
    \Tri(\R) &\iso& \{0\}    \\
    \Tri(\C) &\iso& \u(1)^2  \\
    \Tri(\H) &\iso& \symp(1)^3 \\
    \Tri(\O) &\iso& \so(8)  .
\end{array} 
\label{tri}
\ee
To express the magic square in terms of these Lie algebras, we need
three facts.  First, it is easy to see that 
\[     \sh_3(\K) \iso \K^3 \oplus \R^2 .\]
Second, Barton and Sudbery show that as vector spaces,
\[    \Der(\h_3(\K)) \iso \Tri(\K) \oplus \K^3 . \]
This follows in a case--by--case way from equation
(\ref{tri}), but they give a unified proof that covers all cases.
Third, they show that as vector spaces,
\[     \Tri(\K) \iso \Der(\K) \oplus \Im(\K)^2 . \]
Now starting with Tits' definition of the magic square,
applying the first two facts, regrouping terms, and applying 
the third fact, we obtain Barton and Sudbery's version of the 
magic square:
\[
\begin{array}{lcl}
\M(\K,\K') 
&\iso& \Der(\K) \oplus \Der(\h_3(\K')) \, \oplus \, 
(\Im(\K) \tensor \sh_3(\K'))   \\
&\iso & \Der(\K) \oplus \Tri(\K') \oplus \K'^3 \oplus 
\Im(\K) \! \tensor \! (\K'^3 \oplus \R^2)  \\
&\iso & \Der(\K) \oplus \Im(\K)^2 \oplus \Tri(\K') \oplus (\K \tensor \K')^3 \\
&\iso &  \Tri(\K) \oplus \Tri(\K') \oplus (\K \tensor \K')^3  .
\end{array}
\]

In the next three sections we use all these different versions of the
magic square to give lots of octonionic descriptions of $\e_6$, $\e_7$
and $\e_8$.  To save space, we usually omit the formulas for the Lie
bracket in these descriptions. However, the patient reader can
reconstruct these with the help of Barton and Sudbery's paper, which 
is packed with useful formulas.

As we continue our tour through the exceptional Lie algebras, we shall
make contact with Adams' work \cite{Adams} constructing
$\f_4,\e_6,\e_7,$ and $\e_8$ with the help of spinors and rotation group
Lie algebras:
\[
\begin{array}{lcl} 
\f_4 &\iso& \so(9) \oplus S_9    \\
\e_6 &\iso& \so(10) \oplus \u(1) \oplus S_{10} \\
\e_7 &\iso& \so(12) \oplus \symp(1) \oplus S_{12}^+ \\
\e_8 &\iso& \so(16) \oplus S_{16}^+ 
\end{array}
\]
as vector spaces.  Note that the numbers 9, 10, 12 and 16 are 8 more
than the dimensions of $\R,\C,\H$ and $\O$.  As usual, this is no
coincidence!  In terms of the octonions, Bott periodicity implies that
\[              S_{n+8} \iso S_n \tensor \O^2  .\]
This gives the following description of spinors in dimensions $\le 16$:
  
\medskip  
{\vbox{   
\begin{center}   
\renewcommand{\baselinestretch}{1.3}  
{\small   
\begin{tabular}{|l|l|}                                  \hline   
$S_1 = \R$           &  $S_9 = \O^2$                       \\ \hline   
$S_2 = \C$           &  $S_{10} = (\C \tensor \O)^2$       \\ \hline  
$S_3 = \H$           &  $S_{11} = (\H \tensor \O)^2$       \\ \hline     
$S_4^\pm = \H$       &  $S_{12}^\pm = (\H \tensor \O)^2$   \\ \hline     
$S_5 = \H^2$         &  $S_{13} = (\H^2 \tensor \O)^2$      \\ \hline     
$S_6 = \C^4$         &  $S_{14} = (\C^4 \tensor \O)^2$      \\ \hline     
$S_7 = \O$           &  $S_{15} = (\O \tensor \O)^2$       \\ \hline     
$S_8^\pm = \O$       &  $S_{16}^\pm = (\O \tensor \O)^2$   \\ \hline     
\end{tabular}}  
\vskip 1em 
\centerline{Table 6 --- Spinor Representations Revisited}  
\end{center}   
}}   
\medskip   

\noindent
Since spinors in dimensions 1,2,4 and 8 are isomorphic to the division
algebras $\R,\C,\H$ and $\O$, spinors in dimensions 8 higher are 
isomorphic to the `planes' $\O^2, (\C \tensor \O)^2, (\H \tensor \O)^2$ 
and $(\O \tensor \O)^2$ --- and are thus closely linked to $\f_4$, $\e_6$, 
$\e_7$ and $\e_8$, thanks to the magic square.

\subsection{$\E_6$}   \label{E6}   

We begin with the 78-dimensional exceptional Lie group $\E_6$.     
As we mentioned in Section \ref{OP2}, there is a nice description of a 
certain noncompact real form of $\E_6$ as the group of collineations
of $\OP^2$, or equivalently, the group of determinant-preserving 
linear transformations of $\h_3(\O)$.  But before going into these,
we consider the magic square constructions of the Lie algebra $\e_6$.
Vinberg's construction gives 
\[
\e_6 = \Der(\O) \oplus \sa_3(\C \tensor \O)  .
\label{e6.1}
\]
Tits' construction, which is asymmetrical, gives
\[
\e_6 \iso \Der(\h_3(\O)) \oplus \sh_3(\O)  
\label{e6.2}
\]
and also
\[
\e_6 \iso \Der(\O) \oplus \Der(\h_3(\C)) \oplus    
(\Im(\O) \!\tensor \! \sh_3(\C))  .
\label{e6.3}
\]
The Barton-Sudbery construction gives
\[
\e_6 \iso \Tri(\O) \oplus \Tri(\C) \oplus (\C \tensor \O)^3 .
\label{e6.4}
\]
We can use any of these to determine the dimension of $\e_6$.   For
example, we have 
\[ \dim(\e_6) = \dim(\Der(\h_3(\O))) + \dim(\sh_3(\O)) = 52 + 26 = 78. \]

Starting from the Barton-Sudbery construction and using the concrete
descriptions of $\Tri(\O)$ and $\Tri(\C)$ from equation
(\ref{tri}), we obtain
\[
\e_6 \iso \so(\O) \oplus \so(\C) \oplus \Im(\C) \oplus (\C \tensor \O)^3 
\label{e6.5}
\]
Using equation (\ref{so(n+m)}), we may rewrite this as
\[
\e_6 \iso \so(\O \oplus \C) \oplus \Im(\C) \oplus (\C \tensor \O)^2
\label{e6.6}
\]
and it turns out that the summand $\so(\O \oplus \C) \oplus \Im(\C)$ is 
actually a Lie subalgebra of $\e_6$.  This result can also be found in 
Adams' book \cite{Adams}, phrased as follows: 
\[
\e_6 \iso \so(10) \oplus \u(1) \oplus S_{10}
\label{e6.7}
\]
In fact, he describes the bracket in $\e_6$ in terms of natural
operations involving $\so(10)$ and its spinor representation $S_{10}$. 
The funny-looking factor of $\u(1)$ comes from the fact that this spinor
representation is complex.  The bracket of an element of $\u(1)$ and
an element of $S_{10}$ is another element of $S_{10}$, defined using
the obvious action of $\u(1)$ on this complex vector space.  

If we define $\E_6$ to be the simply connected group with Lie algebra
$\e_6$, it follows from results of Adams that the subgroup generated by
the Lie subalgebra $\so(10) \oplus \u(1)$ is isomorphic to  $(\Spin(10)
\times \U(1))/\Z_4$.  This lets us define the {\bf bioctonionic
projective plane} by
\[     (\C \tensor \O)\P^2 = \E_6\, / \, ((\Spin(10) \times \U(1))/\Z_4) \]
and conclude that the tangent space at any point of this manifold is
isomorphic to $S_{10} \iso (\C \tensor \O)^2$.  

Since $\E_6$ is compact, we can put an $\E_6$-invariant Riemannian metric on the bioctonionic
projective plane by averaging any metric with respect to the action
of this group.  It turns out \cite{Besse} that the isometry group of this 
metric is exactly $\E_6$, so we have
\[        \E_6 \iso \Isom((\C \tensor \O)\P^2).
\]
It follows that
\[
\e_6 \iso \isom((\C \tensor \O)\P^2)  .
\label{e6.8}
\]

Summarizing, we have 6 octonionic descriptions of $\e_6$:   
\begin{thm} \et \label{e6-description}  The compact real form of    
$\e_6$ is given by    
\ban   
 \e_6 &\iso& \isom((\C \tensor \O)\P^2) \\   
 &\iso& \Der(\O) \oplus \Der(\h_3(\C)) \oplus    
        (\Im(\O) \! \tensor\!  \sh_3(\C)) \\   
 &\iso& \Der(\h_3(\O)) \oplus \sh_3(\O)  \\   
 &\iso& \Der(\O) \oplus \sa_3(\C \tensor \O) \\   
 &\iso& \so(\O \oplus \C) \oplus \Im(\C) \oplus (\C \tensor \O)^2  \\   
 &\iso& \so(\O) \oplus \so(\C) \oplus \Im(\C) \oplus (\C \tensor \O)^3   
\ean   
where in each case the Lie bracket of $\e_6$ is built from    
natural bilinear operations on the summands.     
\end{thm}   

The smallest nontrivial representations of $\E_6$ are 27-dimensional:
in fact it has two inequivalent representations of this dimension, which
are dual to one another.  Now, the exceptional Jordan algebra is also
27-dimensional, and in 1950 this clue led Chevalley and Schafer \cite{CS}
to give a nice description of $\E_6$ as symmetries of this algebra.
These symmetries do not preserve the product, but only the determinant.

More precisely, the group of determinant-preserving linear
transformations of $\h_3(\O)$ turns out to be a noncompact real form of
$\E_6$.  This real form is sometimes called $\E_{6(-26)}$, because its
Killing form has signature $-26$.  To
see this, note that any automorphism of $\h_3(\O)$ preserves the
determinant, so we get an inclusion
\[            \F_4 \hookrightarrow \E_{6(-26)}  .\] 
This means that $\F_4$ is a compact subgroup of $\E_{6(-26)}$.  In fact
it is a maximal compact subgroup, since if there were a larger one, we
could average a Riemannian metric group on $\OP^2$ with respect to this
group and get a metric with an isometry group larger than $\F_4$, but no
such metric exists.  It follows that the Killing form on the Lie algebra
$\e_{6(-26)}$ is negative definite on its 52-dimensional maximal compact
Lie algebra, $\f_4$, and positive definite on the complementary
26-dimensional subspace, giving a signature of $26 - 52 = -26$.

We saw in Section \ref{OP2} that the projective plane structure of
$\OP^2$ can be constructed starting only with the determinant function
on the vector space $\h_3(\O)$.   It follows that $\E_{6(-26)}$
acts as {\bf collineations} on $\OP^2$, that is, line-preserving
transformations.  In fact, the group of collineations of $\OP^2$ is
precisely $\E_{6(-26)}$:
\[   \E_{6(-26)}  \iso     {\rm Coll}(\OP^2). \]
Moreover, just as the group of isometries of $\OP^2$ fixing a specific
point is a copy of $\Spin(9)$, the group of collineations fixing a
specific point is $\Spin(9,1)$.  This fact follows with some work
starting from equation (\ref{jordan.10d}), and it gives us a commutative
square of inclusions:
\[
\begin{array}{ccl} 
   \Spin(9) &\longrightarrow& \Isom(\OP^2) \iso \F_4 \\
   \downarrow & &     \;\;\;\;\;  \downarrow  \\
   \Spin(9,1)    & \longrightarrow & {\rm Coll}(\OP^2) \iso \E_{6(-26)} 
\end{array}
\]
where the groups on the top are maximal compact subgroups of those on
the bottom.  Thus in a very real sense, $\F_4$ is to 9-dimensional
Euclidean geometry as $\E_{6(-26)}$ is to 10-dimensional Lorentzian
geometry.

\subsection{$\E_7$}  \label{E7}   

Next we turn to the 133-dimensional exceptional Lie group $\E_7$.
In 1954, Freudenthal \cite{Freudenthal2} described this group as
the automorphism group of a 56-dimensional octonionic structure now
called a `Freudenthal triple system'.  We sketch this idea
below, but first we give some magic square constructions.
Vinberg's version of the magic square gives
\[
\e_7 = \Der(\H) \oplus \Der(\O) \oplus \sa_3(\H \tensor \O)  .
\label{e7.1}
\]
Tits' version gives
\be
\e_7 \iso \Der(\H) \oplus \Der(\h_3(\O)) \oplus    
(\Im(\H) \!\tensor \! \sh_3(\O))  
\label{e7.2}
\ee
and also
\[
\e_7 \iso \Der(\O) \oplus \Der(\h_3(\H)) \oplus    
(\Im(\O) \!\tensor \! \sh_3(\H))  
\label{e7.3}
\]
The Barton-Sudbery version gives
\be
\e_7 \iso \Tri(\O) \oplus \Tri(\H) \oplus (\H \tensor \O)^3 
\label{e7.4}
\ee
   
Starting from equation (\ref{e7.2}) and using the fact that
$\Der(\H) \iso \Im(\H)$ is 3-dimensional, we obtain the elegant
formula
\[ 
\e_7 \iso \Der(\h_3(\O)) \, \oplus \, \h_3(\O)^3  .
\label{e7.5}
\]
This gives an illuminating way to compute the dimension of $\e_7$:
\[ \dim(\e_7) = \dim(\Der(\h_3(\O))) + 3 \dim(\h_3(\O)) = 52 + 3 \cdot 27 = 
133 .\]
Starting from equation (\ref{e7.4}) and using the concrete 
descriptions of $\Tri(\H)$ and $\Tri(\O)$ from equation
(\ref{tri}), we obtain
\[
\e_7 \iso \so(\O) \oplus \so(\H) \oplus \Im(\H) \oplus (\H \tensor \O)^3 
\label{e7.6}
\]
Using equation (\ref{so(n+m)}), we may rewrite this as
\[
\e_7 \iso \so(\O \oplus \H) \oplus \Im(\H) \oplus (\H \tensor \O)^2.
\label{e7.7}
\]
Though not obvious from what we have done, the direct summand $\so(\O
\oplus \H) \oplus \Im(\H)$ here is really a Lie subalgebra of $\e_7$. In
less octonionic language, this result can also be found in Adams' book
\cite{Adams}:
\[
\e_7 \iso \so(12) \oplus \symp(1) \oplus S_{12}^+ 
\label{e7.8}
\]
He describes the bracket in $\e_7$ in terms of natural operations
involving $\so(12)$ and its spinor representation $S_{12}^+$.   The
funny-looking factor of $\symp(1)$ comes from the fact that this
representation is quaternionic.  The bracket of an element of $\symp(1)$ and
an element of $S_{12}^+$ is the element of $S_{12}^+$ defined using
the natural action of $\symp(1)$ on this space.

If we let $\E_7$ be the simply connected group with Lie 
algebra $\e_7$, it follows from results of Adams \cite{Adams} that the
subgroup generated by the Lie subalgebra $\so(12) \oplus \symp(1)$ is
isomorphic to $(\Spin(12) \times \Sp(1))/\Z_2$.  
This lets us define the {\bf quateroctonionic projective plane} by
\[    (\H \tensor \O)\P^2 = \E_7\, / \,((\Spin(12) \times \Sp(1))/\Z_2) \]
and conclude that the tangent space at any point of this manifold is
isomorphic to $S_{12}^+ \iso (\H \tensor \O)^2$.   We can put 
an $\E_7$-invariant Riemannian metric on this manifold by the technique
of averaging over the group action.  It then turns out \cite{Besse} that
\[        \E_7 \iso \Isom((\H \tensor \O)\P^2)
\]
and thus
\[
\e_7 \iso \isom((\H \tensor \O)\P^2)  .
\label{e7.9}
\]

Summarizing, we have the following 7 octonionic descriptions of $\e_7$:   
   
\begin{thm} \et \label{e7-description}  The compact real form of    
$\e_7$ is given by    
\ban   
\e_7  &\iso& \isom((\H \tensor \O)\P^2)    \\   
&\iso& \Der(\h_3(\O)) \oplus \h_3(\O)^3   \\   
&\iso& \Der(\O) \oplus \Der(\h_3(\H)) \oplus    
       (\Im(\O) \!\tensor \! \sh_3(\H)) \\   
&\iso& \Der(\H) \oplus \Der(\h_3(\O)) \oplus    
        (\Im(\H) \! \tensor\! \sh_3(\O)) \\   
&\iso& \Der(\O) \oplus \Der(\H) \oplus \sa_3(\H \tensor \O) \\   
&\iso& \so(\O \oplus \H) \oplus \Im(\H) \oplus (\H \tensor \O)^2 \\   
&\iso& \so(\O) \oplus \so(\H) \oplus \Im(\H) \oplus (\H \tensor \O)^3  
\ean   
where in each case the Lie bracket of $\e_7$ is built from    
natural bilinear operations on the summands.     
\end{thm}   

Before the magic square was developed, Freudenthal \cite{Freudenthal2} 
used another octonionic construction to study $\E_7$.   The smallest
nontrivial representation of this group is 56-dimensional.  Freudenthal
realized we can define a 56-dimensional space 
\[ F = \{ \left( \begin{array}{cc}  
                         \alpha  & x     \\  
                          y      & \beta \\ 
\end{array} \right) : \;
x,y \in \h_3(\O) , \; \alpha , \beta \in \R \} 
\] 
and equip this space with a symplectic structure 
\[     \omega \maps F \times F \to \R  \]
and trilinear product 
\[     \tau \maps F \times F \times F \to F \]
such that the group of linear transformations preserving both these
structures is a certain noncompact real form of $\E_7$, namely
$\E_{7(-25)}$.  The symplectic structure and trilinear product on
$F$ satisfy some relations, and algebraists have made these into the
definition of a `Freudenthal triple system' \cite{Brown,Faulkner,Meyberg}.
The geometrical significance of this rather complicated sort of structure 
has recently been clarified by some physicists working on string theory.
At the end of the previous section, we
mentioned a relation between 9-dimensional Euclidean geometry and
$\F_4$, and a corresponding relation between 10-dimensional Lorentzian
geometry and $\E_{6(-26)}$.  Murat G\"unaydin \cite{Gunaydin} has 
extended this to a relation between 10-dimensional {\sl conformal} 
geometry and $\E_{7(-25)}$, and in work with Kilian Koepsell 
and Hermann Nikolai \cite{GKN} has explicated how this is connected 
to Freudenthal triple systems.  

\subsection{$\E_8$}   \label{E8}   
   
With 248 dimensions, $\E_8$ is the biggest of the exceptional Lie
groups, and in some ways the most mysterious.  The easiest way to
understand a group is to realize it as as symmetries of a structure one
already understands.  Of all the simple Lie groups, $\E_8$ is the only
one whose smallest nontrivial representation is the adjoint
representation.  This means that in the context of linear algebra,
$\E_8$ is most simply described as the group of symmetries of its own
Lie algebra!  One way out of this vicious circle would be to describe
$\E_8$ as isometries of a Riemannian manifold.  As already mentioned,
$\E_8$ is the isometry group of a 128-dimensional manifold called $(\O
\tensor \O)\P^2$.  But alas, nobody seems to know how to define $(\O
\tensor \O)\P^2$ without first defining $\E_8$.  Thus this group remains
a bit enigmatic.

At present, to get our hands on $\E_8$ we must start with its Lie
algebra.   We can define this using any of the three equivalent magic
square constructions explained in Section \ref{magic}.  Vinberg's
construction gives
\[
\e_8 = \Der(\O) \oplus \Der(\O) \oplus \sa_3(\O \tensor \O)  .
\label{e8.1}
\]
Tits' construction gives
\[
\e_8 \iso \Der(\O) \oplus \Der(\h_3(\O)) \oplus    
(\Im(\O) \!\tensor \! \sh_3(\O)) .
\label{e8.2}
\]
The Barton-Sudbery construction gives
\be
\begin{array}{lcl}
 \e_8 &\iso& \Tri(\O) \oplus \Tri(\O) \oplus (\O \tensor \O)^3 \\
      &\iso& \so(\O) \oplus \so(\O) \oplus (\O \tensor \O)^3 
\end{array}
\label{e8.3}
\ee
We can use any of these to count the dimension of $\e_8$; for example,
the last one gives
\[   \dim \e_8 = 28 + 28 + 3 \cdot 8^2 = 248.\]   

To emphasize the importance of triality, we can rewrite equation 
(\ref{e8.3}) as:
\be
\e_8  \iso \so(8) \oplus \so(8) \oplus (V_8 \tensor V_8) \oplus
             (S_8^+ \tensor S_8^+) \oplus (S_8^- \tensor S_8^-).
\label{e8.4}
\ee
Here the Lie bracket is built from natural maps relating $\so(8)$
and its three 8-dimensional irreducible representations.  In particular,
$\so(8) \oplus \so(8)$ is a Lie subalgebra, and the first copy of
$\so(8)$ acts on the first factor in $V_8 \tensor V_8$, $S_8^+ \tensor
S_8^+$, and $S_8^- \tensor S_8^-$, while the second copy acts on the
second factor in each of these.   The reader should compare this to
the description of $\f_4$ in equation (\ref{f4.5}).

Now, equation (\ref{so(n+m)}) implies that
\[     \so(16) \iso \so(8) \oplus \so(8) \oplus (V_8 \tensor V_8) .\]
Together with equation (\ref{e8.4}), this suggests that $\e_8$ contains
$\so(16)$ as a Lie subalgebra. In fact this is true!  Even better, if we
restrict the right-handed spinor representation of $\so(16)$ to $\so(8)
\oplus \so(8)$, it decomposes as
\[    S^+_{16} \iso (S_8^+ \tensor S_8^+) \oplus (S_8^- \tensor S_8^-),\]
so we obtain
\be   \e_8 \iso \so(16) \oplus S^+_{16}     \label{e8.5} \ee   
or in more octonionic language,
\[   
\e_8 \iso \so(\O \oplus \O) \oplus (\O \otimes \O)^2 
\label{e8.6} \]   
where we use $\so(V)$ to mean the Lie algebra of skew-adjoint real-linear
transformations of the real inner product space $V$.  

The really remarkable thing about equation (\ref{e8.5}) is that the Lie
bracket in $\e_8$ is entirely built from natural maps involving
$\so(16)$ and $S^+_{16}$:  
\[         \so(16) \tensor \so(16) \to \so(16) , \qquad   
           \so(16) \tensor S^+_{16} \to S^+_{16}  , \qquad   
           S^+_{16} \tensor S^+_{16} \to \so(16) .\]   
The first of these is the Lie bracket in $\so(16)$, the second is the
action of $\so(16)$ on its right-handed spinor representation, and the
third is obtained from the second by duality, using the natural inner
product on $\so(16)$ and $S^+_{16}$ to identify these spaces with their
duals.   In fact, this is a very efficient way to {\it define} $\e_8$. 
If we take this approach, we must verify the Jacobi identity:   
\[       [[a,b],c] = [a,[b,c]] - [b,[a,c]]  .\]   
When all three of $a,b,c$ lie in $\so(16)$ this is just the Jacobi   
identity for $\so(16)$.   When two of them lie in $\so(16)$, it boils   
down to fact that spinors indeed form a representation of $\so(16)$.   
Thanks to duality, the same is true when just one lies in $\so(16)$.  It
thus suffices to consider the case when $a,b,c$ all lie in $S_{16}^+$. 
This is the only case that uses anything special about the number 16. 
Unfortunately, at this point a brute-force calculation seems to be
required.  For two approaches that minimize the pain, see the books by
Adams \cite{Adams} and by Green, Schwarz and Witten \cite{GSW}.   It
would be nice to find a more conceptual approach.
   
Starting from $\e_8$, we can define $\E_8$ to be the simply-connected
Lie group with this Lie algebra.  As shown by Adams \cite{Adams}, the 
subgroup of $\E_8$ generated by the Lie subalgebra $\so(16) \subset
\e_8$ is $\Spin(16)/\Z_2$.  This lets us define the {\bf octooctonionic
projective plane} by  
\[            (\O \tensor \O)\P^2 = \E_8\,/\,(\Spin(16)/\Z_2) . \]  
By equation (\ref{e8.5}), the tangent space at any point 
of this manifold is isomorphic to $S_{16}^+ \iso (\O \tensor \O)^2$.
This partially justifies calling it `octooctonionic projective plane',
though it seems not to satisfy the usual axioms for a projective plane.

We can put an $\E_8$-invariant Riemannian metric on the octooctonionic
projective plane by the technique of averaging over the group action.
It then turns out \cite{Besse} that 
\[        \E_8 \iso \Isom((\O \tensor \O)\P^2)   \]
and thus 
\[       \e_8 \iso \isom((\O \tensor \O)\P^2)  . \label{e8.7} \]
   
Summarizing, we have the following octonionic descriptions of 
$\E_8$:
\begin{thm} \et \label{e8-description}  The compact real form of    
$\e_8$ is given by    
\ban   
\e_8 &\iso& \isom((\O \tensor \O)\P^2)    \\  
     &\iso& \Der(\O) \oplus \Der(\h_3(\O)) \oplus    
            (\Im(\O) \!\tensor\! \sh_3(\O)) \\   
     &\iso&  \Der(\O) \oplus \Der(\O) \oplus \sa_3(\O \tensor \O) \\   
     &\iso& \so(\O \oplus \O) \oplus (\O \tensor \O)^2 \\   
     &\iso& \so(\O) \oplus \so(\O) \oplus (\O \tensor \O)^3     
\ean   
where in each case the Lie bracket on $\e_8$ is built from    
natural bilinear operations on the summands.     
\end{thm}   

\section{Conclusions} 

It should be clear by now that besides being a fascinating mathematical
object in their own right, the octonions link together many important
phenomena whose connections would otherwise be completely mysterious.  
Indeed, the full story of these connections is deeper and more elaborate
than I have been able to explain here!   It also includes:
\begin{itemize} 
\item Attempts to set up an octonionic analogue of the theory of 
analytic functions (see \cite{GT} and the references therein). 
\item The role of Jordan pairs, Jordan triple systems and 
Freudenthal triple systems in the construction of exceptional Lie groups 
\cite{Brown,Faulkner,FF,GKN,GT,McCrimmon,Meyberg}.
\item Constructions of the $\E_8$ lattice and Leech lattice using
integral octonions \cite{Coxeter,EG}.
\item Tensor-categorical approaches to normed division algebras
and the invariant of framed trivalent graphs coming from the
quantum group associated to $\G_2$ \cite{Boos,Bremner,Kuperberg,Rost}.
\item Octonionic constructions of vertex operator algebras \cite{FFrenkel}. 
\item Octonionic constructions of the exceptional simple Lie superalgebras 
\cite{Sudbery2}. 
\item Octonionic constructions of symmetric spaces \cite{Besse}.
\item Octonions and the geometry of the `squashed 7-spheres', that is, 
the homogeneous spaces $\Spin(7)/\G_2$, $\Spin(6)/\SU(3)$, and 
$\Spin(5)/\SU(2)$, all of which are diffeomorphic to $S^7$ with its 
usual smooth structure \cite{CD}. 
\item The theory of `Joyce manifolds', that is, 7-dimensional Riemannian 
manifolds with holonomy group $\G_2$ \cite{Joyce}. 
\item The octonionic Hopf map and instanton solutions 
of the Yang-Mills equations in 8 dimensions \cite{GKS}. 
\item Octonionic aspects of 10-dimensional superstring theory and 
10-dimensional super-Yang-Mills theory 
\cite{CH,Deligne,Evans,KT,Schray,Sierra}.
\item Octonionic aspects of 11-dimensional supergravity and supermembrane 
theories, and the role of Joyce manifolds in compactifying 11-dimensional 
supergravity to obtain theories of physics in 4 dimensions \cite{Duff}. 
\item Geoffrey Dixon's extension of the Standard Model based on the
algebra $\C \tensor \H \tensor \O$ \cite{Dixon}.  
\item Other attempts to use the octonions in physics 
\cite{CMT,GT,LPS,Okubo}.
\end{itemize}  

\noindent I urge the reader to explore these with the help of the references.

\subsection*{Acknowledgements}  

I thank John Barrett, Toby Bartels, Robert Bryant, Geoffrey Dixon, James 
Dolan, Tevian Dray, Bertram Kostant, Linus Kramer, Pertti Lounesto,
Corinne Manogue, John McKay, David Rusin, Tony Smith, Anthony Sudbery, and 
Matthew Wiener for useful discussions.

\end{document}